\documentclass[10pt,a4paper]{article}

\usepackage[a4paper,margin=2.5cm]{geometry}

\usepackage{amsmath,amssymb,amsthm,mathtools}
\usepackage{mathabx,stmaryrd,bm}
\input{Ustmry.fd}
\DeclareFontShape{U}{stmry}{b}{n}{<->ssub*stmry/m/n}{}

\usepackage{array,booktabs,enumitem}
\usepackage{graphicx,calc,xcolor,tikz}
\usetikzlibrary{arrows.meta}

\usepackage{xurl}
\usepackage[breaklinks,bookmarks=false]{hyperref}
\hypersetup{colorlinks=true,linkcolor=blue,citecolor=blue,urlcolor=blue,
  plainpages=false,pdfwindowui=false,pdfstartview={FitH}}

\newtheorem{lemma}{Lemma}[section]
\newtheorem{theorem}[lemma]{Theorem}
\newtheorem{proposition}[lemma]{Proposition}
\newtheorem{corollary}[lemma]{Corollary}
\newtheorem{assumption}[lemma]{Assumption}
\newtheorem{remark}[lemma]{Remark}
\newtheorem{example}[lemma]{Example}

\newcommand{\email}[1]{\href{mailto:#1}{#1}}
\newcommand{\HHOabstractstart}{\noindent}
\newcommand{\HHOmetadata}{
  \bigskip\noindent\textbf{Keywords: }\KW\par
  \noindent\textbf{MSC2020:} 35J60, 65N12, 65N15, 65N30.}

\numberwithin{equation}{section}

\def\bae#1\eae{\begin{equation}\begin{aligned}#1\end{aligned}\end{equation}}
\def\ba#1\ea{\begin{align}#1\end{align}}
\def\ban#1\ean{\begin{align*}#1\end{align*}}
\newcommand{\bse}{\begin{subequations}}
\newcommand{\ese}{\end{subequations}}
\newcommand{\bt}{\begin{theorem}}
\newcommand{\et}{\end{theorem}}
\newcommand{\bpr}{\begin{proposition}}
\newcommand{\epr}{\end{proposition}}
\newcommand{\bl}{\begin{lemma}}
\newcommand{\el}{\end{lemma}}
\newcommand{\bc}{\begin{corollary}}
\newcommand{\ec}{\end{corollary}}
\newcommand{\bp}{\begin{proof}}
\newcommand{\ep}{\end{proof}}
\newcommand{\br}{\begin{remark}}
\newcommand{\er}{\end{remark}}
\newcommand{\bas}{\begin{assumption}}
\newcommand{\eas}{\end{assumption}}

\newcommand\Om{\Omega}

\newcommand\eq{\coloneqq}

\newcommand\tL{\bm{L}}
\newcommand\tW{\bm{W}}

\newcommand\F{\mathcal{F}}

\newcommand\T{\mathcal{T}}

\newcommand\RR{{\mathbb R}}
\newcommand\NN{{\mathbb N}}

\newcommand\PP{{\mathbb P}}

  \newcommand{\restri}[1]{\ \!\!_{|_{#1}}}

\DeclareMathOperator*{\essinf}{ess\,inf}
\DeclareMathOperator*{\esssup}{ess\,sup}

\def\b{\boldsymbol}
\newcommand*\cst[1]{\mathrm{#1}}

\newcommand\dU[2]{{\u U}_{#1}^{#2}}

\renewcommand\u{\underline}
\newcommand{\Fb}{\F_h^\cst{b}}

\newcommand{\res}[1]{\ \!\!_{|_{#1}}}
\newcommand{\dgrad}[2]{\b{\cst{G}}^{#1}_{#2}}
\newcommand{\drec}[2]{\cst{R}^{#1}_{#2}}
\newcommand{\dfbres}[2]{\Delta^{#1}_{#2}}
\newcommand{\GRAD}{\b\nabla}
\newcommand{\DIV}{\b\nabla{\cdot}}
\newcommand\Interp[2]{\u{I}_{#1}^{#2}}
\newcommand\proj[2]{\pi_{#1}^{#2}}
\newcommand\PROJ[2]{\b{\pi}_{#1}^{#2}}
\newcommand\stress{\b\sigma}

\newcommand\stabfl{\sigma}

\newcommand\sob{p}
\newcommand\sobmax{\widehat{\sob}}
\newcommand\sobmin{\widecheck{\sob}}
\newcommand\petp{{\scalebox{0.5}{$+$}}}
\newcommand\petm{{\scalebox{0.5}{$-$}}}

\newcommand\Plog{\mathcal P^\cst{log}}
\newcommand\Plogm{\mathcal P^\cst{log}_-}
\newcommand\Plogp{\mathcal P^\cst{log}_+}
\newcommand\Plogpm{\mathcal P^\cst{log}_{\pm}}

\title{A priori error estimates for a Hybrid High-Order method on variable-exponent Leray--Lions problems}
\author{Andr\'e Harnist\footnotemark[1]$~^{1}$\and Romain Mottier\footnotemark[2]$~^{2}$}
\def\KW{variable-exponent Leray--Lions equations, Hybrid High-Order methods, a priori error estimates, nonlinear stabilization, polytopal meshes}

\begin{document}
\maketitle

\footnotetext[1]{\email{andre.harnist@utc.fr}, \url{https://andreharnist.fr}, corresponding author}
\footnotetext[2]{\email{romain.mottier@outlook.com}, \url{https://romainmottier.github.io}}
\renewcommand{\thefootnote}{\fnsymbol{footnote}}
\footnotetext[1]{Université de technologie de Compiègne, Alliance Sorbonne Université, LMAC, Compiègne, France}
\footnotetext[2]{University of Basel, Department of Mathematics and Computer Science, Basel, Switzerland}
\renewcommand{\thefootnote}{\arabic{footnote}}

\begin{abstract}
\HHOabstractstart
{
We derive a priori error estimates for the Hybrid High-Order (HHO) discretization of variable-exponent Leray--Lions problems introduced in a companion convergence study~\cite{harnist-romain-conv}. HHO methods belong to the family of hybrid nonconforming methods, based on cell and face polynomial unknowns of arbitrary degree. Under additional broken regularity of the exact solution and its flux, we establish a localized error estimate that retains the individual element sizes and a global estimate whose rate distinguishes the purely subquadratic, mixed subquadratic--superquadratic, and purely superquadratic regimes. In the subquadratic range, we also derive a regime-dependent refinement that compares the local approximation scale with the local non-degeneracy scale and recovers the optimal polynomial order when that scale is resolved by the mesh. The analysis relies on variable-exponent approximation properties of polynomial projectors, mean-exponent tools for discrete polynomial spaces, and weighted monotonicity estimates. The framework also covers regularized and mildly non-symmetric anisotropic Carreau--Yasuda-type fluxes. Numerical experiments on Cartesian and Voronoi meshes support the predicted rates and show that the regularity assumption on the flux cannot be dispensed with.
}
\end{abstract}

\HHOmetadata

\section{Introduction}
\label{intro}

We study the homogeneous Dirichlet boundary-value problem
\bse\label{eq:intro.problem}
\begin{alignat}{2}
-\operatorname{div}\bigl(\stress(\cdot,\nabla u)\bigr)&=f &&\quad\text{in }\Omega,\\
u&=0 &&\quad\text{on }\partial\Omega,
\end{alignat}
\ese
where $\Omega\subset\mathbb R^d$ is a bounded polytopal Lipschitz domain and $\stress:\Omega\times\mathbb R^d\to\mathbb R^d$ is a Leray--Lions flux with spatially varying growth $p(\cdot)$. The $p(\cdot)$-Laplacian is the basic prototype, and the class considered below also includes regularized and anisotropic Carreau--Yasuda-type laws with spatially dependent parameters. Such variable-growth models arise, in particular, in electrorheological and generalized Newtonian fluids and in variational models with nonstandard growth; see, e.g., \cite{Ruzicka2000,diening2011lebesgue}. In the companion paper~\cite{harnist-romain-conv}, we introduced a Hybrid High-Order (HHO) discretization of~\eqref{eq:intro.problem} and proved its well-posedness and its convergence to the weak solution under minimal regularity. The purpose of the present paper is to develop a quantitative a priori error analysis of the same method, under additional regularity of the solution and of its flux.

A priori error analysis for nonlinear elliptic problems with constant power growth has a long history. Early finite element analyses for monotone quasilinear equations and the $p$-Laplacian include \cite{Glowinski.Marrocco:75,Chow1989,Barrett.Liu:93}. The introduction of problem-adapted quasi-norms was decisive for obtaining sharp estimates in degenerate regimes \cite{LiuBarrett1994Degenerate,Liu.Yan:01,EbmeyerLiu2005}. Mixed formulations and nonconforming methods provide complementary approaches, with error estimates covering the full range $1<p<\infty$ \cite{FarhloulManouzi2000,Liu.Yan:01}, and adaptive finite element methods for the $p$-Laplacian were subsequently analyzed in \cite{belenki:12}. Discontinuous Galerkin (dG) methods for systems with $p$-structure were analyzed in \cite{DieningKronerRuzickaToulopoulos2014}, and recent quasi-optimal dG analyses have substantially weakened the regularity requirements on the exact solution \cite{BlechtaGazcaKaltenbachRuzicka2026}. These works show that the natural error measure and the convergence rate depend strongly on the degeneracy of the constitutive law, and on whether the error is measured in the energy norm, in a weighted quadratic distance, or in a quasi-norm.

Closely related estimates arise in the numerical analysis of non-Newtonian flows. For power-law and Carreau-type constitutive relations, conforming finite element error estimates were established in \cite{BarrettLiu1993Flow,BarrettLiu1994Flow}, and velocity and pressure estimates in suitable nonlinear metrics were developed for $p$-Stokes and generalized Stokes systems in \cite{BelenkiBerselliDieningRuzicka2012,Hirn:13}. The same monotonicity and continuity structure underlies numerical analyses of glaciological and generalized Newtonian models, such as \cite{Glowinski.Rappaz:03}. This literature is particularly relevant here because, in the shear-thinning range, the sharp consistency analysis exploits the weighted quadratic structure of the stress rather than a direct H\"older estimate of the Nemytskii map.

Polytopal methods have brought these ideas to general meshes. Examples include dG~\cite{CangianiDongGeorgoulisHouston2017} and virtual element (VEM)~\cite{Beirao-da-Veiga.Brezzi.ea:13} methods, as well as hybrid nonconforming methods, which attach polynomial unknowns both to the mesh cells and to the mesh faces, such as the hybridizable discontinuous Galerkin (HDG)~\cite{CockburnGopalakrishnanLazarov2009}, weak Galerkin (WG)~\cite{WangYe2013}, and Hybrid High-Order (HHO)~\cite{DiPietroErnLemaire2014,DiPietroErn2015} methods; we refer to the monographs~\cite{Di-Pietro.Droniou:20,CicuttinErnPignet2021} for HHO methods. In HHO methods, a cellwise gradient reconstruction is built from the polynomial cell and face unknowns, while a stabilization term weakly controls the discrepancy between cell traces and face unknowns. This construction supports arbitrary polynomial degrees and general polytopal cells and is locally conservative; the cell unknowns can moreover be eliminated by static condensation, leaving only face unknowns coupled in the global system. For constant-exponent Leray--Lions equations, the HHO method was introduced and proved convergent in \cite{Di-Pietro.Droniou:17}, and the $W^{s,p}$ approximation properties of the local projectors and the first HHO a priori estimates were obtained in \cite{Di-Pietro.Droniou:17*1}, with the robust rates $h^{(k+1)(p-1)}$ for $p<2$ and $h^{(k+1)/(p-1)}$ for $p>2$. For $1<p\le2$, the regime-dependent refinement of \cite{harnist2} shows how the rate improves up to the optimal order $k+1$ when the local approximation scale is small compared with the non-degeneracy scale. Other HHO approaches to degenerate convex minimization and to the $p$-Laplace problem include \cite{carstensen:20,RenLiu2024}, and HHO discretizations of non-Newtonian creeping and convective flows were analyzed in \cite{Botti.Castanon-Quiroz.ea:20,harnist3}.

Among related polytopal discretizations, mimetic schemes were developed for quasilinear elliptic equations \cite{Antonietti.Bigoni.ea:14}. For incompressible non-Newtonian Stokes flows, the Virtual Element analysis of \cite{AntoniettiBeiraoBottiVaccaVerani2024} establishes a velocity error rate of $h^{kr/2}$ in the broken $W^{1,r}$ norm for the degenerate shear-thinning regime $1<r\le2$, where $k$ is the polynomial degree of the VEM discretization. Its shear-thickening extension \cite{harnist5} proves, for $r>2$, velocity and pressure rates of $h^{k/(r-1)}$ for $\delta\ge0$ and the improved rate $h^{2k/r}$ for $\delta>0$ under additional regularity assumptions. These latter estimates are expressed in a discrete norm, which yields the standard $W^{1,r}$ velocity estimate on convex meshes. The exponent dependence of the robust rates parallels that of the present HHO estimates, up to the different polynomial-degree conventions and model problems; the non-degenerate VEM refinement further highlights the role of $\delta$. More recently, a mixed Virtual Element Method for the scalar $p$-Laplace problem with a constant exponent $1<p<\infty$ was analyzed in \cite{HaileVacca2026}, including well-posedness and a priori error estimates for both the primal variable and the flux. These works motivate keeping track of the local degeneracy regime when passing to spatially varying growth.

For spatially varying growth, the analytical setting is provided by variable-exponent Lebesgue and Sobolev spaces \cite{Ruzicka2000,diening2011lebesgue}, and the numerical literature is considerably more recent; we refer to \cite{harnist-romain-conv} for a broader review. A priori error estimates for conforming finite element approximations of the $p(\cdot)$-Laplacian were obtained in \cite{BreitDieningSchwarzacher2015,DelPezzoMartinez2015}: the former analyzes the error in a nonlinear quasi-norm, while the latter establishes convergence rates in the subquadratic range on two-dimensional convex domains. A medius error estimate and ensuing a priori convergence rates for a Crouzeix--Raviart discretization of the variable-exponent Dirichlet problem were derived in \cite{BalciKaltenbach2025}. Related quantitative velocity and pressure estimates for finite element discretizations of electrorheological-flow models were developed in \cite{BerselliBreitDiening2016}. Other contributions include an interior-penalty dG method with convergence of discrete minimizers \cite{DelPezzoLombardiMartinez2012}, nonconforming discretizations of non-autonomous variational problems \cite{BalciOrtnerStorn2022}, finite-volume and control-volume schemes \cite{AndreianovQuenjel2023}, and solution strategies based on decomposition--coordination and quasi-Newton minimization \cite{DelPezzoMartinez2014,CaliariZuccher2017,AragonBonderRubio2023}. To the best of our knowledge, however, quantitative a priori estimates for a hybrid nonconforming discretization, such as HDG, WG, or HHO, that keeps the pointwise exponent throughout both the bulk and face contributions have not previously been established.

The passage from constant to variable exponent creates several additional difficulties. Polynomial inverse, trace, and approximation estimates are naturally formulated with a fixed exponent on each cell, while the continuous energy and the constitutive inequalities involve the pointwise function $p(\cdot)$. The log-H\"older condition is therefore used to compare local polynomial norms at the mean exponent with variable-exponent modulars, uniformly in the mesh size. Moreover, the strong-monotonicity and consistency arguments must distinguish the subquadratic and superquadratic regions without replacing the pointwise powers by elementwise mean exponents. The bulk degeneracy $\delta$ and the face degeneracies $(\delta_F)_{F\in\F_h}$ are retained throughout, and enter the estimates through the same discrete data quantity as in \cite{harnist-romain-conv}.

The contributions of this paper are as follows. First, we establish the variable-exponent approximation theory needed by the HHO analysis: the bounded-projector approximation lemma, the uniform $L^{p(\cdot)}$-stability of local $L^2$-orthogonal polynomial projectors, and the resulting optimal cell and trace estimates. A geometric ingredient of this analysis is a uniform John-domain construction for patches made of a bounded number of shape-regular simplices, such as the mesh cells, which allows us to invoke the variable-exponent Poincar\'e--Wirtinger inequality without reverting to star-shaped cells. These approximation results rely on the polynomial mean-exponent equivalence established independently of these approximation estimates in the companion convergence paper~\cite{harnist-romain-conv}. The corresponding approximation statements are also recalled there as tools for the convergence analysis. As a complementary discrete tool, we establish an HHO discrete Lebesgue embedding in Subsection~\ref{subsec:discrete.lebesgue}, transferring polynomial Lebesgue embeddings to local HHO seminorms. 
This organization reflects the complementary objectives of the two papers: the companion work establishes well-posedness and qualitative convergence under minimal regularity assumptions, while the present paper develops the approximation theory and derives quantitative a priori error estimates under additional regularity assumptions.

Second, preserving the three-term consistency decomposition of the HHO error equation, we derive a localized a priori estimate in the natural discrete $W^{1,p(\cdot)}$-type norm, in which the approximation contributions retain the element sizes $h_T$ rather than being globalized through $h=\max_T h_T$. Under uniform control of the regularity seminorms displayed in the estimate, its global counterpart gives
\bae
 \|\u u_h-\Interp{h}{k}u\|_{p(\cdot),h}
 =O\!\left(h^{(k+1)\alpha}\right),
 \qquad
 \alpha=
 \begin{cases}
 \dfrac{p_\petm}{2}, & p_\petp\le2,\\[2mm]
 \min\!\left(\dfrac1{p_\petp-1},\dfrac{p_\petm}{p_\petp}\right),
 & p_\petm<2<p_\petp,\\[3mm]
 \dfrac1{p_\petp-1}, & p_\petm\ge2.
 \end{cases}
\eae
The estimate thus keeps a genuinely variable-exponent dependence through $p_\petm$ and $p_\petp$, including in the mixed subquadratic--superquadratic case. For a constant exponent, it gives the rates $(k+1)p/2$ for $1<p\le2$ and $(k+1)/(p-1)$ for $p\ge2$. In the subquadratic range, the robust exponent $p/2$ improves on the exponent $p-1$ of the earlier HHO estimate \cite{Di-Pietro.Droniou:17*1}, under the regularity hypotheses assumed here. The factor $p/2$ is also consistent with the $r/2$ factor in the VEM velocity estimate of \cite{AntoniettiBeiraoBottiVaccaVerani2024} for pseudoplastic Stokes flows, for a different model and error quantity. This comparison does not supersede the regime-dependent estimates of \cite{harnist2}, which can attain order $k+1$ under suitable non-degeneracy conditions; in the superquadratic range, the present exponent agrees with that of \cite{Di-Pietro.Droniou:17*1}. The cited variable-exponent finite element analyses use different error metrics and regularity assumptions, whereas our results cover arbitrary fixed polynomial degrees on polytopal meshes and Leray--Lions fluxes satisfying the structural hypotheses, with rates expressed explicitly in terms of $p_\petm$ and $p_\petp$ in each regime. The present analysis is also complementary to the interior-penalty dG analysis of \cite{DelPezzoLombardiMartinez2012}, whose emphasis is on the convergence of discrete minimizers rather than on quantitative error bounds.

Third, in the subquadratic regime $p_\petp\le2$, we derive a regime-dependent refinement governed by a local indicator $\eta_T$ that compares the approximation scale with the bulk and face non-degeneracy scales. The proof uses the HHO discrete Lebesgue embedding of Subsection~\ref{subsec:discrete.lebesgue} (Lemma~\ref{lem:error:local-change-exp}), based on the polynomial embedding in \cite{harnist-romain-conv} to change exponent only on polynomial discrete errors. When the non-degeneracy scale is resolved uniformly by the mesh, the refined estimate recovers the optimal rate
\bae
 \|\u u_h-\Interp{h}{k}u\|_{p(\cdot),h}
 =O(h^{k+1});
\eae
when it is not, the robust estimate above remains available.

Fourth, the analysis covers a mildly non-symmetric anisotropic Carreau--Yasuda law, for which we verify in the appendix that the structural assumptions of the analysis hold.

Numerical experiments on Cartesian and Voronoi meshes confirm these results. For smooth solutions, the observed convergence rates reach at least the orders predicted by Remark~\ref{rem:error:rates-general} in the subquadratic, mixed, and superquadratic regimes. The numerical experiments also illustrate the role of the flux-regularity assumption: when the flux loses regularity at a critical point of the solution where $p\neq2$, increasing the polynomial degree eventually brings no gain, and the rate falls below the predicted order. Finally, in the non-degenerate case $\delta=1$, the optimal order $k+1$ is observed in all three regimes, which, for $p_\petp\le2$, agrees with the regime-dependent refinement.

The paper is organized as follows. Section~\ref{sec:continuous.setting} recalls the continuous model and the structural assumptions. Section~\ref{sec:scheme} details the HHO discrete setting, the reconstructions, the stabilization, and the discrete problem. Section~\ref{sec:discrete.norms} develops the discrete seminorm equivalences, the HHO discrete Lebesgue embedding, and a local trace inequality, while Section~\ref{sec:properties} establishes the approximation and stability properties of the discrete operators. Section~\ref{sec:error.estimate} states the localized and global a priori error bounds and the regime-dependent refinement, together with their proofs. Section~\ref{sec:numerics-a-priori} presents numerical experiments, and Section~\ref{sec:conclusion} concludes. The appendices summarize the variable-exponent background and polynomial tools and verify the anisotropic Carreau--Yasuda example.

\section{Continuous problem}\label{sec:continuous.setting}
\label{sec:setting}

{Throughout, $\Omega \subset \RR^d$, $d\in\NN^*$, is a bounded, connected, polytopal domain with Lipschitz boundary $\partial\Omega$, and the model is given by \eqref{eq:intro.problem}.} {Here, $f : \Omega \to \RR$ denotes the volume source, whereas} $\stress : \Omega \times \RR^d \to \RR^d$ is the flux function; it acts on the gradient of the solution and may depend on the position.

{
\br[Nonhomogeneous Dirichlet data and manufactured solutions]\label{rem:nonhomogeneous}
{The theoretical results below are formulated for zero Dirichlet data, which simplifies both the continuous and discrete spaces.} {The manufactured tests also include solutions whose boundary trace $g$ is nonzero.} {If $g$ has a lifting} $G\in W^{1,p(\cdot)}(\Omega)$ with ${\rm tr}\,G=g$, the associated continuous weak problem is posed in the affine space $G+W^{1,p(\cdot)}_0(\Omega)$ and tested against $W^{1,p(\cdot)}_0(\Omega)$. {Equivalently, setting $u=G+v$ leads to a nonlinear equation for the zero-trace unknown $v$;} the flux then depends on $\nabla G+\nabla v$. {In the discrete formulation, boundary face degrees of freedom are set from a suitable approximation of $g$.} {Structural monotonicity should still control differences of admissible solutions,} but a full extension of the convergence and, in particular, of the quantitative error estimates requires control of the chosen discrete lifting and its additional consistency terms. {The results stated below do not include this extension.}
\er
}

\subsection{Flux function}\label{sec:strain.rate.shear.stress.law}

{We impose the following structural conditions on $\stress$; the relevant variable-exponent function spaces are recalled in Appendix~\ref{app:ve-background}.}

\bas[Flux function]\label{ass:stress}
  {The map $\stress : \Omega \times \RR^d \to \RR^d$ is measurable and there are} $\sob \in \Plogpm(\Omega)$, $\delta \in L^{p(\cdot)}(\Om,[0,\infty))$ and $\sigma_\cst{c},\sigma_\cst{m} \in (0,\infty)$ such that, for all $\b\tau,\b\eta \in \RR^d$ with $\b \tau \ne \b \eta$, and a.e. $\b x \in \Omega$, {it satisfies the H\"older continuity, strong monotonicity, and zero-at-rest conditions below, respectively:}
 \bse\label{eq:ass:stress}
 \ba
       |\stress(\b x,\b\tau)-\stress(\b x,\b\eta)| &\le \sigma_\cst{c} \left(\delta(\b x)^{\sob(\b x)}+|\b\tau|^{\sob(\b x)}+|\b\eta|^{\sob(\b x)}\right)^\frac{{\sob(\b x)}-2}{{\sob(\b x)}}| \b\tau-\b\eta |,
    \label{eq:ass:stress:cont} 
\\
      \left(\stress(\b x,\b\tau)-\stress(\b x,\b\eta)\right)\cdot\left(\b\tau-\b\eta\right) &\ge \sigma_\cst{m}\left(\delta(\b x)^{\sob(\b x)}+|\b\tau|^{\sob(\b x)}+|\b\eta|^{\sob(\b x)}\right)^\frac{{\sob(\b x)}-2}{{\sob(\b x)}}|\b\tau-\b\eta|^{2},\label{eq:ass:stress:mono}\\
      \stress(\b x,\b 0) &= \b 0.\label{eq:ass:stress:zero}
\ea
\ese
\eas
{In the above inequalities and throughout the paper, a non-negative quantity raised to the power zero is interpreted as one, including at a zero base; expressions involving coincident arguments in the structural estimates are understood by the corresponding limiting value.}

{
{To limit the number of named constants in this section, we use} $a\lesssim b$ (respectively, $a\gtrsim b$) for non-negative quantities $a,b$ whenever $a\le Cb$ (respectively, $a\ge Cb$), where $C\in(0,\infty)$ depends only on $p_\petm$, $p_\petp$, $p_\cst{log}$, and $\Omega$ (including the fixed spatial dimension $d$). {In addition, $a\eqsim b$ means} $C^{-1}b\le a\le Cb$, equivalently $a\lesssim b\lesssim a$. {These implicit factors are, in particular, independent of} $\sigma_\cst{m}$ or $\sigma_\cst{c}$, which are kept explicit when relevant. {Any further dependencies are indicated where relevant.}
}

Some remarks are in order.

\br[Residual shear stress]
  {The condition \eqref{eq:ass:stress:zero} may be weakened to} $\stress(\cdot,\b 0)  \in L^{{\sob}'(\cdot)}(\Omega,\RR^d)$. {The corresponding adjustments to the analysis are straightforward and are omitted.}
\er

\br[Singular exponents]
  Inequalities \eqref{eq:ass:stress:cont}--\eqref{eq:ass:stress:mono} {imply the following estimates in a form suitable for subsequent H\"older arguments.}
  {For a.e. $\b x \in \Omega$ and any distinct $\b\tau,\b\eta \in \RR^d$,}
  \bse\label{eq:ass:stress:sing}
  \ba
      \left|
        \stress(\b x,\b\tau)-\stress(\b x,\b\eta)
        \right| &\lesssim\sigma_\cst{c} \left(\delta(\b x)^{\sob(\b x)}+|\b\tau|^{\sob(\b x)}+|\b\eta|^{\sob(\b x)}\right)^{\frac{{\sob(\b x)}-{\sobmin(\b x)}}{{\sob(\b x)}}}| \b\tau-\b\eta |^{{\sobmin(\b x)}-1},\label{eq:ass:stress:sing:cont}
      \\
        \left(\stress(\b x,\b\tau)-\stress(\b x,\b\eta)\right)\cdot(\b\tau-\b\eta)  &\gtrsim\sigma_\cst{m}\left(\delta(\b x)^{\sob(\b x)}+|\b\tau|^{\sob(\b x)}+|\b\eta|^{\sob(\b x)}\right)^{\frac{{\sob(\b x)-\sobmax(\b x)}}{{\sob(\b x)}}}|\b\tau-\b\eta|^{{\sobmax(\b x)}}.\label{eq:ass:stress:sing:mono}
\ea
\ese
{The notation $\sobmin$ and $\sobmax$ is introduced in \eqref{eq:exp:conj.lower.uper}.}
{These bounds follow from \eqref{eq:ass:stress:cont}--\eqref{eq:ass:stress:mono} and the triangle inequality $|\b\tau-\b\eta|^p\le 2^{p-1}(|\b\tau|^p+|\b\eta|^p)$, applied separately for $p<2$ and $p\ge2$. The hidden constants depend only on $p_\petm$ and $p_\petp$.}
\er

\br[Order between the continuity and monotonicity constants]
  {The relation $\sigma_\cst{m} \le \sigma_\cst{c}$ follows on choosing} $\b\tau,\b\eta \in \RR^d$ such that $|\b\tau| = 1$ and $\b \eta = \b 0$ in \eqref{eq:ass:stress} and writing a.e. in $\Om$,
  \bae
 \label{eq:power-framed:constants.bound}
    \sigma_\cst{m} \le  \frac{\stress(\cdot,\b\tau)\cdot\b\tau}{\left(\delta^{\sob}+|\b\tau|^{\sob}\right)^\frac{{\sob}-2}{{\sob}}|\b\tau|^{2}}\le \frac{|\stress(\cdot,\b\tau)||\b\tau|}{\left(\delta^{\sob}+|\b\tau|^{\sob}\right)^\frac{{\sob}-2}{{\sob}}|\b\tau|^{2}} \le \sigma_\cst{c}.
\eae
\er

{We illustrate Assumption~\ref{ass:stress} with three admissible constitutive fluxes.}

\begin{example}[$p(\cdot)$-Laplace flux function]\label{ex:laplace}
{For the $p(\cdot)$-Laplacian} $-\DIV(|\GRAD u|^{p(\cdot)-2}\GRAD u)=f$, {the constitutive flux reads} $\stress(\b x,\b\tau)=|\b\tau|^{p(\b x)-2}\b\tau$ and $\delta(\b x)=0$. {The classical vector-valued power inequalities, applied pointwise} with $p=p(\b x)$ and using $1<p_-\le p(\b x)\le p_+<\infty$, {show that} Assumption~\ref{ass:stress} holds with positive constants depending only on $p_-$ and $p_+$; see, e.g., \cite[Appendix~A]{Botti.Castanon-Quiroz.ea:20} for the corresponding power-framed inequalities.
\end{example}

\begin{example}[Carreau--Yasuda flux function]\label{ex:Carreau--Yasuda}
  {A further flux satisfying Assumption~\ref{ass:stress}, motivated by Carreau--Yasuda rheology, is defined for} almost every $\b x\in\Omega$ and all $\b\tau \in \RR^d$,
  \bae\label{eq:Carreau--Yasuda}
    \stress(\b x,\b\tau) = \mu(\b x)\left(\delta(\b x)^{\alpha(\b x)}+|\b\tau|^{\alpha(\b x)}\right)^\frac{\sob(\b x)-2}{\alpha(\b x)}\b\tau,
  \eae
  {Here} $\mu : \Omega \to [\mu_\petm,\mu_\petp] \subset (0,+\infty)$ {is the measurable \emph{local flow consistency index}}, $\delta \in L^{p(\cdot)}(\Omega,[0,+\infty))$ is the \emph{degeneracy parameter}, $\alpha : \Omega \to [\alpha_\petm,\alpha_\petp] \subset (0,+\infty)$ {is the measurable \emph{local transition flow behavior index}}, and $\sob \in \Plogpm(\Omega)$ is the \emph{flow behavior index}.
  {The power-framed inequalities of \cite[Appendix~A]{Botti.Castanon-Quiroz.ea:20}, used pointwise} in $\b x$ and uniformly over the prescribed ranges of $p(\b x)$, $\alpha(\b x)$, and $\mu(\b x)$, {verify Assumption~\ref{ass:stress} for $\stress$} with positive continuity and monotonicity constants depending only on these ranges.
{For $\delta \equiv 0$ and $\mu \equiv 1$, this reduces to} the $p(\cdot)$-Laplace example \ref{ex:laplace}.
\end{example}

\begin{example}[Non-symmetric anisotropic Carreau--Yasuda flux]\label{ex:Carreau--Yasuda;mat}
For a.e. $\b x\in\Omega$ and every $\b\tau\in\mathbb R^d$, set
\bae\label{eq:Carreau--Yasuda:mat}
\stress(\b x,\b\tau)=
\left(\delta(\b x)^{\alpha(\b x)}+|\b A(\b x)\b\tau|^{\alpha(\b x)}\right)^{\frac{p(\b x)-2}{\alpha(\b x)}}
\b B(\b x)\b\tau.
\eae
Here $\alpha:\Omega\to[\alpha_\petm,\alpha_\petp]\subset(0,\infty)$ is measurable, $\b A$ is measurable and uniformly invertible,
\bae
 A_\petm|\b\xi|^2\le |\b A(\b x)\b\xi|^2\le A_\petp|\b\xi|^2,
\eae
and
\bae\label{eq:anisotropic:compatibility}
\b B(\b x)=\b A(\b x)^\text{T}(\mu(\b x)\b I_d+\b K(\b x))\b A(\b x),
\qquad
\|\b K(\b x)\|_{\rm op}\le\varepsilon_0\mu(\b x),
\eae
with $0<\mu_\petm\le\mu\le\mu_\petp$ and $\varepsilon_0$ sufficiently small in terms of $p_\petm,p_\petp$. Appendix~\ref{sec:properties.stress} verifies Assumption~\ref{ass:stress}. The case $\b K\ne0$ allows a genuinely non-symmetric constitutive law.
\end{example}

\subsection{Weak formulation}\label{sec:weak.formulation}

{The homogeneous Dirichlet condition is incorporated into the potential space} 
\bae
U \eq W_0^{1,\sob(\cdot)}(\Omega) \eq \{v \in W^{1,\sob(\cdot)}(\Omega) : v\restri{\partial\Omega} = 0\}.
\eae
{For $f \in L^{\sob'(\cdot)}(\Omega)$, the variational form of \eqref{eq:intro.problem} is:}
Find $u \in U$ such that
\bae\label{eq:lr.weak}
     \cst{a}(u,v) = (f, v)_\Om \qquad \forall v \in U.
\eae
{The \emph{diffusion function} $\cst{a} : U \times U \to \RR$ is then specified, for all $v,w \in U$, by}
\bae\label{eq:a.b}
  \cst{a}(w,v) \eq ( \stress(\cdot,\GRAD w) ,\GRAD v)_\Om.
\eae

The well-posedness of problem \eqref{eq:lr.weak} and the following a priori estimate are established in \cite{harnist-romain-conv}.

\bt[{ Continuous well-posedness and variable-exponent a priori estimate}]\label{prop:a-priori}
{Assumption~\ref{ass:stress} ensures that the continuous formulation \eqref{eq:lr.weak} has exactly one solution} $u\in U$ satisfying
\bae\label{eq:continuous.solution:bounds:uh}
\rho_{\tL^{\sob(\cdot)}(\Om)}(\GRAD u)
\lesssim\max(1,\sigma_\cst{m}^{-p'_\petm},\sigma_\cst{m}^{-p'_\petp})\max\!\left(
\rho_{L^{\sob'(\cdot)}(\Om)}(f)^{p'_\petp/p'_\petm},
\rho_{L^{\sob'(\cdot)}(\Om)}(f)^{p'_\petm/p'_\petp},
\rho_{L^{\sob(\cdot)}(\Om)}(\delta)
\right).
\eae
{The hidden constant depends only on $p_\petm$, $p_\petp$, $p_\cst{log}$, and $\Omega$.}
\et

\section{HHO discretization}\label{sec:scheme}

We specify the discrete unknowns and the local reconstruction and stabilization operators before stating the HHO problem.

\subsection{Discrete setting}\label{sec:discrete.setting}

\paragraph{Mesh.}

{For a subset $U\subset\RR^d$, its diameter is written $h_U$.}
{We represent a polytopal mesh by the pair} $\mathcal M_h\coloneq(\T_h,\F_h)$, where $\T_h$ is a finite collection of polytopal elements with mesh size $h\coloneq\max_{T\in\T_h}h_T$, while $\F_h$ is a finite collection of hyperplanar faces.
{Throughout, $\mathcal M_h$ is required to satisfy the geometric conditions of} \cite[Definition 1.7]{Di-Pietro.Droniou:20}, {and $\kappa$ denotes its regularity parameter.}
{Faces belonging to $\partial\Omega$ and interior faces in}
$\Omega$ are collected in the sets $\F_h^{\rm b}$ and $\F_h^{\rm i}$, respectively.
{For each $T\in\T_h$, $\F_T$ comprises those faces in $\F_h$ that form $\partial T$.} {Likewise, $\T_F$ is the set of mesh cells adjacent to a face $F \in \F_h$, comprising} the one (if $F\in\F_h^{\rm b}$) or two (if $F\in\F_h^{\rm i}$) elements on whose boundary $F$ lies.
{Given a cell $T\in\T_h$ and $F\in\F_T$,} $\b n_{TF}$ {is the constant outward unit normal to $T$ on $F$}, and we denote $(\b n_{\partial T})\res{F} \eq \b n_{TF}$.

{For the $h$-analysis we work with a refining mesh sequence, regular} in the sense of \cite[Definition 1.9]{Di-Pietro.Droniou:20}, with regularity parameter uniformly bounded away from zero. We write $\mathfrak H\subset(0,\infty)$ for the set of mesh sizes of this sequence, whose only accumulation point is $0$.
These mesh assumptions ensure uniform comparability of cell diameters and characteristic face scales, and a mesh-independent upper bound on the number of faces per cell.

\br[One-dimensional meshes]\label{rm:one-dimensional-mesh}
In one space dimension, the cells are intervals and their faces consist of single endpoints. We integrate on faces using $\mathcal H^0$ (counting measure), so that $\PP^k(F)=\RR$ and integrals over $\partial T$ reduce to two-point sums. Since the geometric diameter of $F$ vanishes, any face length $h_F$ used in the analysis is interpreted as the characteristic scale
\bae
 h_F\coloneq\frac{1}{|\T_F|}\sum_{T\in\T_F}h_T.
\eae
The shape-regularity assumptions give $h_F\eqsim h_T$ for every incident pair $F\in\F_T$ (cf.~\cite[Remark~1.13]{Di-Pietro.Droniou:20}). All reconstruction and stabilization formulas are already scaled with $h_T$, so the one-dimensional case needs no separate formulation.
\er

\paragraph{Projectors and broken spaces.}

{For $U \in \T_h \cup \F_h$ and $\ell \in \NN$, the notation $\PP^\ell(U)$ refers to restrictions to $U$ of} scalar-valued, $d$-variate polynomials of total degree $\le \ell$.
{The elementwise $L^2$-orthogonal projection} $\proj{U}{\ell} : L^{1}(U) \to \PP^\ell(U)$ {is characterized, for every $v \in L^{1}(U)$, by}
\bae\label{eq:proj}
  (\proj{U}{\ell}v, w)_U = (v,w)_U \qquad \forall w \in  \PP^{\ell}(U).
\eae
{On vector fields in $L^1(U)^d$,} the $L^2$-orthogonal projector mapping on $\PP^\ell(U)^d$ {operates component by component and is written in bold as} $\PROJ{U}{\ell}$. {We also employ the symbols} $\proj{h}{\ell}$ and $\proj{\partial T}{\ell}$ such that for all $T \in \T_h$ and all $F \in \F_T$, $(\proj{h}{\ell})\res{T} = \proj{T}{\ell}$ and $(\proj{\partial T}{\ell})\res{F} = \proj{F}{\ell}$.
The elliptic projector $\proj{T}{1,k+1}:W^{1,1}(T)\to\PP^{k+1}(T)$ is characterized by $(\GRAD(\proj{T}{1,k+1}v-v),\GRAD w)_T=0$ for all $w\in\PP^{k+1}(T)$, together with $(\proj{T}{1,k+1}v-v,1)_T=0$.

{To state the approximation estimates, we use the broken Sobolev space}
\bae
W^{n,\sob(\cdot)}(\T_h)\coloneq\left\{ v\in L^{\sob(\cdot)}(\Omega)\ : \ v\res{T}\in W^{n,\sob(\cdot)}(T)\quad\forall T\in\T_h\right\}.
\eae
{Its seminorm follows \eqref{eq:luxembourg}, with the associated modular given by} 
\bae 
\rho_{\dot W^{n,\sob(\cdot)}(\T_h)}(v) &\eq \sum_{T\in\T_h}\rho_{\dot W^{n,\sob(\cdot)}(T)}(v\res{T}).
\eae

For a piecewise Sobolev function $v$, the broken gradient $\nabla_h v$ is defined cellwise by $(\nabla_h v)|_T=\GRAD(v|_T)$ for each $T\in\T_h$.

\subsection{Hybrid space and norms}

{Fix $k\ge 0$ and choose $\widetilde k \in \{k,k+1\}$. } 
{Taking $\widetilde k=k$ leads to the equal-order case, whereas $\widetilde k=k+1$ gives the mixed-order case.} {Both choices are treated simultaneously by the notation adopted below.} We use $\PP^{\widetilde k}(\T_h)$ for elementwise polynomials and $\PP^k(\F_T):=\prod_{F\in\F_T}\PP^k(F)$ for the family of face polynomials on $\partial T$.  
We define the HHO space by
\bae
\dU{h}{k} \eq \left\{
\u v_h = ((v_T)_{T \in \T_h},(v_F)_{F\in \F_h}) \ : \ v_T \in \PP^{\widetilde k}(T)\ \ \forall T \in \T_h\ \mbox{ and }\ v_F \in \PP^k(F)\ \ \forall F \in \F_h \right\}.
\eae
{We impose the homogeneous boundary values by restricting the discrete potentials to the subspace of $\dU{h}{k}$:}
\bae
\dU{h,0}{k} \eq \left\{ \u v_h = ((v_T)_{T \in \T_h},(v_F)_{F\in \F_h}) \in \dU{h}{k} \ : \ v_F = 0 \quad \forall F \in \Fb \right\}.
\eae
{The interpolator $\Interp{h}{k} : W^{1,1}(\Omega) \to  \dU{h}{k} $ assigns to} $v \in W^{1,1}(\Omega)$ {the discrete vector $\Interp{h}{k}v$ specified below:}
\bae\label{eq:I:def}
  \Interp{h}{k} v \eq ((\proj{T}{\widetilde k} v\res{T})_{T \in \T_h},(\proj{F}{k} v\res{F})_{F \in \F_h}).
\eae
{On each $T \in \T_h$, $\dU{T}{k}$ and $\Interp{T}{k}$ stand for the local restrictions of} $\dU{h}{k}$ and $\Interp{h}{k}$ to $T$, respectively, and, for all $\u v_h \in \dU{h}{k}$, we let $\u v_T \eq (v_T,(v_F)_{F\in \F_T}) \in \dU{T}{k}$ {collect the degrees of freedom associated with $T$ and its faces.}
{For $\u v_h \in \dU{h}{k}$, we also assemble the broken polynomial field} $v_h\in\PP^{\widetilde k}(\T_h)$ {by joining its cell components, namely,}
\bae
(v_h)\res{T} \eq v_T\qquad\forall T \in \T_h.
\eae
{For each} $\u v_T \in \dU{T}{k}$, we define the broken polynomial field $v_{\partial T} \in \PP^{k}(\F_T)$ obtained by patching face unknowns, that is, 
\bae
(v_{\partial T})\res{F} \eq v_F\qquad\forall F \in \F_T.
\eae

\subsection{Operator reconstructions}

We construct the discrete operators locally from the cell and face unknowns: a reconstructed gradient and a higher-degree potential, complemented by a stabilization that controls discrepancies between cell traces and face unknowns. These local contributions are then assembled over the mesh.

\paragraph{Local gradient reconstruction.}

For all $T \in \T_h$, we define 
$\dgrad{k}{T} : \dU{T}{k} \to \PP^{k}(T)^d$ such that, for all $\u v_T \in \dU{T}{k}$,
\bae\label{eq:G}
  (\dgrad{k}{T} \u v_T,\b\tau)_{T} = (\GRAD v_T,\b\tau)_{T} + (v_{\partial T}-v_T,\b\tau \cdot \b n_{\partial T})_{\partial T}\qquad \forall \b\tau \in  \PP^{k}(T)^d.
\eae
The defining relation~\eqref{eq:G} reproduces integration by parts at the discrete level, with the face unknowns replacing the trace on $\partial T$. Taking the reconstructed gradient in the whole space $\PP^k(T)^d$ yields the commuting property~\eqref{eq:G:proj} proved below (see also~\cite{Di-Pietro.Droniou:17}).

\paragraph{Local potential reconstruction.}

We define 
$\drec{k+1}{T} : \dU{T}{k} \to \PP^{k+1}(T)$ such that, for all $\u v_T \in \dU{T}{k}$,
\bse \label{eq:rT}
\ba
(\GRAD \drec{k+1}{T} \u v_T,\GRAD w)_{T} &= (\GRAD v_T,\GRAD w)_{T} + (v_{\partial T}-v_T,\GRAD w \cdot \b n_{\partial T})_{\partial T} \qquad \forall w \in  \PP^{k+1}(T),\\
  (\drec{k+1}{T} \u v_T-v_T,1)_{T} &= 0.
\ea
\ese
{Using the elliptic projector identity (see \cite[(2.14)]{Di-Pietro.Droniou:20}), we obtain for every} $v \in W^{1,1}(T)$,
\bae
\drec{k+1}{T}\Interp{T}{k}v = \proj{T}{1,k+1} v.
\eae
The higher-degree potential reconstruction is used in the stabilization through the boundary residual, where it contributes to polynomial consistency of order $k+1$.

\paragraph{Boundary residual operator.}

We define 
$\dfbres{k}{\partial T} : \dU{T}{k} \rightarrow \PP^{k}(\F_T)$ such that, for all $\u v_T \in \dU{T}{k}$,
   \bae\label{eq:dfbres}
     \dfbres{k}{\partial T} \u v_T \eq
     \frac{1}{h_T}\proj{\partial T}{k}\left[
       v_T-v_{\partial T}+(I-\proj{T}{\widetilde k})(\drec{k+1}{T} \u v_T)
       \right].
   \eae
\br[Merging of the equal- and mixed-order formulations]
 {In the mixed-order case $\widetilde k=k+1$,} for all $T \in \T_h$ and all $\u v_T \in \dU{T}{k}$, {since $\drec{k+1}{T} \u v_T \in \PP^{k+1}(T)$,} {the contribution $(I-\proj{T}{\widetilde k})(\drec{k+1}{T} \u v_T)$ in \eqref{eq:dfbres} is zero} and we get the usual mixed-order boundary residual operator. 
\er

The boundary residual quantifies the mismatch between cell and face values, with the factor $h_T^{-1}$ giving it the scaling of a gradient. It vanishes when applied to the interpolant of any polynomial of degree at most $k+1$; see~\eqref{eq:sT:polynomial.consistency}.

{We write} $\dgrad{k}{h}$ and $\dfbres{k}{\partial h}$ {for the global counterparts of the local operators} $\dgrad{k}{T}$ and $\dfbres{k}{\partial T}$. 

\subsection{Stabilization function}

\paragraph{Global skeleton scalar product.}
For any $q\in\Plogpm(\Om)$ and $w,v$ defined separately on each cell boundary (the broken skeleton $\partial\mathcal T_h^{\rm br}:=\bigsqcup_{T\in\T_h}\partial T$), with $w|_{\partial T}\in L^{q'(\cdot)}(\partial T)$ and $v|_{\partial T}\in L^{q(\cdot)}(\partial T)$, we set
\bae
(w,v)_{\partial h} \eq \sum_{T \in \T_h}h_T(w,v)_{\partial T}.
\eae

\paragraph{Stabilization function.}
{The stabilization form is introduced as} $\cst{s}_h : \dU{h}{k} \times \dU{h}{k} \to \RR$ such that, for all $\u v_h,\u w_h \in \dU{h}{k}$,
\bse
\bae\label{eq:sh}
  \cst{s}_h(\u w_h, \u v_h) 
  \eq (\stabfl_{\partial h}(\cdot,\dfbres{k}{\partial h}\u w_h),\dfbres{k}{\partial h}\u v_h)_{\partial h} = \sum_{T\in\T_h}\cst{s}_T(\u w_T,\u v_T),
\eae
with for all $T \in \T_h$, 
\bae\label{eq:sT}
\cst{s}_T(\u w_T,\u v_T) \eq h_T(\stabfl_{\partial T}(\cdot,\dfbres{k}{\partial T}\u w_T),\dfbres{k}{\partial T}\u v_T)_{\partial T},
\eae
\ese
{Here, for every $F \in \F_T$,} $(\stabfl_{\partial h})\res{F} = (\stabfl_{\partial T})\res{F} = \stabfl_{F}$, and we consider for all $F \in \F_h$, the stabilization flux function $\stabfl_F : F \times \RR \to \RR$ satisfying the following assumption. 

 \bas[Stabilization flux function]\label{ass:sT}
{On each $F \in \F_h$, the map} $\stabfl_F : F \times \RR \to \RR$ is measurable and there exists $\delta_F \in L^{p(\cdot)}(F;[0,+\infty))$ such that for all $v,w \in \RR$ with $v \ne w$ and almost every $\b x \in F$,
   \bse\label{eq:S:holder.continuity.strong.monotonicity}
     \ba 
       |\stabfl_F(\b x,w)-\stabfl_F(\b x,v)| &\lesssim \sigma_\cst{c}\left(\delta_F(\b x)^{\sob(\b x)}+|w|^{\sob(\b x)}+|v|^{\sob(\b x)}\right)^\frac{\sob(\b x)-2}{\sob(\b x)}| w-v |,\label{eq:S:holder.continuity}  \\
       \left(\stabfl_F(\b x,w)-\stabfl_F(\b x,v)\right) \left(w-v\right) &\gtrsim \sigma_\cst{m}\left(\delta_F(\b x)^{\sob(\b x)}+|w|^{\sob(\b x)}+|v|^{\sob(\b x)}\right)^\frac{\sob(\b x)-2}{\sob(\b x)}|w-v|^{2},\label{eq:S:strong.monotonicity}\\
       \stabfl_F(\b x,0) &= 0. \label{eq:S:zero}
  \ea
\ese
\eas

 \begin{example}[Stabilization flux function]
   {An admissible stabilization flux under Assumption~\ref{ass:sT} is specified, for all $F \in \F_h$, $w \in \RR$, and $\b x \in F$, by}
   \bae\label{eq:sF}
     \stabfl_F(\b x,w) \eq \mu_F(\b x)\left(\delta_F(\b x)^{a(\b x)}+|w|^{a(\b x)}\right)^\frac{\sob(\b x)-2}{a(\b x)}w,
   \eae
   {where $\mu_F:F\to(0,\infty)$ is measurable and satisfies $\sigma_\cst{m}\lesssim\mu_F(\b x)\lesssim\sigma_\cst{c}$ for almost every $\b x\in F$. Here $a:F\to[a_{\min},a_{\max}]\subset(0,\infty)$ is measurable, with uniform positive lower and finite upper bounds; choosing $a\equiv2$ recovers the numerical Carreau--Yasuda stabilization.}
   {The function $\stabfl_F$ then has the $p(\cdot)$-power-framed structure required in} \eqref{eq:S:holder.continuity.strong.monotonicity}; see \cite[Appendix A]{Botti.Castanon-Quiroz.ea:20} for a proof.
 \end{example}

 \br[Weight and degeneracy of the stabilization flux function]
 {For each $F \in \F_h$, the weight of $\stabfl_F$ may be adjusted, for example by taking} $\mu_F$ equal to a convex combination $\theta \sigma_\cst{m}+ (1-\theta)\sigma_\cst{c}$ with $\theta \in [0,1]$. A factor depending on the polynomial degree $k$ may also be included in $\mu_F$ to improve numerical convergence; the choice $\mu_F=2(k+1)(k+2)$ is used in our experiments. The implicit comparison constants may then depend on $k$, which is fixed in the analysis.
 {The face degeneracy parameter may always be set to} $\delta_F = 1$, even if $\delta = 0$. Indeed, the stabilization is a numerical device, and $\stabfl_F$ need not reproduce the constitutive flux $\stress$ used in the bulk.
 \er

 \subsection{Discrete problem}

 {The discrete formulation consists in finding} $\u u_h \in \dU{h,0}{k}$ such that
 \bae
   \label{eq:lr.discrete} \cst{a}_h(\u u_h,\u v_h) = (f,v_h)_\Om \qquad \forall \u v_h \in \dU{h,0}{k}.
 \eae
 {The \emph{discrete diffusion function} $\cst{a}_h : \dU{h}{k} \times \dU{h}{k} \to \RR$ is specified for all $\u v_h,\u w_h \in \dU{h}{k}$ by}
\bse\label{eq:ah}
\bae
  \cst{a}_h(\u w_h, \u v_h) 
  \eq \sum_{T\in\T_h}\cst{a}_T(\u w_T,\u v_T),
\eae
with for all $T \in \T_h$,
\bae\label{eq:aT}
\cst{a}_T(\u w_T,\u v_T) \eq (\stress(\cdot,\dgrad{k}{T} \u w_T),\dgrad{k}{T} \u v_T)_T+\cst{s}_T(\u w_T,\u v_T).
\eae
\ese

For each cell $T\in\T_h$, let $(\delta_{\partial T})\res{F}=\delta_F$ for every $F\in\F_T$, and define the degeneracy quantity
\bae\label{eq:degeneracy:Dh}
\mathcal D_h\eq \rho_{L^{\sob(\cdot)}(\Omega)}(\delta) + \sum_{T\in\T_h}h_T\rho_{L^{\sob(\cdot)}(\partial T)}(\delta_{\partial T}).
\eae
\br[Uniform boundedness of $\mathcal D_h$]\label{rem:degeneracy:Dh}
The condition $\sup_{h\in\mathfrak H}\mathcal D_h<+\infty$ used in Section~\ref{sec:error.estimate} is mild. In particular, the admissible choice $\delta_F\equiv1$ on every face gives
\bae 
 \sum_{T\in\T_h}h_T\rho_{L^{p(\cdot)}(\partial T)}(1)
 =\sum_{T\in\T_h}h_T|\partial T|\lesssim|\Omega|,
\eae
so that $\sup_{h\in\mathfrak H}\mathcal D_h<+\infty$ follows directly from $\delta\in L^{p(\cdot)}(\Omega)$ and mesh regularity. This constant choice is therefore sufficient whenever no specific stabilization degeneracy is required.
\er
{ \paragraph{Discrete solvability.} The following result, proved in the companion convergence analysis, records the well-posedness and stability of the discrete formulation.
}\bt[{ Discrete well-posedness and variable-exponent stability bound}]\label{thm:well-posedness}
{Assumptions~\ref{ass:stress} and~\ref{ass:sT} ensure existence and uniqueness of a solution} $\u u_h\in\dU{h,0}{k}$ to \eqref{eq:lr.discrete}. Moreover,
\bae\label{eq:discrete.solution:bounds}
&\rho_{\sob(\cdot),\cst{G},h}(\u u_h)\lesssim\max(1,\sigma_\cst{m}^{-p'_\petm},\sigma_\cst{m}^{-p'_\petp})
\max\!\left(
\rho_{L^{\sob'(\cdot)}(\Om)}(f)^{p'_\petp/p'_\petm},
\rho_{L^{\sob'(\cdot)}(\Om)}(f)^{p'_\petm/p'_\petp},
\mathcal D_h
\right).
\eae
{The estimate is uniform in the mesh size and in all degeneracy parameters, and does not hide any dependence on $\sigma_\cst{m}$ or $\sigma_\cst{c}$. Its implicit factor may depend on $p_\petm$, $p_\petp$, $p_\cst{log}$, $\Omega$, $k$, and $\kappa$.}
\et

\section{Variable-exponent discrete norms and inequalities}\label{sec:discrete.norms}

We introduce the discrete seminorms associated with the hybrid unknowns, recall their equivalence, establish an HHO discrete Lebesgue embedding, and record the local trace inequality needed for the error analysis.

\paragraph{Variable-exponent discrete seminorm.}

{Given $q \in \Plogp(\Om)$, we equip $\dU{h}{k}$ with} the $W^{1,q(\cdot)}(\Omega)$-like modular $\rho_{q(\cdot),h}$ such that, for all $\u v_h \in \dU{h}{k}$,
\bse\label{eq:dseminorm:T:mod}
 \bae
    \rho_{q(\cdot),h}(\u v_h) \eq \sum_{T \in \T_h}\rho_{q(\cdot),T}(\u v_T),
 \eae
  with for all $T \in \T_h$,
  \bae
      \rho_{q(\cdot),T}(\u v_T) \eq \rho_{\tL^{q(\cdot)}(T)}(\GRAD v_T)+ h_T\rho_{L^{q(\cdot)}(\partial T)}(h_T^{-1}(v_{\partial T} - v_T)),
\eae
\ese 
{with the corresponding Luxemburg seminorms defined via \eqref{eq:luxembourg}} $\| {\cdot} \|_{q(\cdot),h}$ and $\| {\cdot} \|_{q(\cdot),T}$ for all $T\in\T_h$.

\paragraph{Reconstructed-gradient discrete seminorm.}

{For $q\in\Plogp(\Om)$, the reconstructed-gradient HHO modular is given by}
\bse\label{eq:dseminorm:G:mod}
\bae
 \rho_{q(\cdot),\cst{G},h}(\u v_h)
 \eq \sum_{T\in\T_h}\rho_{q(\cdot),\cst{G},T}(\u v_T),
\eae
with, for all $T\in\T_h$,
\bae
 \rho_{q(\cdot),\cst{G},T}(\u v_T)
 \eq
 \rho_{\tL^{q(\cdot)}(T)}(\dgrad{k}{T}\u v_T)
 +h_T\rho_{L^{q(\cdot)}(\partial T)}
 \big(\dfbres{k}{\partial T}\u v_T\big),
\eae
\ese
{with the corresponding Luxemburg seminorms defined via \eqref{eq:luxembourg}} $\|{\cdot}\|_{q(\cdot),\cst{G},h}$ and $\|{\cdot}\|_{q(\cdot),\cst{G},T}$ for all $T\in\T_h$.

\paragraph{Element-mean discrete seminorm.}

{For $q\in\Plogp(\Om)$, the element-mean discrete modular is introduced as}
\bse
 \bae
    \rho_{q_h,h}(\u v_h) \eq \sum_{T\in\T_h}\|\u v_T\|_{q_T,T}^{q_T},
 \eae
  {where, on every} $T \in \T_h$, $q_T \eq |T|^{-1}\int_T q(\b x)\,\cst d\b x$, and 
  \bae
      \|\u v_T\|_{q_T,T}^{q_T} &\eq \|\GRAD v_T\|_{\tL^{q_T}(T)}^{q_T}+ h_T\|h_T^{-1}(v_{\partial T} - v_T)\|_{L^{q_T}(\partial T)}^{q_T},
\eae
\ese
{with the corresponding Luxemburg seminorms defined via \eqref{eq:luxembourg}} $\| {\cdot} \|_{q_h,h}$ and $\| {\cdot} \|_{q_T,T}$ for all $T \in \T_h$.

\subsection{Equivalence between discrete seminorms}

\bpr[Local and global equivalence of the discrete seminorms]\label{pr:discrete_norms_equiv}
For all $q\in\Plogp(\Om)$, all $T \in \T_h$ and all $\u v_T \in \dU{T}{k}$,
\bse
\ba
     \| \u v_T \|_{q(\cdot),T}
     \eqsim \| \u v_T \|_{q(\cdot),\cst{G},T}
     \eqsim \| \u v_T \|_{q_T,T},
\ea
\ese
and, for all $\u v_h \in \dU{h}{k}$,
\bse\label{eq:disc:mod-norm:equiv}
\ba
    \| \u v_h \|_{q(\cdot),h}
    \eqsim \| \u v_h \|_{q(\cdot),\cst{G},h}
    \eqsim \| \u v_h \|_{q_h,h}.
\ea
\ese
{The implicit factors are determined only by $\Om$, $k$, $d$, $\kappa$,} $q_\cst{log}$, $q_\petm$, and $q_\petp$.
\epr

\noindent\emph{Proof.} See \cite{harnist-romain-conv}.

\subsection{HHO discrete Lebesgue embedding}\label{subsec:discrete.lebesgue}

\bl[HHO discrete Lebesgue embedding]\label{lem:error:local-change-exp}
Let $p,q\in\Plogpm(\Omega)$ and, for $T\in\T_h$, set
\bae
 p_T=\frac{1}{|T|}\int_T p(\b x)\,\cst{d}\b x,
 \qquad
 q_T=\frac{1}{|T|}\int_T q(\b x)\,\cst{d}\b x.
\eae
Then, for every $\u v_T\in\dU{T}{k}$,
\bae\label{eq:error:local-change-exp}
 \|\u v_T\|_{q(\cdot),T}
 \eqsim
 h_T^{d(1/q_T-1/p_T)}
 \|\u v_T\|_{p(\cdot),T}.
\eae
The hidden constants depend only on $\Omega$, $k$, $d$, $\kappa$, $p_\cst{log}$, $p_\petm$, $p_\petp$, $q_\cst{log}$, $q_\petm$, and $q_\petp$.
\el

\bp
This is the HHO counterpart of Corollary~\ref{cor:discrete-lebesgue-emb}.
By Proposition~\ref{pr:discrete_norms_equiv}, it is enough to compare the
local seminorms with the frozen exponents $p_T$ and $q_T$.
Applying Corollary~\ref{cor:discrete-lebesgue-emb} with the constant
exponents $p_T$ and $q_T$ to the polynomial $\GRAD v_T$ gives
\bae
 \|\GRAD v_T\|_{\tL^{q_T}(T)}
 \eqsim h_T^{d(1/q_T-1/p_T)}
 \|\GRAD v_T\|_{\tL^{p_T}(T)}.
\eae
On every $F\in\F_T$, the same polynomial Lebesgue embedding on the affine
hull of $F$ (equivalently, \eqref{eq:lebesgue.emb:fixed} in dimension
$d-1$) gives
\bae
 h_T^{1/q_T}\|h_T^{-1}(v_F-v_T)\|_{L^{q_T}(F)}
 \eqsim
 h_T^{d(1/q_T-1/p_T)}
 h_T^{1/p_T}\|h_T^{-1}(v_F-v_T)\|_{L^{p_T}(F)}.
\eae
Using the uniform mesh regularity (in particular, the uniform control of the number of faces and $h_F\eqsim h_T$), combining the cell and face contributions and using again
Proposition~\ref{pr:discrete_norms_equiv} proves
\eqref{eq:error:local-change-exp}.
\ep

\subsection{Continuous local trace inequality}

\bpr[Continuous local trace inequality]\label{pr:cont:trace:ineq}
Let $q\in\Plogpm(\Om)$. We have for all $T\in\T_h$ and all  $v \in W^{1,q(\cdot)}(T)$,
\bae\label{eq:cont:trace:ineq}
\|v\|_{L^{q(\cdot)}(\partial T)} \lesssim \|h_T^{-1/q(\cdot)}v\|_{L^{q(\cdot)}(T)}+h_T\|h_T^{-1/q(\cdot)}\GRAD v\|_{\tL^{q(\cdot)}(T)}.
\eae
\epr

\noindent\emph{Proof.} See \cite{harnist-romain-conv}.

By the discrete Sobolev--Poincar\'e inequality proved in \cite{harnist-romain-conv},
$\|{\cdot}\|_{q(\cdot),h}$ is a norm on $\dU{h,0}{k}$ for every $q\in\Plogpm(\Omega)$; by Proposition~\ref{pr:discrete_norms_equiv}, the same holds for the equivalent seminorms $\|{\cdot}\|_{q_h,h}$ and $\|{\cdot}\|_{q(\cdot),\cst{G},h}$.
In particular, we equip $\dU{h,0}{k}$ with the norm $\|{\cdot}\|_{p(\cdot),h}$.

\section{Projector approximation and discrete operators}\label{sec:properties}

{To streamline the estimates, we suppress generic multiplicative constants; for any} $a,b\ge0$, we write henceforth $a\lesssim b$ (resp. $a\gtrsim b$) for the inequality $a\le Cb$ (resp. $a\ge Cb$) where $C \in (0,\infty)$ is independent of $h$, the functions $u$, $\stress$, $f$ in \eqref{eq:intro.problem}, the parameters $\delta,\sigma_\cst{c},\sigma_\cst{m}$ in Assumption \ref{ass:stress}, and, for local inequalities, of any mesh element $T \in \T_h$ or face $F \in \F_h$ on which the inequality holds, and $\delta_F$ in Assumption \ref{ass:sT}.
{We further use $a\eqsim b$ to express} $C^{-1}b \le a \le C b$, which is equivalent to $b \lesssim a \lesssim b$.
{Further dependencies of the implicit constants are stated locally.}

\subsection{Variable-exponent approximation properties of the projectors}
\label{sec:projector.approx}

The projector estimates needed in the subsequent analysis are established in this section, building on the abstract projector approach of \cite[Lemmas~2.1--2.2]{Di-Pietro.Droniou:17*1} and the $L^2$-projector approximation result of \cite[Theorem~1.45]{Di-Pietro.Droniou:20}. The proof of Lemma~\ref{lem:bounded-projector-approx} uses Lemma~\ref{lem:john-patches} and is given after it; every other statement is followed by its proof.

\bl[Approximation for bounded projectors]\label{lem:bounded-projector-approx}
Suppose that $U$ is either star-shaped with respect to a ball of radius $\vartheta h_U$, or is the union of a uniformly bounded number of shape-regular conforming simplices of diameter comparable to $h_U$, with a connected full-facet adjacency graph. Let $p(\cdot)\in\Plogpm(U)$, and let $\Pi_\ell:L^{p(\cdot)}(U)\to\PP^\ell(U)$ be a projector reproducing $\PP^\ell(U)$ and such that
\bae\label{eq:abstract-projector-Lp-stability}
 \|\Pi_\ell w\|_{L^{p(\cdot)}(U)}\le C_\Pi\|w\|_{L^{p(\cdot)}(U)}
 \qquad\forall w\in L^{p(\cdot)}(U).
\eae
Assume moreover that the polynomial inverse inequality $|w|_{W^{m,p(\cdot)}(U)}\le C_{\rm I}h_U^{-m}\|w\|_{L^{p(\cdot)}(U)}$ holds for all $w\in\PP^\ell(U)$ and all integers $m\ge1$.
{Under these assumptions, for each integer $s\in\{1,\ldots,\ell+1\}$,} every $m\in\{0,\ldots,s\}$, and every $v\in W^{s,p(\cdot)}(U)$,
\bae\label{eq:abstract-projector-approx}
 |v-\Pi_\ell v|_{W^{m,p(\cdot)}(U)}
 \lesssim h_U^{s-m}|v|_{W^{s,p(\cdot)}(U)}.
\eae
The hidden constant depends on $d$, $\ell$, $s$, $p_\petm$, $p_\petp$, $p_\cst{log}$, $C_\Pi$, $C_{\rm I}$, and the relevant uniform geometric parameter (either $\vartheta$ or the simplicial submesh regularity bounds).
\el

\bl[Uniform John geometry of simplicial patches]\label{lem:john-patches}
Let $U\subset\mathbb R^d$ be an open set of the form $U=\operatorname{int}\bigl(\bigcup_{K\in\mathcal S_U}K\bigr)$, where $\mathcal S_U$ comprises at most $N$ conforming closed simplices whose full-facet adjacency graph is connected. Assume that these simplices are uniformly shape-regular and satisfy $h_K\eqsim h_U$ for every $K\in\mathcal S_U$, with constants independent of $U$. Then $U$ is a John domain, with a John constant depending only on $d$, $N$, and these uniform geometric bounds.
\el

\begin{figure}[h]
\centering
\begin{tikzpicture}[scale=3.4,line join=round,line cap=round,every node/.style={font=\scriptsize}]
  \coordinate (a) at (0,0);
  \coordinate (b) at (1,0);
  \coordinate (c) at (2,0);
  \coordinate (d) at (3,0);
  \coordinate (e) at (0,1);
  \coordinate (f) at (1,1);
  \coordinate (g) at (2,1);
  \coordinate (h) at (3,1);
  \coordinate (i) at (0,2);
  \coordinate (j) at (1,2);
  \coordinate (k) at (2,2);
  \coordinate (l) at (3,2);

  \fill[gray!12] (g)--(k)--(l)--cycle;
  \fill[gray!22] (g)--(h)--(l)--cycle;
  \fill[gray!12] (g)--(h)--(d)--cycle;
  \fill[gray!22] (g)--(c)--(d)--cycle;
  \fill[gray!12] (b)--(c)--(g)--cycle;
  \fill[gray!22] (b)--(f)--(g)--cycle;
  \fill[gray!12] (a)--(b)--(f)--cycle;
  \fill[gray!22] (a)--(e)--(f)--cycle;
  \fill[gray!12] (e)--(f)--(i)--cycle;
  \fill[gray!22] (f)--(j)--(i)--cycle;
  \draw[black,very thick] (a)--(d)--(l)--(k)--(g)--(f)--(j)--(i)--cycle;
  \draw[black,thin] (a)--(b) (a)--(e) (a)--(f) (b)--(c) (b)--(f) (b)--(g) (c)--(d) (c)--(g)
                    (d)--(g) (d)--(h) (e)--(f) (e)--(i) (f)--(g) (f)--(i) (f)--(j) (g)--(h)
                    (g)--(k) (g)--(l) (h)--(l) (i)--(j) (k)--(l);

  \coordinate (z0) at (barycentric cs:g=1,k=1,l=1);
  \coordinate (z1) at (barycentric cs:g=1,h=1,l=1);
  \coordinate (z2) at (barycentric cs:g=1,h=1,d=1);
  \coordinate (z3) at (barycentric cs:g=1,c=1,d=1);
  \coordinate (z4) at (barycentric cs:b=1,c=1,g=1);
  \coordinate (z5) at (barycentric cs:b=1,f=1,g=1);
  \coordinate (z6) at (barycentric cs:a=1,b=1,f=1);
  \coordinate (z7) at (barycentric cs:a=1,e=1,f=1);
  \coordinate (z8) at (barycentric cs:e=1,f=1,i=1);
  \coordinate (z9) at (barycentric cs:f=1,j=1,i=1);
  \coordinate (y0) at (barycentric cs:g=1,l=1);
  \coordinate (y1) at (barycentric cs:g=1,h=1);
  \coordinate (y2) at (barycentric cs:d=1,g=1);
  \coordinate (y3) at (barycentric cs:c=1,g=1);
  \coordinate (y4) at (barycentric cs:b=1,g=1);
  \coordinate (y5) at (barycentric cs:b=1,f=1);
  \coordinate (y6) at (barycentric cs:a=1,f=1);
  \coordinate (y7) at (barycentric cs:e=1,f=1);
  \coordinate (y8) at (barycentric cs:f=1,i=1);
  \coordinate (x) at (barycentric cs:g=0.15,k=0.6,l=0.25);

  \draw[black,line width=1.1pt,-{Stealth[length=2mm,width=1.6mm]}]
    (x)--(z0)--(y0)--(z1)--(y1)--(z2)--(y2)--(z3)--(y3)--(z4)--(y4)--(z5)--(y5)--(z6)--(y6)--(z7)--(y7)--(z8)--(y8)--(z9);
  \fill[red] (x) circle (0.6pt);
  \foreach \P in {z0,z1,z2,z3,z4,z5,z6,z7,z8,z9} {\fill[red] (\P) circle (0.6pt);}
  \foreach \P in {y0,y1,y2,y3,y4,y5,y6,y7,y8} {\filldraw[fill=white,draw=red,line width=0.5pt] (\P) circle (0.6pt);}

  \node at (1.50,1.50) {{\Large $U$}};
  \node at (2.70,1.50) {$\b y_{K_{0},K_{1}}$};
  \node at (2.36,1.08) {$\b y_{K_{1},K_{2}}$};
  \node at (2.70,0.50) {$\b y_{K_{2},K_{3}}$};
  \node at (2.16,0.59) {$\b y_{K_{3},K_{4}}$};
  \node at (1.70,0.50) {$\b y_{K_{4},K_{5}}$};
  \node at (1.18,0.40) {$\b y_{K_{5},K_{6}}$};
  \node at (0.70,0.50) {$\b y_{K_{6},K_{7}}$};
  \node at (0.64,1.08) {$\b y_{K_{7},K_{8}}$};
  \node at (0.70,1.50) {$\b y_{K_{8},K_{9}}$};
  \node at (2.27,1.59) {$\b z_{K_{0}}$};
  \node at (2.84,1.20) {$\b z_{K_{1}}$};
  \node at (2.73,0.59) {$\b z_{K_{2}}$};
  \node at (2.30,0.26) {$\b z_{K_{3}}$};
  \node at (1.64,0.26) {$\b z_{K_{4}}$};
  \node at (1.36,0.74) {$\b z_{K_{5}}$};
  \node at (0.64,0.26) {$\b z_{K_{6}}$};
  \node at (0.27,0.59) {$\b z_{K_{7}}$};
  \node at (0.28,1.42) {$\b z_{K_{8}}$};
  \node at (0.70,1.74) {$\b x_\ast{=}\b z_{K_\ast}$};
  \node at (2.33,1.85) {$\b x$};
  \node at (1.86,0.30) {{\large $\b\gamma_{\b x}$}};
\end{tikzpicture}
\caption{A non-star-shaped, U-shaped simplicial patch $U$ made of five congruent squares, each split into two triangles ($K_0,\dots,K_9$), and the polygonal John curve $\b\gamma_{\b x}$ from $\b x\in K_0$ to $\b x_\ast=\b z_{K_\ast}$, $K_\ast=K_9$. It passes through the interior points $\b z_{K_i}$ (filled dots) and the barycenters $\b y_{K_i,K_{i+1}}$ of the shared full facets (open circles).}
\label{fig:john-patches}
\end{figure}
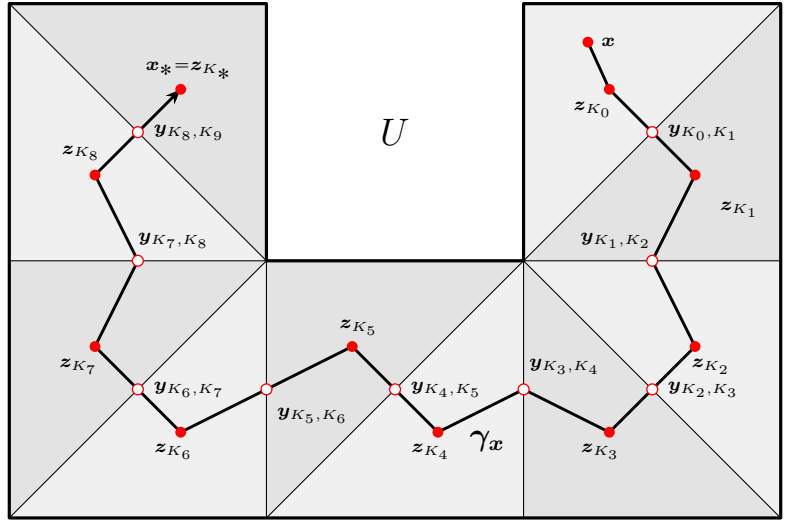

\bp
If $d=1$, the connected open set $U$ is an interval $(a,b)$. Taking $\b x_\ast=(a+b)/2$ and joining any $\b x\in U$ to $\b x_\ast$ by a segment parametrized by arc length from $\b x$ gives $\operatorname{dist}(\b\gamma_{\b x}(t),\partial U)\ge t$ for $0\le t\le|\b x-\b x_\ast|$. Thus $U$ is a John domain with John constant $1$. In the rest of the proof, assume $d\ge2$.

Scaling by $h_U^{-1}$, we may assume $h_U=1$. Let $\mathcal G$ be the full-facet adjacency graph of the simplices, and choose a spanning tree rooted at a simplex $K_\ast$. Its number of edges is at most $N-1$.

\emph{Step 1: uniform interior balls.} Shape regularity and $h_K\eqsim1$ give a number $r_0>0$, depending only on the stated geometric parameters, and points $\b z_K\in\operatorname{int}K$, $K\in\mathcal S_U$, such that
\bae
B(\b z_K,r_0)\subset\operatorname{int}K\subset U.
\eae
For every pair of simplices $K,K'\in\mathcal S_U$ sharing a full facet $F$, let $\b y_{K,K'}$ be the barycenter of $F$. The uniform lower bounds on the simplex altitudes and on the relative interior size of their common facet yield a number $r_1>0$, independent of the pair, for which
\bae
B(\b y_{K,K'},r_1)\subset\operatorname{int}(K\cup K')\subset U.
\eae
In particular, the interface points $\b y_{K,K'}$ belong to $U$, even though they lie on the individual simplex boundaries.

\emph{Step 2: a quantitative transfer estimate.} Put $\theta:=\min\{\tfrac12,\tfrac{r_1}{2}\}$ and $\delta_0:=\min\{\theta r_0,\tfrac{r_1}{2}\}>0$. For adjacent $K,K'$, consider the segment $\b\xi(s)=(1-s)\b z_K+s \b y_{K,K'}$, $0\le s\le1$. If $s\le1-\theta$, convexity of $K$ gives
\bae
B\bigl(\b\xi(s),(1-s)r_0\bigr)\subset K,
\qquad
\operatorname{dist}(\b\xi(s),\partial U)\ge\theta r_0.
\eae
If $s\ge1-\theta$, then $|\b z_K-\b y_{K,K'}|\le h_K\le1$ and the distance-to-the-boundary function is $1$-Lipschitz; hence
\bae
\operatorname{dist}(\b\xi(s),\partial U)
\ge\operatorname{dist}(\b y_{K,K'},\partial U)-|\b\xi(s)-\b y_{K,K'}|
\ge r_1-\theta\ge\frac{r_1}{2}.
\eae
Consequently, every point on $[\b z_K,\b y_{K,K'}]$ is at distance at least $\delta_0$ from $\partial U$. The same conclusion holds for $[\b y_{K,K'},\b z_{K'}]$.

\emph{Step 3: construction and John estimate.} Fix $\b x\in U$, choose $K_0\in\mathcal S_U$ with $\b x\in K_0$, and follow the unique spanning-tree chain $K_0,K_1,\ldots,K_m=K_\ast$. Join $\b x$ to $\b z_{K_0}$ by a straight segment and then, for each $i<m$, concatenate the segments
\bae
[\b z_{K_i},\b y_{K_i,K_{i+1}}]\quad\text{and}\quad[\b y_{K_i,K_{i+1}},\b z_{K_{i+1}}].
\eae
This produces a polygonal path $\b\gamma_{\b x}\subset U$ ending at the fixed point $\b x_\ast:=\b z_{K_\ast}$; see Figure~\ref{fig:john-patches}. Since each segment has length at most $1$, its total length $L_{\b x}$ satisfies $L_{\b x}\le1+2m\le2N-1$. Parametrize $\b\gamma_{\b x}$ by arc length $t\in[0,L_{\b x}]$, starting at $\b x$, and write $t_0:=|\b x-\b z_{K_0}|$. On the first segment, for $0\le t\le t_0$, set $s=t/t_0$ when $t_0>0$. Convexity of $K_0$ implies
\bae
B\bigl((1-s)\b x+s \b z_{K_0},s r_0\bigr)\subset K_0,
\qquad
\operatorname{dist}(\b\gamma_{\b x}(t),\partial U)\ge s r_0\ge r_0 t,
\eae
where the last inequality uses $t_0\le h_{K_0}\le1$. If $t_0=0$, there is no first segment. On the remainder of the curve, Step~2 yields $\operatorname{dist}(\b\gamma_{\b x}(t),\partial U)\ge\delta_0$. Since $t\le L_{\b x}\le2N-1$, we obtain in both cases
\bae
\operatorname{dist}(\b\gamma_{\b x}(t),\partial U)\ge c_J t
\quad(0\le t\le L_{\b x}),
\qquad
c_J:=\min\left(r_0,\frac{\delta_0}{2N-1}\right)>0.
\eae
The constant $c_J$ is independent of $U$ and of the mesh size. Rescaling completes the proof.
\ep

\bp[Proof of Lemma~\ref{lem:bounded-projector-approx}]
The argument is the variable-exponent counterpart of \cite[Lemma~2.1]{Di-Pietro.Droniou:17*1}. The scaled Poincar\'e--Wirtinger inequality follows from \cite[Theorem~8.2.4(b)]{diening2011lebesgue}, since a domain star-shaped with respect to a uniformly sized ball is uniformly John and the simplicial alternative is covered by Lemma~\ref{lem:john-patches}. In both cases, the John constants are uniform with respect to the mesh cell. It gives
\bae\label{eq:local-PW-projectors}
 \|w-\langle w\rangle_U\|_{L^{p(\cdot)}(U)}
 \lesssim h_U\|\GRAD w\|_{\tL^{p(\cdot)}(U)},
 \qquad
 \langle w\rangle_U \eq \frac{1}{|U|}\int_U w\,\cst{d}\b x.
\eae
Let $Q_U^{s-1}v\in\PP^{s-1}(U)$ be the moment polynomial determined by
\bae\label{eq:moment-Taylor-projector}
 \int_U \partial_\alpha(v-Q_U^{s-1}v)\,\cst{d}\b x=0
 \qquad (|\alpha|\le s-1).
\eae
Applying \eqref{eq:local-PW-projectors} successively to the derivatives of
$z\eq v-Q_U^{s-1}v$ yields
\bae\label{eq:variable-BH-projectors}
 |z|_{W^{j,p(\cdot)}(U)}
 \lesssim h_U^{s-j}|v|_{W^{s,p(\cdot)}(U)},
 \qquad 0\le j\le s.
\eae
Indeed, the moment conditions give zero mean for $\partial_\alpha z$ whenever
$|\alpha|\le s-1$, while the derivatives of order $s$ of $Q_U^{s-1}v$ vanish.
Since $s-1\le\ell$ and $\Pi_\ell$ reproduces $\PP^\ell(U)$,
$v-\Pi_\ell v=z-\Pi_\ell z$. For $m=0$, the stability
\eqref{eq:abstract-projector-Lp-stability} and \eqref{eq:variable-BH-projectors} give the claim directly. For $m\ge1$, the assumed polynomial inverse inequality gives
\bae
 |\Pi_\ell z|_{W^{m,p(\cdot)}(U)}
 \lesssim h_U^{-m}\|\Pi_\ell z\|_{L^{p(\cdot)}(U)}
 \lesssim h_U^{-m}\|z\|_{L^{p(\cdot)}(U)},
\eae
and \eqref{eq:variable-BH-projectors} with $j=0$ and $j=m$ concludes the proof.
\ep

\bl[Boundedness of orthogonal projectors on polynomial subspaces]\label{bl:boundedness:proj}
{Suppose that $U$ satisfies \eqref{eq:mean.exp.regularity} and} $h_U\le h_0$ for some fixed $h_0<\infty$, let $p(\cdot)\in\Plogpm(U)$, and let $\mathcal P\subset\PP^\ell(U)$. The $L^2$-orthogonal projector $\pi_{\mathcal P}:L^1(U)\to\mathcal P$ satisfies
\bae\label{eq:L2proj-variable-bounded}
 \|\pi_{\mathcal P}\varphi\|_{L^{p(\cdot)}(U)}
 \lesssim \|\varphi\|_{L^{p(\cdot)}(U)}.
\eae
{The implicit factor is uniform with respect to the choice of $\mathcal P$} and depends only on $d$, $\ell$, the regularity constant in \eqref{eq:mean.exp.regularity}, $p_\petm$, $p_\petp$, $p_\cst{log}$, and $h_0$.
\el

\bp
We adapt \cite[Lemma~2.2]{Di-Pietro.Droniou:17*1}. Let
$(\psi_i)_{i=1}^{N_{\mathcal P}}$ be an $L^2(U)$-orthonormal basis of
$\mathcal P$. The fixed-exponent
Lebesgue embedding \eqref{eq:lebesgue.emb:fixed} gives
\bae
 \|\psi_i\|_{L^\infty(U)}\lesssim |U|^{-1/2},
\eae
and therefore, since $N_{\mathcal P}\le\dim\PP^\ell(U)$, the kernel
$K_{\mathcal P}(\b x,\b y)\eq\sum_i\psi_i(\b x)\psi_i(\b y)$ satisfies
$|K_{\mathcal P}(\b x,\b y)|\lesssim |U|^{-1}$. Hence
\bae\label{eq:L2proj-pointwise-L1}
 |(\pi_{\mathcal P}\varphi)(\b x)|\lesssim |U|^{-1}\|\varphi\|_{L^1(U)}.
\eae
By the variable-exponent H\"older inequality~{\eqref{eq:holder}},
\bae\label{eq:L2proj-before-indicator}
 \|\pi_{\mathcal P}\varphi\|_{L^{p(\cdot)}(U)}
 \lesssim |U|^{-1}
 \|1\|_{L^{p(\cdot)}(U)}\|1\|_{L^{p'(\cdot)}(U)}
 \|\varphi\|_{L^{p(\cdot)}(U)}.
\eae
It remains only to bound the prefactor. Let $p_U$ and $(p')_U$ be the mean
exponents. Theorem~\ref{thm:mean_exp_equiv}, applied to the constant
polynomial $1$, gives
\bae
 \|1\|_{L^{p(\cdot)}(U)}\eqsim |U|^{1/p_U},
 \qquad
 \|1\|_{L^{p'(\cdot)}(U)}\eqsim |U|^{1/(p')_U}.
\eae
Choose any $\b x_U\in U$. Using \eqref{eq:mean.exp.diam.equiv} for $p$ and
$p'$, together with $|U|\eqsim h_U^d$, we obtain
\bae
 |U|^{1/p_U+1/(p')_U}
 \eqsim h_U^{d/p(\b x_U)+d/p'(\b x_U)}
 =h_U^d\eqsim |U|.
\eae
Thus the prefactor in \eqref{eq:L2proj-before-indicator} is uniformly
bounded, which proves \eqref{eq:L2proj-variable-bounded}.
\ep

\bpr[Variable-exponent approximation properties]\label{pr:variable-projector-approx}
Let $q\in\Plogpm(\Om)$ (in particular, $q=p$ or $q=p'$). {Fix $T\in\T_h$ and integers $1\le n\le\ell+1$ and $0\le m\le n$.} For every $v\in W^{n,q(\cdot)}(T)$,
\bae\label{eq:proj:app:norm:T}
 |v-\proj{T}{\ell}v|_{W^{m,q(\cdot)}(T)}
 \lesssim h_T^{n-m}|v|_{W^{n,q(\cdot)}(T)}.
\eae
Moreover, the scaled trace estimate
\bae\label{eq:proj:app:norm:F}
 \big\|h_T^{1/q(\cdot)}(v-\proj{T}{\ell}v)\big\|_{L^{q(\cdot)}(\partial T)}
 \lesssim h_T^{n}|v|_{W^{n,q(\cdot)}(T)}
\eae
holds. {The implicit bounds are uniform over the mesh} and depend only on $d$, $\kappa$, $\ell$, $n$, $q_\petm$, $q_\petp$, and $q_\cst{log}$.
\epr

\bp
Lemma~\ref{bl:boundedness:proj}, with $U=T$ and
$\mathcal P=\PP^\ell(T)$, gives the uniform $L^{q(\cdot)}(T)$-stability
of $\proj{T}{\ell}$. Moreover, by \cite[Definitions~1.8--1.9]{Di-Pietro.Droniou:20}, each cell $T$ is the union of a uniformly bounded number of shape-regular simplices of diameter comparable to $h_T$, and these simplices form a matching (conforming) simplicial submesh. Its full-facet adjacency graph is connected: since $T$ has connected interior, two simplex interiors can be joined by a path in $T$ that avoids the skeleton of codimension at least two, so the path crosses only shared full facets. The required polynomial inverse inequality follows from \cite[Lemma~1.28]{Di-Pietro.Droniou:20} together with the mean-exponent equivalence of Theorem~\ref{thm:mean_exp_equiv}. Lemma~\ref{lem:bounded-projector-approx} with $s=n$ therefore yields \eqref{eq:proj:app:norm:T} immediately.
For the trace estimate, set $e_T\eq v-\proj{T}{\ell}v$ and define $q_T\eq|T|^{-1}\int_T q(\b x)\,\cst{d}\b x$. By \eqref{eq:mean.exp.diam.equiv}, the following powers of $h_T$ are uniformly equivalent on $\overline T$:
\bae
h_T^{\pm1/q(\cdot)}\eqsim h_T^{\pm1/q_T}.
\eae
Hence Proposition~\ref{pr:cont:trace:ineq} gives
\bae
 \|h_T^{1/q(\cdot)}e_T\|_{L^{q(\cdot)}(\partial T)}
 &\lesssim h_T^{1/q_T}\|e_T\|_{L^{q(\cdot)}(\partial T)}\lesssim \|e_T\|_{L^{q(\cdot)}(T)}
 +h_T\|\GRAD e_T\|_{\tL^{q(\cdot)}(T)}.
\eae
Applying \eqref{eq:proj:app:norm:T} with $m=0$ and $m=1$ proves
\eqref{eq:proj:app:norm:F}.
\ep

\bl[From local approximation bounds to broken norms]\label{lem:local-to-global-projector}
Let $q\in\Plogpm(\Omega)$, let $n\ge1$ be an integer, and let
$v\in W^{n,q(\cdot)}(\T_h)$. Fix positive integers $m_1,m_2$ and
non-negative exponents $\alpha,\beta$. For every $T\in\T_h$,
let $a_T\in L^{q(\cdot)}(T;\RR^{m_1})$ and
$b_{\partial T}\in L^{q(\cdot)}(\partial T;\RR^{m_2})$ satisfy
\bse
\ba
 \|a_T\|_{L^{q(\cdot)}(T)}
 &\lesssim h_T^{\alpha}|v|_{W^{n,q(\cdot)}(T)},
 \label{eq:local-to-global-cell}\\
 \|h_T^{1/q(\cdot)}b_{\partial T}\|_{L^{q(\cdot)}(\partial T)}
 &\lesssim h_T^{\beta}|v|_{W^{n,q(\cdot)}(T)}.
 \label{eq:local-to-global-face}
\ea
\ese
Define the broken fields $a_h$ on $\Omega$ and $b_{\partial h}$ on the
broken skeleton by $a_h|_T=a_T$ and
$b_{\partial h}|_{\partial T}=b_{\partial T}$ for each $T\in\T_h$.
The weighted skeleton Luxemburg norm is defined in
\eqref{eq:error:skeleton:norm}. Then
\bse
\ba
 \|a_h\|_{L^{q(\cdot)}(\Omega)}
 &\lesssim h^{\alpha}|v|_{W^{n,q(\cdot)}(\T_h)},
 \label{eq:local-to-global-cell-bound}\\
 \|b_{\partial h}\|_{q(\cdot),\partial h}
 &\lesssim h^{\beta}|v|_{W^{n,q(\cdot)}(\T_h)}.
 \label{eq:local-to-global-face-bound}
\ea
\ese
\el

\noindent\emph{Proof.}
If $|v|_{W^{n,q(\cdot)}(\T_h)}=0$, both conclusions follow from
\eqref{eq:local-to-global-cell}--\eqref{eq:local-to-global-face}.
Otherwise, set for all $T \in \T_h$,
\bae
 t_T\eq\frac{|v|_{W^{n,q(\cdot)}(T)}}
 {|v|_{W^{n,q(\cdot)}(\T_h)}}.
\eae
The definition of the broken Luxemburg seminorm implies $0\le t_T\le1$ and
\bae
 \sum_{T\in\T_h}t_T^{q_{T,+}}
 \le\sum_{T\in\T_h}
 \rho_{\dot W^{n,q(\cdot)}(T)}
 \left(\frac{v}{|v|_{W^{n,q(\cdot)}(\T_h)}}\right)\le1.
\eae
Since $h_T\le h$, the local estimates
\eqref{eq:local-to-global-cell}--\eqref{eq:local-to-global-face}
ensure that, for a uniform $C\ge1$, the Luxemburg norms of
\bae
 \frac{a_T}{C h^{\alpha}|v|_{W^{n,q(\cdot)}(\T_h)}}
 \quad\text{and}\quad
 \frac{h_T^{1/q(\cdot)}b_{\partial T}}
 {C h^{\beta}|v|_{W^{n,q(\cdot)}(\T_h)}}
\eae
are at most $t_T$. By the modular--norm inequalities
\eqref{eq:modular-norm-equiv}, their local modulars are therefore
bounded by $t_T^{q_{T,-}}$. For the face terms, the identity
\bae
 \rho_{L^{q(\cdot)}(\partial T)}
 \left(\frac{h_T^{1/q(\cdot)}b_{\partial T}}
 {C h^{\beta}|v|_{W^{n,q(\cdot)}(\T_h)}}\right)
 =h_T\rho_{L^{q(\cdot)}(\partial T)}
 \left(\frac{b_{\partial T}}
 {C h^{\beta}|v|_{W^{n,q(\cdot)}(\T_h)}}\right)
\eae
recovers the skeleton weight in \eqref{eq:error:skeleton:norm}.
It remains to control the sum of $t_T^{q_{T,-}}$. If
$t_T\ge|T|^{1/q_{T,+}}$, the exponent-oscillation estimate
\eqref{eq:mean.exp.diam.equiv} and $|T|\eqsim h_T^d$ imply
\bae
 t_T^{q_{T,-}}
 \le |T|^{-(q_{T,+}-q_{T,-})/q_{T,+}}t_T^{q_{T,+}}
 \lesssim t_T^{q_{T,+}}.
\eae
For the remaining elements, the same oscillation bound gives
$t_T^{q_{T,-}}\lesssim|T|$. Consequently,
\bae
 \sum_{T\in\T_h}t_T^{q_{T,-}}
 \lesssim\sum_{T\in\T_h}t_T^{q_{T,+}}+\sum_{T\in\T_h}|T|
 \le1+|\Omega|.
\eae
Both normalized global modulars are thus uniformly bounded.
Increasing $C$ by a fixed factor makes each at most one, and the
Luxemburg norm definition gives
\eqref{eq:local-to-global-cell-bound}--\eqref{eq:local-to-global-face-bound}.
\hfill$\square$

\subsection{Properties of the discrete operators}

{We next collect the properties needed for the local discrete operators.}

\paragraph{Boundedness of the interpolator.}
{For constant exponents, the boundedness of $\Interp{T}{k}$ was proved in} \cite[Proposition 6.24]{Di-Pietro.Droniou:20}.

\bpr[Boundedness of the interpolator]\label{pr:interpolator-bounded}
For all $q \in \Plogpm(\Om)$, all $T \in \T_h$, and all $v \in W^{1,q(\cdot)}(T)$,
\bae\label{eq:I:boundedness:norm}
  \|\Interp{T}{k} v\|_{q(\cdot),T} \lesssim |v|_{W^{1,q(\cdot)}(T)}.
\eae
{Only $k$, $d$, $\kappa$, $q_-$, $q_+$, and $q_\cst{log}$ enter the hidden constant.}
\epr

\noindent\emph{Proof.} See \cite{harnist-romain-conv}.

\paragraph{Commutation of the gradient reconstruction.}
{The construction directly yields the identity} (see \cite[Section 7.2.5]{Di-Pietro.Droniou:20}):
for all $T \in \T_h$ and all $v\in W^{1,1}(T)$,
\bae\label{eq:G:proj}
  \dgrad{k}{T} \Interp{T}{k} v = \PROJ{T}{k}(\GRAD v).
\eae

\paragraph{Consistency of the boundary residual operator.}

{The boundary residual enjoys polynomial consistency} (cf. \cite{Di-Pietro.Droniou:20}), i.e., for all $T \in \T_h$,
 \bae\label{eq:sT:polynomial.consistency}
   \dfbres{k}{\partial T}(\Interp{T}{k} w) = 0\qquad\forall w\in\PP^{k+1}(T).
 \eae
The polynomial consistency of this residual is the counterpart of the difference-operator consistency in \cite[Proposition~2.6]{Di-Pietro.Droniou:20}.

\section{A priori error estimate}\label{sec:error.estimate}

The main estimates are stated first, together with their convergence-rate interpretation. Their proofs and the auxiliary lemmas are collected in Subsection~\ref{subsec:error.proofs}.

\subsection{Error estimates and convergence rates}\label{subsec:error.results}

{ For a measurable exponent $q(\cdot)>0$ and $a\ge0$, set}
\bae\label{eq:error:mqT}
 {\mathfrak m_{q,T}(a)}
 \eq {\max\left(a^{q_{T,-}},a^{q_{T,+}}\right)},
 \qquad
 q_{T,-}\eq\essinf_Tq,
 \quad q_{T,+}\eq\esssup_Tq.
\eae
{ For any measurable exponent $q$ satisfying $1<q_-\le q_+<\infty$, and for a non-negative family ${(r_T)_{T\in\T_h}}$, define}
\bae\label{eq:error:Ah}
 {\mathcal A_{q,h}({(r_T)_{T\in\T_h}})}
 \eq
 \inf\left\{\lambda>0:\ 
 \sum_{T\in\T_h}
 {\mathfrak m_{q,T}}\!\left(\frac{r_T}{\lambda}\right)\le1
 \right\}.
\eae

\bt[Localized discrete-norm error estimate]\label{thm:error.estimate}
Under Assumptions~\ref{ass:stress} and~\ref{ass:sT}, let $(\mathcal M_h)_{h\in\mathfrak H}$ be a regular mesh family such that $\sup_{h\in\mathfrak H}\mathcal D_h<+\infty$ (see Remark~\ref{rem:degeneracy:Dh}). Let $u$ and $\u u_h$ solve \eqref{eq:lr.weak} and \eqref{eq:lr.discrete}, respectively, and assume $u\in W^{k+2,p(\cdot)}(\T_h)$, $\stress(\cdot,\GRAD u)
 \in\tW^{1,p'(\cdot)}(\Omega)\cap\tW^{k+1,p'(\cdot)}(\T_h)$.
Then
\bae\label{eq:error.estimate.localized}
 \|\u u_h-\Interp{h}{k}u\|_{p(\cdot),h}
 \le C\Bigg[&
 \mathcal A_{p',h}\!\left(\left(h_T^{k+1}|\stress(\cdot,\GRAD u)|_{\tW^{k+1,p'(\cdot)}(T)}\right)_{T\in\T_h}\right)^{\frac1{\sobmax_\petp-1}}\\
 &+
 \mathcal A_{p,h}\!\left(\left(h_T^{k+1}|u|_{W^{k+2,p(\cdot)}(T)}\right)_{T\in\T_h}\right)^{\frac{\sobmin_\petm}{\sobmax_\petp}}
 \Bigg].
\eae
Moreover, we have the global form
\bae\label{eq:error.estimate.global}
 \|\u u_h-\Interp{h}{k}u\|_{p(\cdot),h}
 \le C\Bigg[&
 h^{\frac{k+1}{\sobmax_\petp-1}}
 |\stress(\cdot,\GRAD u)|_{\tW^{k+1,p'(\cdot)}(\T_h)}^{\frac1{\sobmax_\petp-1}}
+ h^{(k+1)\frac{\sobmin_\petm}{\sobmax_\petp}}
 |u|_{W^{k+2,p(\cdot)}(\T_h)}^{\frac{\sobmin_\petm}{\sobmax_\petp}}
 \Bigg].
\eae
The constant $C\in(0,\infty)$ is independent of $h$ and of the individual elements $T\in\T_h$ and faces $F\in\F_h$. Its dependence is restricted to $d$, $k$, $\Omega$, $\kappa$, $p_\petm$, $p_\petp$, $p_\cst{log}$, $\sigma_\cst{c}$, $\sigma_\cst{m}$, and the data bound
\bae
 \max\!\left\{1,
 \rho_{L^{p'(\cdot)}(\Omega)}(f)^{p'_\petp/p'_\petm},
 \rho_{L^{p'(\cdot)}(\Omega)}(f)^{p'_\petm/p'_\petp},
 \sup_{h\in\mathfrak H}\mathcal D_h
 \right\}.
\eae
\et

\br[Convergence rates in the general case]\label{rem:error:rates-general}
If the broken Sobolev seminorms on the right-hand side of~\eqref{eq:error.estimate.global} remain uniformly bounded along the mesh family, then Theorem~\ref{thm:error.estimate} implies
\bae
 \|\u u_h-\Interp{h}{k}u\|_{p(\cdot),h}
 =O\!\left(h^{(k+1)\alpha}\right),
 \quad
 \alpha\eq\min\left(
 \frac1{\sobmax_\petp-1},
 \frac{\sobmin_\petm}{\sobmax_\petp}
 \right)
 =
 \begin{cases}
 \dfrac{p_\petm}{2}, & p_\petp\le2,\\[2mm]
 \min\!\left(\dfrac1{p_\petp-1},\dfrac{p_\petm}{p_\petp}\right),
 & p_\petm<2<p_\petp,\\[3mm]
 \dfrac1{p_\petp-1}, & p_\petm\ge2.
 \end{cases}
\eae
For a constant exponent $p$, this gives the rates $(k+1)p/2$ for $1<p\le2$ and $(k+1)/(p-1)$ for $p\ge2$.
\er

\bt[Regime-dependent refinement for $p(\cdot)\le2$]\label{thm:error.estimate:regime}
Assume the hypotheses of Theorem~\ref{thm:error.estimate} and $p_\petp\le2$. For every $T\in\T_h$, set
\bae\label{eq:regime:D}
 \mathfrak D_T
 \eq \min\left\{
 1,
 \essinf_{\b x\in T}\big(\delta(\b x)+|\GRAD u(\b x)|\big),
 \min_{F\in\F_T}\essinf_{\b x\in F}\delta_F(\b x)
 \right\},
\eae
and define
\bae\label{eq:regime:etaT}
 \eta_T\eq
 \frac{h_T^{k+1}|u|_{W^{k+2,p(\cdot)}(T)}}
 {|T|^{1/p_T}\mathfrak D_T}
 \qquad\text{if }\mathfrak D_T>0,
\eae
with the conventions $\eta_T=+\infty$ if $\mathfrak D_T=0$ and
$|u|_{W^{k+2,p(\cdot)}(T)}>0$, and $\eta_T=0$ if both quantities vanish. Then
\bae\label{eq:error.estimate.regime}
 \|\u u_h-\Interp{h}{k}u\|_{p(\cdot),h}
 \le C\Bigg[&
 \mathcal A_{p',h}\!\left(\left(h_T^{k+1}|\stress(\cdot,\GRAD u)|_{\tW^{k+1,p'(\cdot)}(T)}\right)_{T\in\T_h}\right)\\
 &+\mathcal A_{p',h}\!\left(\left(
 \min(1,\eta_T)^{2-p_{T,-}}
 {\mathfrak m_{p-1,T}}\!\left(h_T^{k+1}|u|_{W^{k+2,p(\cdot)}(T)}\right)
 \right)_{T\in\T_h}\right)
 \Bigg].
\eae
The constant $C$ has the same dependence as in Theorem~\ref{thm:error.estimate} and is independent of $h$, $T$, $F$, and of the values of $(\mathfrak D_T)_T$ appearing explicitly above.
\et

\br[Resolved non-degenerate regime]\label{rem:error:rate}
The indicator $\eta_T$ compares the local approximation scale with the local non-degeneracy scale. In other words, $\eta_T\lesssim1$ means that the non-degenerate branch is resolved on $T$.
For instance, take $\delta_F\equiv1$ and assume
\bae
 \mathfrak D_0\eq
 \essinf_{\b x\in\Omega}
 \big(\delta(\b x)+|\GRAD u(\b x)|\big)>0.
\eae
If we additionally have
\bae
 \sup_{h\in\mathfrak H}\max_{T\in\T_h}
 |T|^{-1/p_T}|u|_{W^{k+2,p(\cdot)}(T)}<+\infty,
\eae
then $\eta_T\lesssim h_T^{k+1}$ uniformly in $T$ and the optimal order $k+1$ follows. 
If the non-degenerate regime is not resolved, Theorem~\ref{thm:error.estimate} still provides the robust rate $(k+1)p_\petm/2$ when $p_\petp\le2$.
\er

\bc[Broken-gradient error estimate]\label{c:error.estimate:grad}
Under the assumptions of Theorem~\ref{thm:error.estimate},
\bae\label{eq:error:gradient:corrected}
 |u_h-u|_{W^{1,p(\cdot)}(\T_h)}
 &\lesssim \|\u u_h-\Interp{h}{k}u\|_{p(\cdot),h}
 +h^{\widetilde k}|u|_{W^{\widetilde k+1,p(\cdot)}(\T_h)},
\eae
and we deduce convergence rates from Theorems~\ref{thm:error.estimate} and \ref{thm:error.estimate:regime} (the projection term $h^{\widetilde k}$ does not reduce the rate when $\widetilde k=k+1$, but can limit it to order $k$ when $\widetilde k=k$).
\ec

\subsection{Proofs and auxiliary lemmas}\label{subsec:error.proofs}

For a measurable function $z$ on the broken skeleton, we use the weighted modular and its Luxemburg norm
\bae\label{eq:error:skeleton:norm}
 \rho_{q(\cdot),\partial h}(z)
 &\eq \sum_{T\in\T_h}h_T\rho_{L^{q(\cdot)}(\partial T)}(z\res{\partial T}),\\
 \|z\|_{q(\cdot),\partial h}
 &\eq \inf\left\{\lambda>0:\rho_{q(\cdot),\partial h}(z/\lambda)\le1\right\}.
\eae

The following estimate extends the three-point absorption arguments used for constant-exponent nonlinear diffusion in \cite[Lemma~2.1]{BeiraoDiPietroHaile2024} and for non-Newtonian Stokes flows in \cite[Lemmas~22--23]{AntoniettiBeiraoBottiVaccaVerani2024} to spatially varying exponents and to the general fluxes covered by Assumption~\ref{ass:stress}.

\bl[Three-point perturbation estimate]\label{lem:three.point.error}
Under Assumption~\ref{ass:stress}, for every $\varepsilon\in(0,1)$ there exists $C_\varepsilon>0$ such that, for a.e. $\b x\in\Omega$ and all $\b\xi,\b\tau,\b\eta\in\RR^d$,
\bae\label{eq:three.point.bulk}
 &|\stress(\b x,\b\tau)-\stress(\b x,\b\eta)|\,|\b\xi-\b\tau|\\
 &\qquad\le \varepsilon
 \big(\stress(\b x,\b\xi)-\stress(\b x,\b\tau)\big)\cdot(\b\xi-\b\tau)
 +C_\varepsilon
 \big(\delta(\b x)^{\sob(\b x)}+|\b\tau|^{\sob(\b x)}+|\b\eta|^{\sob(\b x)}\big)^{\frac{\sob(\b x)-\sobmin(\b x)}{\sob(\b x)}}
 |\b\tau-\b\eta|^{\sobmin(\b x)}.
\eae
The scalar analogue holds for $\stabfl_F$ on every face $F$, with $\delta$ replaced by $\delta_F$. {For the bulk flux one possible uniform constant is
\bae\label{eq:three.point.constant}
 C_\varepsilon=2^{1/4}\sigma_\cst{c}
 \left(1+\frac{4^{\max(1,p_\petp-2)}\sigma_\cst{c}}{\varepsilon\sigma_\cst{m}}\right)^{\frac{1}{\min(1,p_\petm-1)}}.
\eae
For the face flux the same expression holds with its uniform effective continuity and monotonicity constants from \eqref{eq:S:holder.continuity}--\eqref{eq:S:strong.monotonicity}. Thus $C_\varepsilon$ is independent of $\delta$, $\delta_F$, $\b x$, and the three arguments of $\RR^d$.}
\el

\bp
{
Fix $\b x$ such that Assumption~\ref{ass:stress} holds, and write $p=p(\b x)$. The assertion is immediate if $|\b\xi-\b\tau|=0$ or $|\b\tau-\b\eta|=0$. Otherwise, set
\bae
 K=\left(1+\frac{4^{\max(1,p_\petp-2)}\sigma_\cst{c}}{\varepsilon\sigma_\cst{m}}\right)^{\frac1{\min(1,p_\petm-1)}}>1.
\eae
By the H\"older-continuity and strong-monotonicity assumptions \eqref{eq:ass:stress:cont} and \eqref{eq:ass:stress:mono}, respectively,
\bse
\ba
 |\stress(\b x,\b\tau)-\stress(\b x,\b\eta)|\,|\b\xi-\b\tau| &\le\sigma_\cst{c}\big(\delta(\b x)^p+|\b\tau|^p+|\b\eta|^p\big)^{\frac{p-2}{p}}
 |\b\tau-\b\eta|\,|\b\xi-\b\tau|,\label{eq:three.point.cont} \\
 \big(\stress(\b x,\b\xi)-\stress(\b x,\b\tau)\big)\cdot(\b\xi-\b\tau) &\ge\sigma_\cst{m}\big(\delta(\b x)^p+|\b\xi|^p+|\b\tau|^p\big)^{\frac{p-2}{p}}
 |\b\xi-\b\tau|^2.\label{eq:three.point.mono}
\ea
\ese
First suppose $|\b\xi-\b\tau|>K|\b\tau-\b\eta|$. If $p\ge2$, the triangle inequality implies
\bae
 \delta(\b x)^p+|\b\tau|^p+|\b\eta|^p
 &\le 2^p\big(\delta(\b x)^p+|\b\xi|^p+|\b\tau|^p\big)
 +2^{p-1}|\b\tau-\b\eta|^p\\
 &\le 2^{2p}\big(\delta(\b x)^p+|\b\xi|^p+|\b\tau|^p\big),
\eae
where we used $|\b\tau-\b\eta|<|\b\xi-\b\tau|\le |\b\xi|+|\b\tau|$. Consequently,
\bae
 &\frac{|\b\tau-\b\eta|}{|\b\xi-\b\tau|}
 \left(\frac{\delta(\b x)^p+|\b\tau|^p+|\b\eta|^p}{\delta(\b x)^p+|\b\xi|^p+|\b\tau|^p}\right)^{\frac{p-2}{p}}
 \le 4^{p-2}K^{-1}.
 \label{eq:three.point.case.large}
\eae
If $1<p<2$, the triangle inequality instead yields
\bae
 \delta(\b x)^p+|\b\xi|^p+|\b\tau|^p
 &\le 2^p\big(\delta(\b x)^p+|\b\tau|^p+|\b\eta|^p\big)
 +2^{p-1}|\b\xi-\b\tau|^p\\
 &\le 2^{2p}\left(\frac{|\b\xi-\b\tau|}{|\b\tau-\b\eta|}\right)^p
 \big(\delta(\b x)^p+|\b\tau|^p+|\b\eta|^p\big),
\eae
where $|\b\tau-\b\eta|^p\le2^{p-1}(|\b\tau|^p+|\b\eta|^p)$. Reversing the ratio because $p-2<0$, we obtain
\bae
 &\frac{|\b\tau-\b\eta|}{|\b\xi-\b\tau|}
 \left(\frac{\delta(\b x)^p+|\b\tau|^p+|\b\eta|^p}{\delta(\b x)^p+|\b\xi|^p+|\b\tau|^p}\right)^{\frac{p-2}{p}}\le4^{2-p}\left(\frac{|\b\tau-\b\eta|}{|\b\xi-\b\tau|}\right)^{p-1}
 \le4K^{-(p-1)}.\label{eq:three.point.case.small}
\eae
Using $p\in[p_\petm,p_\petp]$, both \eqref{eq:three.point.case.large} and \eqref{eq:three.point.case.small} give
\bae
 ~&\sigma_\cst{c}\big(\delta(\b x)^p+|\b\tau|^p+|\b\eta|^p\big)^{\frac{p-2}{p}}
 |\b\tau-\b\eta|\,|\b\xi-\b\tau|\\
 &\quad\le\sigma_\cst{c}4^{\max(1,p_\petp-2)}K^{-\min(1,p_\petm-1)}
 \big(\delta(\b x)^p+|\b\xi|^p+|\b\tau|^p\big)^{\frac{p-2}{p}}|\b\xi-\b\tau|^2\\
 &\quad\le\varepsilon\big(\stress(\b x,\b\xi)-\stress(\b x,\b\tau)\big)\cdot(\b\xi-\b\tau),
\eae
where the last step follows from the definition of $K$ and \eqref{eq:three.point.mono}.
It remains to consider $|\b\xi-\b\tau|\le K|\b\tau-\b\eta|$. By \eqref{eq:three.point.cont},
\bae
 ~&|\stress(\b x,\b\tau)-\stress(\b x,\b\eta)|\,|\b\xi-\b\tau| \le\sigma_\cst{c}K
 \big(\delta(\b x)^p+|\b\tau|^p+|\b\eta|^p\big)^{\frac{p-2}{p}}
 |\b\tau-\b\eta|^2.\label{eq:three.point.remaining}
\eae
For $p\ge2$, the final factor already has the claimed form. For $1<p<2$, the triangle inequality and $|\b\tau-\b\eta|^p\le 2^{p-1}(|\b\tau|^p+|\b\eta|^p)$ yield
\bae
 \big(\delta(\b x)^p+|\b\tau|^p+|\b\eta|^p\big)^{\frac{p-2}{p}}|\b\tau-\b\eta|^2 &=|\b\tau-\b\eta|^p
 \left(\frac{|\b\tau-\b\eta|^p}{\delta(\b x)^p+|\b\tau|^p+|\b\eta|^p}\right)^{\frac{2-p}{p}}\\
 &\le 2^{\frac{(p-1)(2-p)}{p}}|\b\tau-\b\eta|^p
 \le 2^{1/4}|\b\tau-\b\eta|^p.
\eae
Combining the two cases proves \eqref{eq:three.point.bulk} with \eqref{eq:three.point.constant}. The same argument applies to the scalar stabilization flux, using \eqref{eq:S:holder.continuity} and \eqref{eq:S:strong.monotonicity} with their effective uniform constants.
}
\ep

\bl[Cell and face perturbation bounds]\label{lem:perturbation.remainder}
Let $T\in\T_h$, let $v:T\to[0,\infty)$ be measurable with $v\in L^1(T)$, and let {$w:T\to\RR^m$ ($m\ge1$)} be measurable with $w\in L^{p(\cdot)}(T,\RR^m)$ and $\|w\|_{L^{p(\cdot)}(T,\RR^m)}\le1$. Then
\bae\label{eq:perturbation.remainder}
 \int_T v^{1-\frac{\sobmin}p}|w|^{\sobmin}\,\cst d\b x
 \lesssim (1+\|v\|_{L^1(T)})\,
 \|w\|_{L^{p(\cdot)}(T,\RR^m)}^{\sobmin_{T,-}}.
\eae
The same estimate holds globally on $\Omega$: for measurable $v\ge0$ in $L^1(\Omega)$ and $w\in L^{p(\cdot)}(\Omega,\RR^m)$ with $\|w\|_{L^{p(\cdot)}(\Omega,\RR^m)}\le1$,
\bae\label{eq:perturbation.remainder.global}
 \int_\Omega v^{1-\frac{\sobmin}p}|w|^{\sobmin}\,\cst d\b x
 \lesssim (1+\|v\|_{L^1(\Omega)})\,
 \|w\|_{L^{p(\cdot)}(\Omega,\RR^m)}^{\sobmin_\petm}.
\eae
This follows from the same proof, using the global Luxemburg norm rather than summing estimates with distinct cellwise normalization factors.
The corresponding bound on the entire weighted broken skeleton is as follows: if $v\ge0$ is measurable with $\sum_{T\in\T_h}h_T\int_{\partial T}v\,\cst ds<\infty$ and $w$ is measurable with $\|w\|_{p(\cdot),\partial h}\le1$, then
\bae\label{eq:perturbation.remainder.face}
 \sum_{T\in\T_h}h_T\int_{\partial T}
 v^{1-\frac{\sobmin}p}|w|^{\sobmin}\,\cst ds
 \lesssim \left(1+\sum_{T\in\T_h}h_T\int_{\partial T}v\,\cst ds\right)
 \|w\|_{p(\cdot),\partial h}^{\sobmin_\petm}.
\eae
\el

\bp
{
If $\|w\|_{L^{p(\cdot)}(T,\RR^m)}=0$, the first assertion is immediate. Otherwise, set $\widehat T\eq T\cap\{\b x\in T:p(\b x)\ge2\}$ and $\widecheck T\eq T\cap\{\b x\in T:p(\b x)<2\}$. The modular--norm inequality~\eqref{eq:modular-norm-equiv} yields
\bae
 \int_{\widecheck T}|w|^p\,\cst d\b x
 \le \|w\|_{L^{p(\cdot)}(T,\RR^m)}^{p_{T,-}}
 \le \|w\|_{L^{p(\cdot)}(T,\RR^m)}^{\sobmin_{T,-}}.
\eae
On $\widehat T$, the scalar Young inequality~\eqref{eq:gener:young:scal} and the Luxemburg normalization give
\bae
 \int_{\widehat T}v^{1-2/p}|w|^2\,\cst d\b x &=\|w\|_{L^{p(\cdot)}(T,\RR^m)}^2
 \int_{\widehat T}v^{1-2/p}
 \left|\frac{w}{\|w\|_{L^{p(\cdot)}(T,\RR^m)}}\right|^2\,\cst d\b x\\
 &\le\|w\|_{L^{p(\cdot)}(T,\RR^m)}^2
 \int_{\widehat T}\left(\frac{p-2}{p}v+
 \frac2p\left|\frac{w}{\|w\|_{L^{p(\cdot)}(T,\RR^m)}}\right|^p\right)\,\cst d\b x\\
 &\le (1+\|v\|_{L^1(T)})\|w\|_{L^{p(\cdot)}(T,\RR^m)}^{\sobmin_{T,-}}.
\eae
Adding these bounds proves \eqref{eq:perturbation.remainder} (absorbing the factor $2$).
{The proof of \eqref{eq:perturbation.remainder.global} is identical, normalizing by $\|w\|_{L^{p(\cdot)}(\Omega,\RR^m)}$ and using $p_-$ in place of $p_{T,-}$.} For \eqref{eq:perturbation.remainder.face}, use the same two inequalities on $\bigcup_{T\in\T_h}\partial T$, sum with weights $h_T$, and normalize directly by the \emph{global} norm $\|w\|_{p(\cdot),\partial h}$ defined in \eqref{eq:error:skeleton:norm}. The modular of the normalized function is at most $1$, giving exactly \eqref{eq:perturbation.remainder.face}.}
\ep

\bl[Properties of the local aggregation]\label{lem:error:Ah-properties}
{ For any exponent $q$ satisfying the hypotheses preceding~\eqref{eq:error:Ah},} the functional $\mathcal A_{q,h}$ is monotone and positively homogeneous:
\bae\label{eq:error:Ah:homogeneity}
 \mathcal A_{q,h}((\lambda r_T)_{T\in\T_h})
 =\lambda\mathcal A_{q,h}({(r_T)_{T\in\T_h}})
 \qquad\forall\lambda\ge0.
\eae
Moreover, if $z_T\in L^{q(\cdot)}(T)$ satisfies $\|z_T\|_{L^{q(\cdot)}(T)}\lesssim r_T$ for every $T$, then the broken function $z_h$ obtained by patching the $z_T$ satisfies
\bae\label{eq:error:Ah:patching}
 \|z_h\|_{L^{q(\cdot)}(\Omega)}
 \lesssim \mathcal A_{q,h}({(r_T)_{T\in\T_h}}).
\eae
The same statement holds on the broken skeleton with the weighted measure $h_T\,\cst ds$.
\el

\bp
Let ${(r_T)_{T\in\T_h}}$ be non-negative. The homogeneity follows directly from the definition. For the patching property, take $\lambda>\mathcal A_{q,h}({(r_T)_{T\in\T_h}})$. Then $r_T/\lambda\le1$ for every $T$, and the local modular--norm relation gives
\bae
 \rho_{L^{q(\cdot)}(T)}\!\left(\frac{z_T}{C\lambda}\right)
 \le {\mathfrak m_{q,T}}\!\left(\frac{r_T}{\lambda}\right).
\eae
Summing over $T$ gives \eqref{eq:error:Ah:patching}. The skeleton version is identical after replacing the underlying measure by $h_T\,\cst ds$.
\ep

\bl[H\"older inequality for local aggregations]\label{lem:error:Ah-holder}
Let ${(r_T)_{T\in\T_h}}$ be a non-negative family and let $\u v_h\in\dU{h}{k}$. Then
\bae\label{eq:error:Ah-holder}
 \sum_{T\in\T_h} r_T\|\u v_T\|_{p(\cdot),T}
 \lesssim
 \mathcal A_{p',h}({(r_T)_{T\in\T_h}})\,\|\u v_h\|_{p(\cdot),h}.
\eae
\el

\bp
By homogeneity, it is enough to consider the normalization
$\mathcal A_{p',h}({(r_T)_{T\in\T_h}})=\|\u v_h\|_{p(\cdot),h}=1$.
If $\|\u v_T\|_{p(\cdot),T}\le1$, the local
modular--norm inequalities~\eqref{eq:modular-norm-equiv} give
$\|\u v_T\|_{p(\cdot),T}^{p_{T,+}}\lesssim\rho_{p(\cdot),T}(\u v_T)$; Young's inequality~\eqref{eq:gener:young:scal}
with the conjugate pair $(p_{T,+},p'_{T,+})$ therefore yields
\bae
 r_T\|\u v_T\|_{p(\cdot),T}\lesssim r_T^{p'_{T,+}}+\rho_{p(\cdot),T}(\u v_T).
\eae
If $\|\u v_T\|_{p(\cdot),T}>1$, the same argument uses the pair
$(p_{T,-},p'_{T,-})$. In both cases,
\bae
 r_T\|\u v_T\|_{p(\cdot),T}\lesssim
 {\mathfrak m_{p',T}}(r_T)+\rho_{p(\cdot),T}(\u v_T).
\eae
Summing over the mesh and using the two normalizations proves
\eqref{eq:error:Ah-holder}.
\ep

\bp[Proof of Theorem~\ref{thm:error.estimate}]
Set $\u e_h\eq\u u_h-\Interp{h}{k}u$,
 $\u{\hat u}_h\eq\Interp{h}{k}u$,
 and $\u{\check u}_h\eq\Interp{h}{k}(\proj{h}{k+1}u)$,
and define the consistency functional
\bae
 \mathcal E_h(\u v_h)
 \eq-(\DIV\stress(\cdot,\GRAD u),v_h)_\Omega-\cst a_h(\u{\hat u}_h,\u v_h).
\eae
Testing with $\u e_h$ and using the discrete equation gives
\bae
 \mathcal E_h(\u e_h)
 =\cst a_h(\u u_h,\u e_h)-\cst a_h(\u{\hat u}_h,\u e_h).
\eae
We next turn the strong-monotonicity energy into a bound for the discrete norm. By strong monotonicity,
\bae
 \mathcal E_h(\u e_h)&\gtrsim \sigma_\cst{m}\Biggl[
 \int_\Omega
 \big(\delta^{p(\cdot)}+|\dgrad{k}{h}\u u_h|^{p(\cdot)}+|\dgrad{k}{h}\u{\hat u}_h|^{p(\cdot)}\big)^{\frac{p(\cdot)-2}{p(\cdot)}}
 |\dgrad{k}{h}\u e_h|^2\,\cst{d}\b x\\
 &\qquad+\sum_{T\in\T_h}h_T\int_{\partial T}
 \big(\delta_{\partial T}^{p(\cdot)}+|\dfbres{k}{\partial T}\u u_T|^{p(\cdot)}+|\dfbres{k}{\partial T}\u{\hat u}_T|^{p(\cdot)}\big)^{\frac{p(\cdot)-2}{p(\cdot)}}
 |\dfbres{k}{\partial T}\u e_T|^2\,\cst{d}s\Biggr].
\eae
By Proposition~\ref{pr:discrete_norms_equiv} and the boundedness of $\Interp{h}{k}$, $ \|\u{\hat u}_h\|_{p(\cdot),\cst{G},h}
 \lesssim \|\GRAD u\|_{\tL^{p(\cdot)}(\Omega)}$, and hence
\bae
 ~&\int_\Omega (\delta^{p(\cdot)}+|\dgrad{k}{h}\u u_h|^{p(\cdot)}+|\dgrad{k}{h}\u{\hat u}_h|^{p(\cdot)})\,\cst{d}\b x
 +\sum_{T\in\T_h}h_T\int_{\partial T}(\delta_{\partial T}^{p(\cdot)}+|\dfbres{k}{\partial T}\u u_T|^{p(\cdot)}+|\dfbres{k}{\partial T}\u{\hat u}_T|^{p(\cdot)})\,\cst{d}s\\
 &\qquad\lesssim
 \mathcal D_h
 +\rho_{p(\cdot),\cst{G},h}(\u u_h)
 +\rho_{p(\cdot),\cst{G},h}(\u{\hat u}_h)\\
 &\qquad\lesssim
 \max\!\left\{1,
 \rho_{L^{p'(\cdot)}(\Omega)}(f)^{p'_\petp/p'_\petm},
 \rho_{L^{p'(\cdot)}(\Omega)}(f)^{p'_\petm/p'_\petp},
 \mathcal D_h
 \right\},
\eae
where we concluded with the a priori bounds of Theorems~\ref{prop:a-priori} and~\ref{thm:well-posedness}. In particular, the resulting constant may depend on $\sigma_\cst{m}$ through these well-posedness estimates; this dependence is not included in the mesh-uniform factors denoted by $\lesssim$.
Now, assume first that $\mathcal E_h(\u e_h)>0$. On the set $\{p<2\}$, the pointwise identity
\bae
 \left(\frac{|\dgrad{k}{h}\u e_h|}
 {\mathcal E_h(\u e_h)^{1/2}}\right)^p
 &= \left(
 \frac{\bigl(\delta^{p(\cdot)}+|\dgrad{k}{h}\u u_h|^{p(\cdot)}
 +|\dgrad{k}{h}\u{\hat u}_h|^{p(\cdot)}\bigr)^{(p-2)/p}
 |\dgrad{k}{h}\u e_h|^2}{\mathcal E_h(\u e_h)}
 \right)^{p/2}\\
 &\qquad\times
 \bigl(\delta^{p(\cdot)}+|\dgrad{k}{h}\u u_h|^{p(\cdot)}
 +|\dgrad{k}{h}\u{\hat u}_h|^{p(\cdot)}\bigr)^{(2-p)/2}
\eae
and the weighted arithmetic--geometric mean inequality yield
\bae
 \rho_{L^{p(\cdot)}(\{p<2\})}\!\left(
 \frac{\dgrad{k}{h}\u e_h}{\mathcal E_h(\u e_h)^{1/2}}\right)
 &\lesssim
 \frac{1}{\mathcal E_h(\u e_h)}\int_{\{p<2\}}
 \bigl(\delta^{p(\cdot)}+|\dgrad{k}{h}\u u_h|^{p(\cdot)}
 +|\dgrad{k}{h}\u{\hat u}_h|^{p(\cdot)}\bigr)^{(p-2)/p}
 |\dgrad{k}{h}\u e_h|^2\,\cst{d}\b x\\
 &\qquad+\int_{\{p<2\}}
 \bigl(\delta^{p(\cdot)}+|\dgrad{k}{h}\u u_h|^{p(\cdot)}
 +|\dgrad{k}{h}\u{\hat u}_h|^{p(\cdot)}\bigr)\,\cst{d}\b x \lesssim 1.
\eae
Hence
$\|\dgrad{k}{h}\u e_h\|_{L^{p(\cdot)}(\{p<2\})}\lesssim\mathcal E_h(\u e_h)^{1/2}$.
On $\{p\ge2\}$, since $|\dgrad{k}{h}\u e_h|^p\lesssim (\delta^{p(\cdot)}+|\dgrad{k}{h}\u u_h|^{p(\cdot)}+|\dgrad{k}{h}\u{\hat u}_h|^{p(\cdot)})$,
\bae
 (\delta^{p(\cdot)}+|\dgrad{k}{h}\u u_h|^{p(\cdot)}+|\dgrad{k}{h}\u{\hat u}_h|^{p(\cdot)})^{(p-2)/p}|\dgrad{k}{h}\u e_h|^2\gtrsim |\dgrad{k}{h}\u e_h|^p,
\eae
so that the modular--norm inequalities~\eqref{eq:modular-norm-equiv} give
$\|\dgrad{k}{h}\u e_h\|_{L^{p(\cdot)}(\{p\ge2\})}
 \lesssim\mathcal E_h(\u e_h)^{1/p_\petp}$ when $\mathcal E_h(\u e_h)\le1$.
The same two arguments apply to the skeleton with the measure $h_T\,\cst ds$. Therefore, using Proposition~\ref{pr:discrete_norms_equiv}, if $\mathcal E_h(\u e_h)\le1$,
\bae
 \|\u e_h\|_{p(\cdot),h}
 \lesssim
 \mathcal E_h(\u e_h)^{1/2}
 +\mathcal E_h(\u e_h)^{1/p_\petp}
 \lesssim
 \mathcal E_h(\u e_h)^{1/\sobmax_\petp}.
\eae
If $\mathcal E_h(\u e_h)>1$, the same conclusion follows from the uniform bound on $\|\u e_h\|_{p(\cdot),h}$ given by the a priori bounds, Proposition~\ref{pr:interpolator-bounded}, and Proposition~\ref{pr:discrete_norms_equiv}. Thus
\bae\label{eq:error.lower.norm}
 \|\u e_h\|_{p(\cdot),h}
 \le C\,\mathcal E_h(\u e_h)^{1/\sobmax_\petp}.
\eae
Owing to this uniform bound, if $\mathcal A_{p',h}((h_T^{k+1}|\stress(\cdot,\GRAD u)|_{\tW^{k+1,p'(\cdot)}(T)})_{T\in\T_h})+\mathcal A_{p,h}((h_T^{k+1}|u|_{W^{k+2,p(\cdot)}(T)})_{T\in\T_h})$ is larger than a fixed structural
constant, \eqref{eq:error.estimate.localized} is immediate after enlarging
$C$. We may therefore assume below that these two quantities are small
enough. By Proposition~\ref{pr:variable-projector-approx} and
Lemma~\ref{lem:error:Ah-properties}, all approximation norms to which
Lemma~\ref{lem:perturbation.remainder} is applied are then at most one.
This reduction only avoids carrying the two modular regimes through the
consistency estimates and does not impose an additional asymptotic
assumption on the theorem.

We now derive the consistency decomposition. Continuity of the normal trace of $\stress(\cdot,\GRAD u)$ and the single-valuedness of the face unknowns imply
\bae
 \sum_{T\in\T_h}(e_{\partial T},\stress(\cdot,\GRAD u)\cdot\b n_{\partial T})_{\partial T}=0.
\eae
Using this identity, an elementwise integration by parts, the definition of the $L^2$ projector, and \eqref{eq:G}, we obtain
\bae\label{eq:consistency:ah:EJ}
 \mathcal E_h(\u e_h)
 &=\sum_{T\in\T_h}\Big[(\stress(\cdot,\GRAD u),\GRAD e_T)_T
 +(e_{\partial T}-e_T,\stress(\cdot,\GRAD u)\cdot\b n_{\partial T})_{\partial T}\Big]
 -\cst a_h(\u{\hat u}_h,\u e_h)\\
 &=\sum_{T\in\T_h}\Big[(\PROJ{T}{k}\stress(\cdot,\GRAD u),\GRAD e_T)_T
 +(e_{\partial T}-e_T,\stress(\cdot,\GRAD u)\cdot\b n_{\partial T})_{\partial T}\Big]
 -\cst a_h(\u{\hat u}_h,\u e_h)\\
 &=\mathcal I_1+(\PROJ{h}{k}\stress(\cdot,\GRAD u),\dgrad{k}{h}\u e_h)_\Omega
 -\cst a_h(\u{\hat u}_h,\u e_h)\\
 &=\mathcal I_1+\mathcal I_2+\mathcal I_3,
\eae
where
\ban
 \mathcal I_1 &\eq \sum_{T\in\T_h}
 (e_{\partial T}-e_T,(\stress(\cdot,\GRAD u)-\PROJ{T}{k}\stress(\cdot,\GRAD u))\cdot\b n_{\partial T})_{\partial T},\\
 \mathcal I_2 &\eq (\stress(\cdot,\GRAD u)-\stress(\cdot,\PROJ{h}{k}\GRAD u),\dgrad{k}{h}\u e_h)_\Omega,\\
 \mathcal I_3 &\eq \cst{s}_h(\u{\check u}_h,\u e_h)-\cst{s}_h(\u{\hat u}_h,\u e_h).
\ean
Here we used $\dgrad{k}{h}\u{\hat u}_h=\PROJ{h}{k}\GRAD u$, the identity
$(\PROJ{h}{k}\stress(\cdot,\GRAD u),\dgrad{k}{h}\u e_h)_\Omega=(\stress(\cdot,\GRAD u),\dgrad{k}{h}\u e_h)_\Omega$, and $\cst{s}_h(\u{\check u}_h,\u e_h)=0$ by polynomial consistency.

For $\mathcal I_1$, H\"older's inequality~\eqref{eq:holder} on the weighted skeleton, \eqref{eq:error.lower.norm}, and Young's inequality~\eqref{eq:gener:young:scal} give
\bae\label{eq:error:I1-simple}
 |\mathcal I_1|
 \le \varepsilon\mathcal E_h(\u e_h)
 +C_\varepsilon
 \left\|\left((\PROJ{T}{k}\stress(\cdot,\GRAD u)-\stress(\cdot,\GRAD u))\res{\partial T}\right)_{T\in\T_h}\right\|_{p'(\cdot),\partial h}^{\frac{\sobmax_\petp}{\sobmax_\petp-1}}.
\eae
{
For $\mathcal I_2$, Lemma~\ref{lem:three.point.error} first gives the explicit perturbation remainder
\bae
 |\mathcal I_2|&\le\varepsilon\mathcal E_h(\u e_h)
 +C_\varepsilon\int_\Omega
 \big(\delta^{p(\cdot)}+|\GRAD u|^{p(\cdot)}
 +|\PROJ{h}{k}\GRAD u|^{p(\cdot)}\big)^{1-\frac{\sobmin}{p(\cdot)}}
 |\PROJ{h}{k}\GRAD u-\GRAD u|^{\sobmin}\,\cst d\b x.
\eae
{Applying the global volume estimate \eqref{eq:perturbation.remainder.global} directly (with $v=\delta^{p(\cdot)}+|\GRAD u|^{p(\cdot)}+|\PROJ{h}{k}\GRAD u|^{p(\cdot)}$ and $w=\PROJ{h}{k}\GRAD u-\GRAD u$) yields}
\bae\label{eq:error:I2-simple}
 |\mathcal I_2|
 \le \varepsilon\mathcal E_h(\u e_h)
 +C_\varepsilon
 \|\PROJ{h}{k}\GRAD u-\GRAD u\|_{\tL^{p(\cdot)}(\Omega)}^{\sobmin_\petm}.
\eae
Here the factor $1+\int_\Omega(\delta^{p(\cdot)}+|\GRAD u|^{p(\cdot)}+|\PROJ{h}{k}\GRAD u|^{p(\cdot)})\,\cst{d}\b x$ is absorbed into the uniform a priori constant.

For $\mathcal I_3$, the scalar version of Lemma~\ref{lem:three.point.error} gives the corresponding explicit remainder on the broken skeleton:
\bae
 |\mathcal I_3|&\le\varepsilon\mathcal E_h(\u e_h)
 +C_\varepsilon\sum_{T\in\T_h}h_T\int_{\partial T}
 \Big(\delta_{\partial T}^{p(\cdot)}
 +|\dfbres{k}{\partial T}\u{\hat u}_T|^{p(\cdot)}
 +|\dfbres{k}{\partial T}\u{\check u}_T|^{p(\cdot)}\Big)^{1-\frac{\sobmin}{p(\cdot)}}
 |\dfbres{k}{\partial T}(\u{\hat u}_T-\u{\check u}_T)|^{\sobmin}
 \,\cst ds.
\eae
Applying the face estimate \eqref{eq:perturbation.remainder.face} (with the sum weighted by $h_T$) gives
\bae\label{eq:error:I3-simple}
 |\mathcal I_3|
 \le \varepsilon\mathcal E_h(\u e_h)
 +C_\varepsilon
 \left\|\dfbres{k}{\partial h}(\u{\hat u}_h-\u{\check u}_h)\right\|_{p(\cdot),\partial h}^{\sobmin_\petm}.
\eae
The face-weight integral is controlled by the uniform stabilization bounds.}
After absorption and \eqref{eq:error.lower.norm},
\bae
 \|\u e_h\|_{p(\cdot),h}
 \le C\Bigg[&
 \left\|\left((\PROJ{T}{k}\stress(\cdot,\GRAD u)-\stress(\cdot,\GRAD u))\res{\partial T}\right)_{T\in\T_h}\right\|_{p'(\cdot),\partial h}^{\frac1{\sobmax_\petp-1}}\\
 &+\left(
 \|\PROJ{h}{k}\GRAD u-\GRAD u\|_{\tL^{p(\cdot)}(\Omega)}
 +\left\|\dfbres{k}{\partial h}(\u{\hat u}_h-\u{\check u}_h)\right\|_{p(\cdot),\partial h}
 \right)^{\frac{\sobmin_\petm}{\sobmax_\petp}}
 \Bigg].
\eae
The projector approximation and trace estimates proved above give directly
\bae
 \left\|\left((\PROJ{T}{k}\stress(\cdot,\GRAD u)-\stress(\cdot,\GRAD u))\res{\partial T}\right)_{T\in\T_h}\right\|_{p'(\cdot),\partial h}
 \lesssim
 \mathcal A_{p',h}\!\left(\left(h_T^{k+1}|\stress(\cdot,\GRAD u)|_{\tW^{k+1,p'(\cdot)}(T)}\right)_{T\in\T_h}\right),
\eae
and
\bae
 \|\PROJ{h}{k}\GRAD u-\GRAD u\|_{\tL^{p(\cdot)}(\Omega)}
 +\left\|\dfbres{k}{\partial h}(\u{\hat u}_h-\u{\check u}_h)\right\|_{p(\cdot),\partial h}
 \lesssim
 \mathcal A_{p,h}\!\left(\left(h_T^{k+1}|u|_{W^{k+2,p(\cdot)}(T)}\right)_{T\in\T_h}\right).
\eae
This proves \eqref{eq:error.estimate.localized}. To obtain the global
form, we verify the local hypotheses
\eqref{eq:local-to-global-cell}--\eqref{eq:local-to-global-face}
for each of the three approximation errors appearing above.
First, applying the trace estimate \eqref{eq:proj:app:norm:F}
with $q=p'$, $\ell=k$, and $n=k+1$ to
$\stress(\cdot,\GRAD u)$ yields, for every $T\in\T_h$,
\bae
 \big\|h_T^{1/p'(\cdot)}
  (\PROJ{T}{k}\stress(\cdot,\GRAD u)-\stress(\cdot,\GRAD u))
 \big\|_{\tL^{p'(\cdot)}(\partial T)}
 \lesssim h_T^{k+1}
 |\stress(\cdot,\GRAD u)|_{\tW^{k+1,p'(\cdot)}(T)}.
\eae
Second, the cell estimate \eqref{eq:proj:app:norm:T}, applied to
$\GRAD u$ with $q=p$, $\ell=k$, $n=k+1$, and $m=0$, gives
\bae
 \|\PROJ{T}{k}\GRAD u-\GRAD u\|_{\tL^{p(\cdot)}(T)}
 \lesssim h_T^{k+1}|u|_{W^{k+2,p(\cdot)}(T)}.
\eae
Finally, polynomial consistency \eqref{eq:sT:polynomial.consistency}
and the definitions of $\u{\hat u}_T$ and $\u{\check u}_T$ show that
\bae
 \dfbres{k}{\partial T}(\u{\hat u}_T-\u{\check u}_T)
 =\dfbres{k}{\partial T}\Interp{T}{k}
 (u-\proj{T}{k+1}u).
\eae
By Proposition~\ref{pr:discrete_norms_equiv}, the interpolator bound
\eqref{eq:I:boundedness:norm}, and the approximation estimate
\eqref{eq:proj:app:norm:T} with $\ell=k+1$, $n=k+2$, and $m=1$,
\bae
 \big\|h_T^{1/p(\cdot)}
 \dfbres{k}{\partial T}(\u{\hat u}_T-\u{\check u}_T)
 \big\|_{L^{p(\cdot)}(\partial T)}
 \lesssim h_T^{k+1}|u|_{W^{k+2,p(\cdot)}(T)}.
\eae
Lemma~\ref{lem:local-to-global-projector} now applies to the flux
face error with $(q,n,\beta)=(p',k+1,k+1)$, and to the solution
cell and face errors with $(q,n,\alpha,\beta)=(p,k+2,k+1,k+1)$.
Using \eqref{eq:local-to-global-cell-bound}--\eqref{eq:local-to-global-face-bound}
then gives

\bae\label{eq:error:global-approx}
 \left\|\left((\PROJ{T}{k}\stress(\cdot,\GRAD u)-\stress(\cdot,\GRAD u))\res{\partial T}\right)_{T\in\T_h}\right\|_{p'(\cdot),\partial h}
 &\lesssim h^{k+1}|\stress(\cdot,\GRAD u)|_{\tW^{k+1,p'(\cdot)}(\T_h)},\\
 \|\PROJ{h}{k}\GRAD u-\GRAD u\|_{\tL^{p(\cdot)}(\Omega)}
 +\left\|\dfbres{k}{\partial h}(\u{\hat u}_h-\u{\check u}_h)\right\|_{p(\cdot),\partial h}
 &\lesssim h^{k+1}|u|_{W^{k+2,p(\cdot)}(\T_h)}.
\eae
Inserting \eqref{eq:error:global-approx} in the consistency estimate above proves \eqref{eq:error.estimate.global}. The factor $h^{k+1}$ is therefore extracted at the level of the Luxemburg norms themselves, before applying the powers $1/(\sobmax_\petp-1)$ and $\sobmin_\petm/\sobmax_\petp$. 
\ep

\bp[Proof of Theorem~\ref{thm:error.estimate:regime}]

Since $p_\petp\le2$, \eqref{eq:error.lower.norm} becomes
\bae\label{eq:regime:quadratic-lower}
 \|\u e_h\|_{p(\cdot),h}^2\le C\,\mathcal E_h(\u e_h).
\eae
The flux-projection term is estimated exactly as in the proof of
Theorem~\ref{thm:error.estimate}, before Young's inequality~\eqref{eq:gener:young:scal}:
\bae\label{eq:regime:I1}
 |\mathcal I_1|
 \lesssim
 \mathcal A_{p',h}\!\left(\left(h_T^{k+1}|\stress(\cdot,\GRAD u)|_{\tW^{k+1,p'(\cdot)}(T)}\right)_{T\in\T_h}\right)
 \|\u e_h\|_{p(\cdot),h}.
\eae
We refine the other two terms element by element. By Proposition~\ref{pr:variable-projector-approx}, the projection error satisfies $\|\PROJ{T}{k}\GRAD u-\GRAD u\|_{\tL^{p(\cdot)}(T)}\lesssim h_T^{k+1}|u|_{W^{k+2,p(\cdot)}(T)}$.
For the degenerate branch, the continuity assumption and $p(\cdot)\le2$
give
$|\stress(\cdot,\PROJ{T}{k}\GRAD u)-\stress(\cdot,\GRAD u)|
\lesssim |\PROJ{T}{k}\GRAD u-\GRAD u|^{p(\cdot)-1}$. When the projection error is nonzero, normalization by its $L^{p(\cdot)}(T)$ norm yields
\bae
 \big\||\PROJ{T}{k}\GRAD u-\GRAD u|^{p(\cdot)-1}\big\|_{\tL^{p'(\cdot)}(T)}
 \lesssim {\mathfrak m_{p-1,T}}(\|\PROJ{T}{k}\GRAD u-\GRAD u\|_{\tL^{p(\cdot)}(T)})
 \lesssim {\mathfrak m_{p-1,T}}(h_T^{k+1}|u|_{W^{k+2,p(\cdot)}(T)}),
\eae
where the harmless constant in the last bound is absorbed using the monotonicity and {uniform scaling estimate
${\mathfrak m_{q,T}}(Ct)\le\max(C^{q_{T,-}},C^{q_{T,+}}){\mathfrak m_{q,T}}(t)$ for $C\ge1$, which follows directly from~\eqref{eq:error:mqT}.}
We now write $\mathcal I_{2,T}$ for the contribution of element $T$ to the sum defining $\mathcal I_2$. Then,
\bae\label{eq:regime:I2-degenerate}
 |\mathcal I_{2,T}|
 \lesssim {\mathfrak m_{p-1,T}}(h_T^{k+1}|u|_{W^{k+2,p(\cdot)}(T)})\,\|\u e_T\|_{p(\cdot),T}.
\eae

Assume now $\mathfrak D_T>0$. Since the exponent $p(\b x)-2$ is
non-positive and $\delta+|\GRAD u|\ge\mathfrak D_T$, the non-degenerate
continuity bound gives
\bae
 |\stress(\b x,\PROJ{T}{k}\GRAD u)-\stress(\b x,\GRAD u)|
 \lesssim \mathfrak D_T^{p_{T,-}-2}|\PROJ{T}{k}\GRAD u-\GRAD u|.
\eae
Here the change of exponent is performed on the \emph{discrete error}, not
on $\PROJ{T}{k}\GRAD u-\GRAD u$. Variable-exponent H\"older{~\eqref{eq:holder}} and the HHO discrete Lebesgue embedding,
Lemma~\ref{lem:error:local-change-exp} with $q=p'$, give
\bae
 |\mathcal I_{2,T}|
 &\lesssim
 \mathfrak D_T^{p_{T,-}-2}
 \|\PROJ{T}{k}\GRAD u-\GRAD u\|_{\tL^{p(\cdot)}(T)}
 \|\u e_T\|_{p'(\cdot),T}\notag\\
 &\lesssim 
 h_T^{d(1/(p')_T-1/p_T)}
 \mathfrak D_T^{p_{T,-}-2}h_T^{k+1}|u|_{W^{k+2,p(\cdot)}(T)}
 \|\u e_T\|_{p(\cdot),T}.
 \label{eq:regime:I2-nondegenerate}
\eae
If $0<\eta_T<1$, let $\b x_T\in\overline T$ be such that
$p(\b x_T)=p_{T,-}$ (recall that $p$ is continuous). Applying
\eqref{eq:mean.exp.diam.equiv} to $p$ and $p'$ gives
\bae
 h_T^{d(1/(p')_T-1/p_T)}
 \eqsim h_T^{d(1/p'(\b x_T)-1/p(\b x_T))}
 =h_T^{d(1-2/p_{T,-})},
 \qquad
 |T|^{1/p_T}\eqsim h_T^{d/p_{T,-}}.
\eae
Using these equivalences in the definition of $\eta_T$ yields
\bae\label{eq:regime:eta-local}
 h_T^{d(1/(p')_T-1/p_T)}
 \mathfrak D_T^{p_{T,-}-2}h_T^{k+1}|u|_{W^{k+2,p(\cdot)}(T)}&\lesssim
 \eta_T^{2-p_{T,-}}\big(h_T^{k+1}|u|_{W^{k+2,p(\cdot)}(T)}\big)^{p_{T,-}-1}\\
 &\le \eta_T^{2-p_{T,-}}
 {\mathfrak m_{p-1,T}}\!\left(h_T^{k+1}|u|_{W^{k+2,p(\cdot)}(T)}\right).
\eae
For $\eta_T\ge1$ (and also when $\mathfrak D_T=0$), we use
\eqref{eq:regime:I2-degenerate}. {Thus, in all cases,}
\bae\label{eq:regime:I2-local}
 |\mathcal I_{2,T}|
 &\lesssim \min(1,\eta_T)^{2-p_{T,-}}
 {\mathfrak m_{p-1,T}}\!\left(h_T^{k+1}|u|_{W^{k+2,p(\cdot)}(T)}\right)\|\u e_T\|_{p(\cdot),T}.
\eae
Summing and applying Lemma~\ref{lem:error:Ah-holder} gives
\bae\label{eq:regime:I2}
 |\mathcal I_2|
 &\lesssim \mathcal A_{p',h}\!\left(\left(
 \min(1,\eta_T)^{2-p_{T,-}}
 {\mathfrak m_{p-1,T}}\!\left(h_T^{k+1}|u|_{W^{k+2,p(\cdot)}(T)}\right)
 \right)_{T\in\T_h}\right)\|\u e_h\|_{p(\cdot),h}.
\eae

The stabilization term is treated with the weighted boundary measure. The stabilization difference is $\u{\hat u}_T-\u{\check u}_T=\Interp{T}{k}(u-\proj{T}{k+1}u)$. Polynomial consistency gives
$\dfbres{k}{\partial T}\u{\check u}_T=0$, while
Propositions~\ref{pr:interpolator-bounded}, \ref{pr:discrete_norms_equiv},
and \ref{pr:variable-projector-approx} yield
\bae
 \|(\u{\hat u}_T-\u{\check u}_T)\|_{p(\cdot),T}\lesssim h_T^{k+1}|u|_{W^{k+2,p(\cdot)}(T)}.
\eae
The degenerate continuity estimate for $\stabfl_F$ gives the analogue of
\eqref{eq:regime:I2-degenerate}. In the non-degenerate branch we put the
stabilization difference in the weighted $L^{p(\cdot)}(\partial T)$ norm
and the residual of $\u e_T$ in the weighted
$L^{p'(\cdot)}(\partial T)$ norm; the HHO discrete Lebesgue embedding, Lemma~\ref{lem:error:local-change-exp},
is again applied to the discrete unknown $\u e_T$. Since
$\delta_F\ge\mathfrak D_T$, the same calculation as
\eqref{eq:regime:I2-nondegenerate}--\eqref{eq:regime:eta-local} gives
{
\bae\label{eq:regime:I3}
 |\mathcal I_3|
 &\lesssim \mathcal A_{p',h}\!\left(\left(
 \min(1,\eta_T)^{2-p_{T,-}}
 {\mathfrak m_{p-1,T}}\!\left(h_T^{k+1}|u|_{W^{k+2,p(\cdot)}(T)}\right)
 \right)_{T\in\T_h}\right)\|\u e_h\|_{p(\cdot),h}.
\eae
}

Combining \eqref{eq:regime:quadratic-lower}, \eqref{eq:regime:I1},
\eqref{eq:regime:I2}, and \eqref{eq:regime:I3}, and cancelling one factor
$\|\u e_h\|_{p(\cdot),h}$ (the zero case being immediate), proves
\eqref{eq:error.estimate.regime}.

\ep

\bp[Proof of Corollary~\ref{c:error.estimate:grad}]
Use the triangle inequality with  Propositions~\ref{pr:discrete_norms_equiv} and \ref{pr:variable-projector-approx}. The projection term is of order $h^{\widetilde k}$ in the broken $W^{1,p(\cdot)}$ seminorm.
\ep

\section{Numerical results}\label{sec:numerics-a-priori}

The aim of this section is to assess numerically the error estimate of
Theorem~\ref{thm:error.estimate}. We compare the convergence rates observed
in practice with the order predicted by the theorem.

\subsection{General setting}
\label{sec:num.res.setting}

\paragraph{Computational domain and meshes.} The
computations are carried out on the unit square $\Omega=(0,1)^2$, with the
Carreau--Yasuda flux of Example~\ref{ex:Carreau--Yasuda} with $\mu\equiv1$,
$\alpha\equiv2$, and $\delta\equiv0$ or $\delta\equiv1$; the value $\delta=0$ yields the
degenerate $p(\cdot)$-Laplacian, while $\delta=1$ removes the degeneracy.
The stabilization relies on a flux of the same type,
$\stabfl_F(\b x,w)=\mu_F\big(\delta^2+|w|^2\big)^{(p(\b x)-2)/2}w$ with
$\mu_F=2(k+1)(k+2)$ and $\delta_F=\delta$. We use the scheme of Section~\ref{sec:scheme} in its
equal-order version $\widetilde k=k$, with $k\in\{0,1,2,3\}$. Except in
Section~\ref{sec:num.res.voronoi}, which is devoted to polygonal meshes,
$\Omega$ is partitioned into uniform Cartesian grids of square cells. We set
$h_T=\mathrm{diam}(T)/\sqrt2$, which is the side length of $T$ when $T$ is a
square. For every cell, including Voronoi cells, $h_T$ equals the geometric diameter up to the constant factor $1/\sqrt2$, so that this choice only rescales the stabilization and the error measure by a uniform constant. The mesh size displayed in the figures is 
$h=\max_{T\in\T_h}\mathrm{diam}(T)$.

\paragraph{Solver.} The discrete nonlinear problem is
solved by a damped Newton method, in which both the flux and the
stabilization are linearized exactly, and the step length is halved until the
residual decreases. This iteration is embedded in a continuation on the
exponent: the problem is first solved for the constant exponent $2$, which
makes it linear, and then for the fields $p_\theta=2+\theta\,(p-2)$, with
$\theta$ increasing from $0$ to $1$. Each value of $\theta$ starts from the
solution computed at the previous one, and the increment in $\theta$ is
halved whenever the iteration stagnates. The Newton method stops once the
relative increment is below $10^{-8}$ and the residual is below $10^{-9}$
times the norm of the right-hand side. For $\delta=0$ and exponents taking values below
$2$, the derivative of the stabilization blows up wherever the face
residuals vanish, and the Newton method no longer converges on the finer
meshes, mostly for $k\geq1$, as already observed
in~\cite{harnist-romain-conv}. On these meshes, we use instead the Picard iteration,
which evaluates the factors $(\delta^2+|\cdot|^2)^{(p-2)/2}$ of the flux
and of the stabilization at the previous iterate and therefore
never involves this derivative. Each Picard update is under-relaxed by
the factor $\min(p_\petm-1,1/(p_\petp-1))^{1/2}$, and the iteration stops once
the relative increment of the face unknowns is below $10^{-10}$; on a
few of the finest meshes with $k\geq2$, where this tolerance is not reached,
it is relaxed to $10^{-9}$.

\paragraph{Exact solutions.} In each test, the source
term $f$ and the Dirichlet data are computed from a prescribed exact solution,
chosen among
\bse\label{eq:num.sol}
\ba
  u(x,y) &= e^{x+y}, \label{eq:num.sol.exp}\\
  u(x,y) &= x^5+y^5+x+y, \label{eq:num.sol.poly}\\
  u(x,y) &= \sin(\pi x)\sin(\pi y), \label{eq:num.sol.sin}\\
  u(x,y) &= \big((x-\tfrac12)^2+(y-\tfrac12)^2\big)^{\gamma/2},
  \qquad \gamma\in\{2,4\}. \label{eq:num.sol.rad}
\ea
\ese
Nonzero boundary traces are imposed through the boundary face unknowns; see Remark~\ref{rem:nonhomogeneous}. The gradients of the exponential and polynomial solutions never vanish: the
flux $\stress(\cdot,\GRAD u)$ is then smooth, and the regularity assumptions
of Theorem~\ref{thm:error.estimate} hold for every $k$. In contrast, the
gradient of the sinusoidal solution vanishes at the center
$\b x_c=(\frac12,\frac12)$ and at the four corners of $\Omega$, and that of
the radial solutions at $\b x_c$ only. At these critical points, the flux
remains smooth for $\delta=1$, but its regularity is limited for $\delta=0$,
as analyzed in Section~\ref{sec:num.res.sharp}.

\paragraph{Exponent fields.} We consider three exponent
fields with values in $[1.5,2.5]$,
\bse\label{eq:num.law}
\ba
  p(x,y) &= 1.5+0.5\sin(\pi x)\sin(\pi y), \label{eq:num.law.low}\\
  p(x,y) &= 2+0.5\sin(\pi x)\sin(\pi y), \label{eq:num.law.high}\\
  p(x,y) &= 1.5+\sin(\pi x)\sin(\pi y), \label{eq:num.law.mid}
\ea
\ese
which range over $[1.5,2]$, $[2,2.5]$ and $[1.5,2.5]$, respectively. We refer
to them as the fields corresponding to the subquadratic, superquadratic, and mixed regimes, respectively. Each of them
reaches its minimum $p_\petm$ on $\partial\Omega$ and its maximum $p_\petp$ at
the center of the domain.

\paragraph{Error measure and predicted order.} For
$\u v_h\in\dU{h}{k}$, we set
$\rho_{\mathrm{e},h}(\u v_h)=\sum_{T\in\T_h}\rho_{\mathrm{e},T}(\u v_T)$, where
\bae
  \rho_{\mathrm{e},T}(\u v_T) \eq
    \rho_{\tL^{p(\cdot)}(T)}(\dgrad{k}{T}\u v_T)
    + \mu_F\,h_T\,\rho_{L^{p(\cdot)}(\partial T)}(\dfbres{k}{\partial T}\u v_T)
\eae
combines the reconstructed gradient and the stabilization term, and we
denote by
$\|\u v_h\|_{\mathrm{e},h}\eq\inf\{\lambda>0:\rho_{\mathrm{e},h}(\u v_h/\lambda)\leq1\}$
the associated Luxemburg norm. The figures display the discrete error
$\|\u u_h-\Interp{h}{k}u\|_{\mathrm{e},h}$, written $\|u-u_h\|_{\mathrm e,h}$
for brevity. The integrals are approximated by a quadrature rule that uses
the value of $p$ at every quadrature point, and the infimum defining the norm
is located by bisection.
Owing to Proposition~\ref{pr:discrete_norms_equiv}, this norm is equivalent to the
norm $\|\cdot\|_{p(\cdot),h}$ of Theorem~\ref{thm:error.estimate}, uniformly
in $h$, so that both decay at the same rate. { By Theorem~\ref{thm:error.estimate} and Remark~\ref{rem:error:rates-general},
provided the regularity seminorms in the estimate remain uniformly bounded,
the error decays at least with the order
\bae\label{eq:num.order}
  \mathcal{O}(k,p)=(k+1)\min\left(
  \frac1{\sobmax_\petp-1},\frac{\sobmin_\petm}{\sobmax_\petp}\right).
\eae
The dashed lines in the figures have slope
$\mathcal{O}(k,p)$; observed rates may be higher, and lower when the
regularity assumption fails. 
}
\br[Higher dimensions]\label{rem:num.dim}
All the tests of this section are carried out in $d=2$. The order
$\mathcal{O}(k,p)$ does not depend on $d$, and the numerical experiments
of~\cite{harnist-romain-conv}, run with the same scheme in $d=3$ and $d=4$
for a smooth solution and sinusoidal exponent fields, also illustrate the convergence behavior in higher dimensions.
\er

\subsection{Convergence rates for smooth solutions}
\label{sec:num.res.smooth}

\begin{figure}[htbp]
  \centering
  \includegraphics[width=0.48\textwidth]{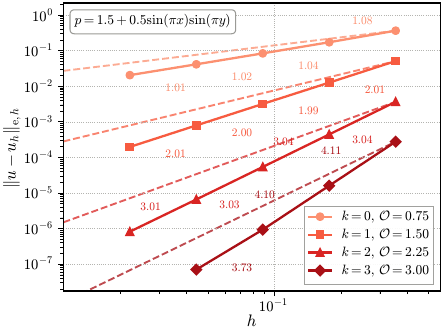}
  \hfill
  \includegraphics[width=0.48\textwidth]{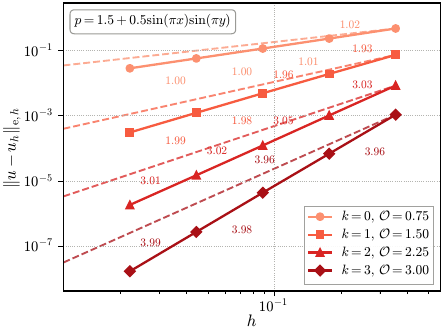}

  \vspace{0.3cm}
  \includegraphics[width=0.48\textwidth]{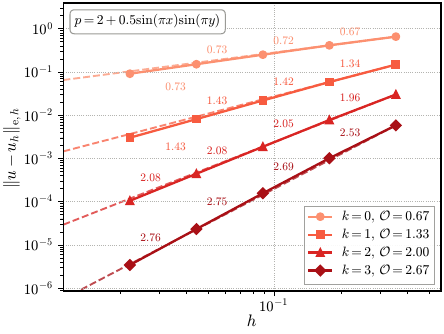}
  \hfill
  \includegraphics[width=0.48\textwidth]{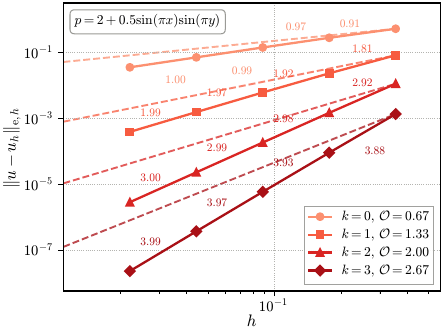}

  \vspace{0.3cm}
  \includegraphics[width=0.48\textwidth]{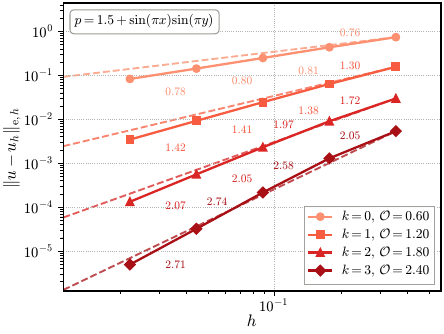}
  \hfill
  \includegraphics[width=0.48\textwidth]{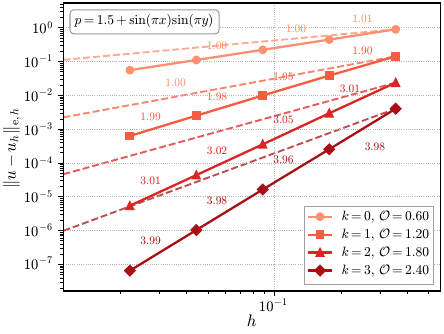}
  \caption{Error $\|u-u_h\|_{\mathrm e,h}$ versus the mesh
  size for the exponential solution~\eqref{eq:num.sol.exp}; dashed lines of
  slope $\mathcal{O}(k,p)$. \textbf{Left}: $\delta=0$. \textbf{Right}:
  $\delta=1$. \textbf{Top}: $p=1.5+0.5\sin(\pi x)\sin(\pi y)$.
  \textbf{Middle}: $p=2+0.5\sin(\pi x)\sin(\pi y)$. \textbf{Bottom}:
  $p=1.5+\sin(\pi x)\sin(\pi y)$.}
  \label{fig:num.exp}
\end{figure}

\begin{figure}[htbp]
  \centering
  \includegraphics[width=0.48\textwidth]{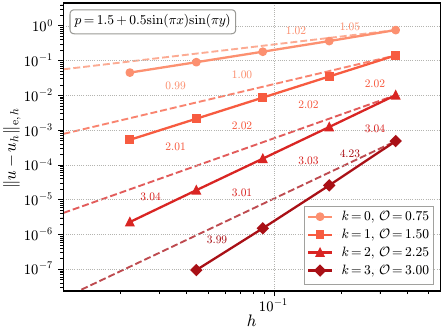}
  \hfill
  \includegraphics[width=0.48\textwidth]{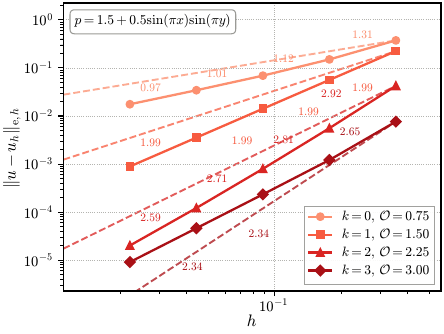}

  \vspace{0.3cm}
  \includegraphics[width=0.48\textwidth]{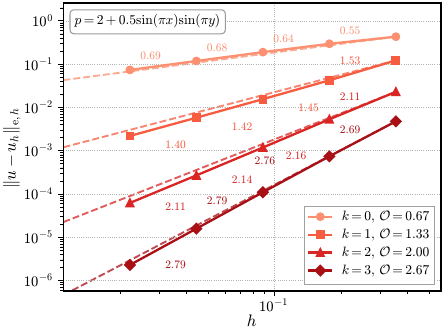}
  \hfill
  \includegraphics[width=0.48\textwidth]{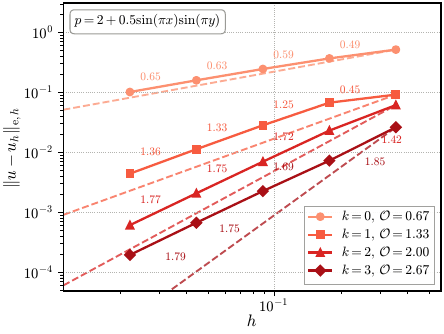}

  \vspace{0.3cm}
  \includegraphics[width=0.48\textwidth]{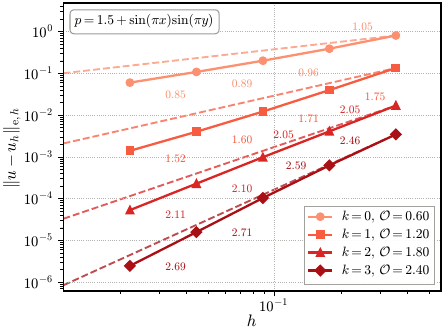}
  \hfill
  \includegraphics[width=0.48\textwidth]{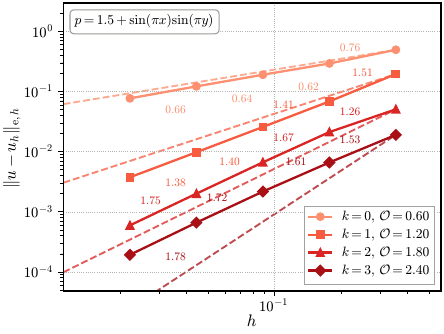}
  \caption{Error $\|u-u_h\|_{\mathrm e,h}$ versus the mesh
  size for $\delta=0$; dashed lines of slope $\mathcal{O}(k,p)$.
  \textbf{Left}: polynomial solution~\eqref{eq:num.sol.poly}.
  \textbf{Right}: sinusoidal solution~\eqref{eq:num.sol.sin}. \textbf{Top}:
  $p=1.5+0.5\sin(\pi x)\sin(\pi y)$. \textbf{Middle}:
  $p=2+0.5\sin(\pi x)\sin(\pi y)$. \textbf{Bottom}:
  $p=1.5+\sin(\pi x)\sin(\pi y)$.}
  \label{fig:num.poly.sin}
\end{figure}

We first consider the smooth solutions~\eqref{eq:num.sol.exp}
and~\eqref{eq:num.sol.poly}. Figure~\ref{fig:num.exp} shows how the error
decreases under mesh refinement for the exponential solution, on the three
exponent fields and for both values of $\delta$, and the left column of
Figure~\ref{fig:num.poly.sin} shows the same quantity for the polynomial
solution with $\delta=0$ (its right column is discussed in
Section~\ref{sec:num.res.sharp}).

For both solutions, the estimate of Theorem~\ref{thm:error.estimate}
holds in the three regimes. For $\delta=0$, the subquadratic regime
converges faster than predicted, at a rate close to the optimal order $k+1$.
{ For $\delta=1$, the experiments exhibit the optimal order $k+1$ in all three regimes. This is a numerical observation: the regime-dependent estimate of Theorem~\ref{thm:error.estimate:regime} only covers $p_\petp\le2$, and does not imply a general $k+1$ bound in the mixed or superquadratic regimes. Since the optimal order is also observed for $\delta=1$ with the other solutions of this section, whose curves are not shown, the other figures on
Cartesian meshes are restricted to $\delta=0$.}

\subsection{Role of the flux regularity assumption}
\label{sec:num.res.sharp}

We next examine whether the flux regularity assumption
$\stress(\cdot,\GRAD u)\in\tW^{k+1,\sob'(\cdot)}(\T_h)$ of Theorem~\ref{thm:error.estimate} is necessary. To this end, we
consider solutions for which this assumption fails beyond a certain degree.
For
$\delta=0$, this assumption can only fail at the critical points $\b x_0$ of
$u$, where $\GRAD u$ vanishes. Assume that $|\GRAD u(\b x)|$ behaves like
$|\b x-\b x_0|^{\gamma-1}$ near such a point, with $\gamma=2$ when the
Hessian of $u$ is invertible at $\b x_0$, and set $p_0=p(\b x_0)$; for the
fields~\eqref{eq:num.law}, $p-p_0=O(|\b x-\b x_0|^2)$ at all the critical
points considered below. If $p_0\neq2$, the flux behaves like
$|\b x-\b x_0|^{(\gamma-1)(p_0-1)}$ and its derivatives of order $m$ like
$|\b x-\b x_0|^{(\gamma-1)(p_0-1)-m}$, which are $p'$-integrable near
$\b x_0$ in dimension $d=2$ if and only if
\bae\label{eq:num.reg}
  m < (p_0-1)\Big(\gamma-1+\frac{2}{p_0}\Big).
\eae
If $p_0=2$, the leading part $\GRAD u$ of the flux is smooth, and the
regularity is limited by the term $(p-2)\ln|\GRAD u|\,\GRAD u$, of order
$|\b x-\b x_0|^{\gamma+1}\ln|\b x-\b x_0|$, which yields $m<\gamma+2$. These
conditions are local and do not depend on the position of $\b x_0$ with
respect to the mesh; they therefore also govern the broken space
$\tW^{k+1,\sob'(\cdot)}(\T_h)$. For $\delta=1$, the flux is smooth and none
of them applies.

\begin{figure}[htbp]
  \centering
  \includegraphics[width=0.48\textwidth]{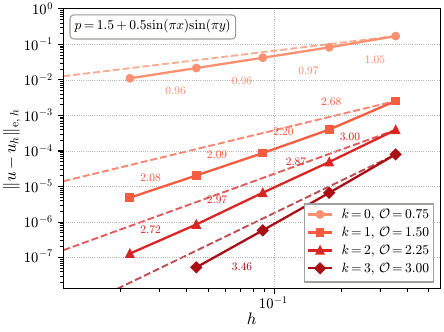}
  \hfill
  \includegraphics[width=0.48\textwidth]{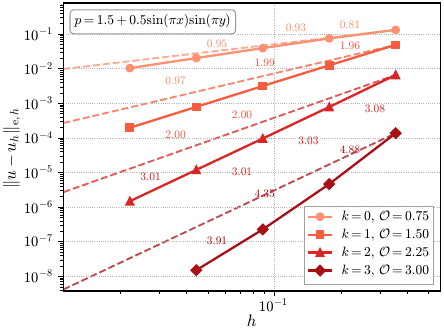}

  \vspace{0.3cm}
  \includegraphics[width=0.48\textwidth]{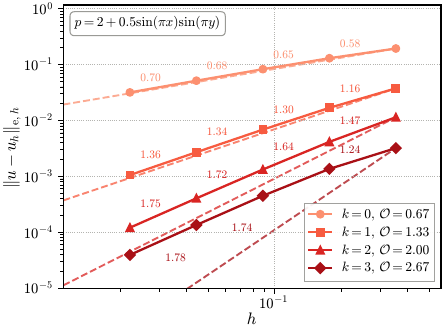}
  \hfill
  \includegraphics[width=0.48\textwidth]{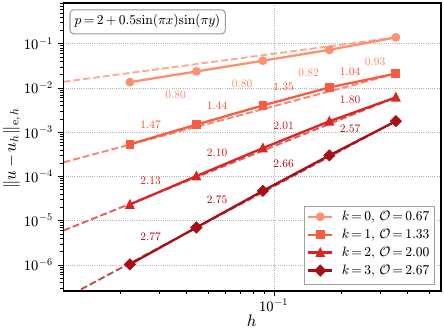}

  \vspace{0.3cm}
  \includegraphics[width=0.48\textwidth]{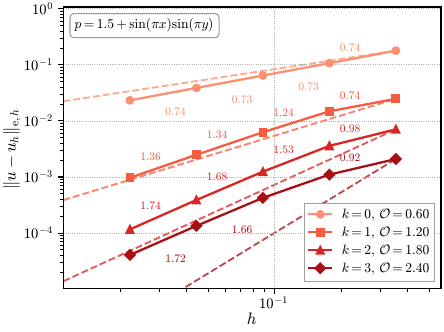}
  \hfill
  \includegraphics[width=0.48\textwidth]{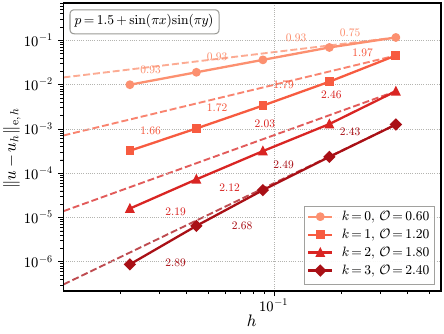}
  \caption{Error $\|u-u_h\|_{\mathrm e,h}$ versus the mesh
  size for the radial solutions~\eqref{eq:num.sol.rad}, $\delta=0$; dashed
  lines of slope $\mathcal{O}(k,p)$. \textbf{Left}: $\gamma=2$.
  \textbf{Right}: $\gamma=4$. \textbf{Top}: $p=1.5+0.5\sin(\pi x)\sin(\pi y)$.
  \textbf{Middle}: $p=2+0.5\sin(\pi x)\sin(\pi y)$. \textbf{Bottom}:
  $p=1.5+\sin(\pi x)\sin(\pi y)$.}
  \label{fig:num.rad}
\end{figure}

\paragraph{Radial solutions.} The radial
solutions~\eqref{eq:num.sol.rad} isolate this mechanism: they are
polynomials, and their gradient vanishes only at $\b x_c$, where the
exponent reaches its maximum $p_\petp$. For $\gamma=4$, the bounds
exceed $4$, and the assumption holds for every
$k\leq3$. For $\gamma=2$, the bound~\eqref{eq:num.reg} equals $2.7$ in the superquadratic and mixed
regimes, where $p_\petp=2.5$, so that the assumption fails from $k=2$ on;
in the subquadratic regime, $p_\petp=2$, and the assumption only fails for
$k=3$, through a logarithmic factor.
Figure~\ref{fig:num.rad} compares the two values of $\gamma$ for
$\delta=0$.

For $\gamma=4$, the estimate holds in the three regimes. For
$\gamma=2$, once the regularity assumption fails, the rate no longer
increases with $k$ in the superquadratic and mixed regimes, and the
estimate is lost for the highest degrees; in the subquadratic regime,
where the deficit is only logarithmic, the estimate still holds.

\paragraph{Sinusoidal solution.} Finally, we consider
the sinusoidal solution~\eqref{eq:num.sol.sin}, which has several critical
points, all with $\gamma=2$: the center $\b x_c$, where $p=p_\petp$, and the
four corners, where $p=p_\petm$. The right column of
Figure~\ref{fig:num.poly.sin} reports the error for $\delta=0$.

Because of the corners and of the center, the regularity assumption
fails from $k=1$ or $k=2$ on, depending on the regime. As for the radial
solution with $\gamma=2$, the rate eventually saturates in the three regimes, and
the estimate is lost for the highest degrees. 

These experiments highlight the role of flux regularity, without establishing optimality of the general rate $\mathcal O(k,p)$ (which is exceeded in the subquadratic tests). On the one hand,
the estimate holds in every configuration in which the regularity assumption
is satisfied. On
the other hand, the regularity assumption cannot be dispensed with: when it
fails at a critical point where $p\neq2$, raising the degree $k$ eventually
brings no gain, and the rate falls below $\mathcal{O}(k,p)$. The convergence
is then limited by the regularity of the flux at the critical points, and no
longer by the polynomial degree.

\subsection{Polygonal meshes}
\label{sec:num.res.voronoi}

Theorem~\ref{thm:error.estimate} holds on general polygonal meshes,
whereas all the meshes used so far are Cartesian. We now check that the
previous observations persist on Voronoi meshes, in two representative
cases: the exponential solution with $\delta=1$, and the sinusoidal solution
with $\delta=0$.

\begin{figure}[htbp]
  \centering
  \includegraphics[width=0.75\textwidth]{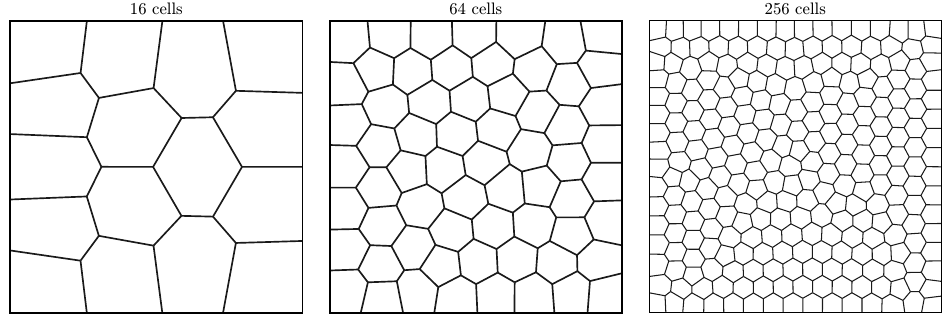}
  \caption{Voronoi meshes of $\Omega$. \textbf{Left}: $N=4$.
  \textbf{Middle}: $N=8$. \textbf{Right}: $N=16$.}
  \label{fig:num.vor.meshes}
\end{figure}

\paragraph{Meshes.} We use Voronoi tessellations of $\Omega$ with
$N^2$ cells, $N\in\{4,8,16,32,64\}$, obtained by Lloyd iterations from
random seeds, then post-processed so that no face is shorter than a
quarter of the largest diameter of its adjacent cells. The three coarsest meshes are shown in
Figure~\ref{fig:num.vor.meshes}.

\begin{figure}[htbp]
  \centering
  \includegraphics[width=0.48\textwidth]{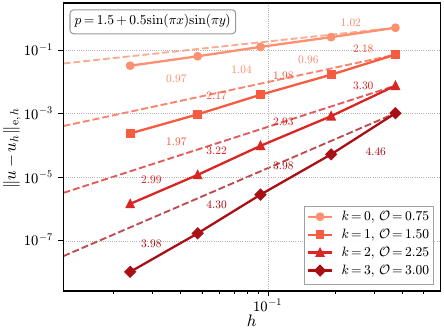}
  \hfill
  \includegraphics[width=0.48\textwidth]{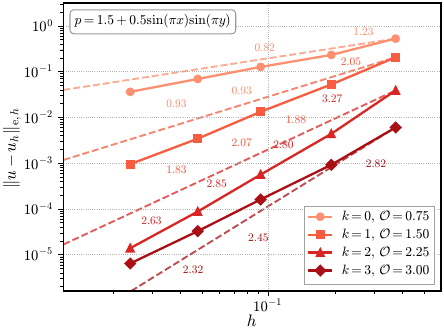}

  \vspace{0.3cm}
  \includegraphics[width=0.48\textwidth]{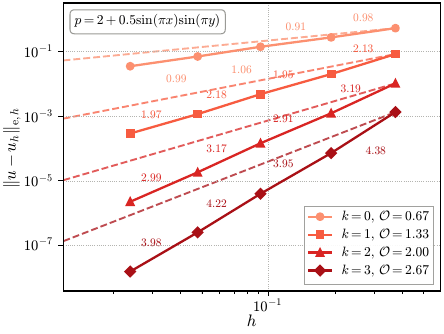}
  \hfill
  \includegraphics[width=0.48\textwidth]{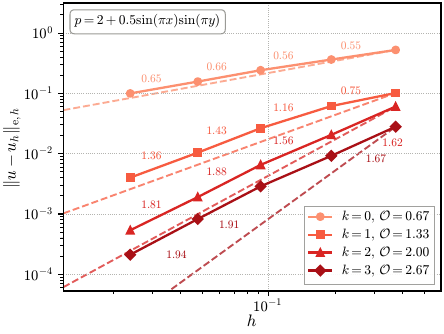}

  \vspace{0.3cm}
  \includegraphics[width=0.48\textwidth]{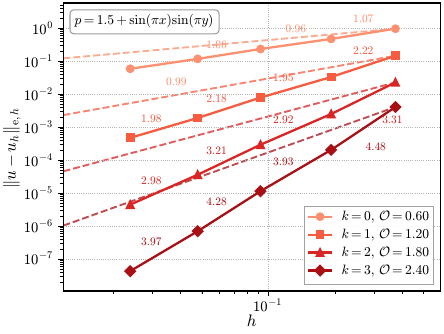}
  \hfill
  \includegraphics[width=0.48\textwidth]{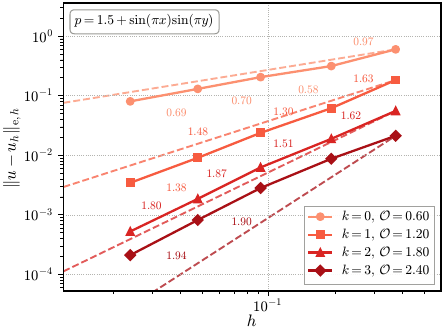}
  \caption{Error $\|u-u_h\|_{\mathrm e,h}$ versus the mesh
  size on the Voronoi meshes; dashed lines of slope $\mathcal{O}(k,p)$.
  \textbf{Left}: exponential solution~\eqref{eq:num.sol.exp}, $\delta=1$.
  \textbf{Right}: sinusoidal solution~\eqref{eq:num.sol.sin}, $\delta=0$.
  \textbf{Top}: $p=1.5+0.5\sin(\pi x)\sin(\pi y)$. \textbf{Middle}:
  $p=2+0.5\sin(\pi x)\sin(\pi y)$. \textbf{Bottom}:
  $p=1.5+\sin(\pi x)\sin(\pi y)$.}
  \label{fig:num.vor}
\end{figure}

\paragraph{Results.} Figure~\ref{fig:num.vor} reports
the error on the Voronoi meshes for the exponential solution with $\delta=1$
and for the sinusoidal solution with $\delta=0$. For the exponential
solution, as on Cartesian meshes, the convergence is optimal, of order $k+1$,
in the three regimes.

For the sinusoidal solution, since the bound~\eqref{eq:num.reg} does
not depend on the mesh, the regularity assumption fails in the same
configurations as on Cartesian meshes. The rates are close to those of
the right column of Figure~\ref{fig:num.poly.sin}, and the loss of rate
observed in Section~\ref{sec:num.res.sharp} is reproduced on the Voronoi
meshes.

Overall, the numerical results confirm the estimate of
Theorem~\ref{thm:error.estimate} and the role of flux regularity, independently of the
Cartesian structure of the meshes.

\section{Conclusion}\label{sec:conclusion}
We have established localized and global a priori estimates for the HHO discretization of variable-exponent Leray--Lions problems. The error bounds retain the pointwise growth structure and distinguish the purely subquadratic, mixed, and superquadratic cases. In the subquadratic case $p_\petp\le2$, the refined estimate recovers the optimal order $k+1$ when the approximation scale is sufficiently small relative to the local non-degeneracy scales. The analysis also accommodates independent degeneracy parameters in the bulk and stabilization and a class of anisotropic constitutive laws.

Numerical tests on Cartesian and Voronoi meshes corroborate the guaranteed rates and often exhibit faster convergence. For smooth non-degenerate tests with $\delta=1$, rates of $k+1$ are observed also in the mixed and superquadratic regimes. Conversely, examples with limited flux regularity show that increasing the approximation degree alone need not improve the asymptotic rate. Obtaining a sharp non-degenerate theory across the full variable-exponent range, and proving error bounds directly for nonhomogeneous Dirichlet data, remain natural extensions.

\section*{Acknowledgments}
A. Harnist acknowledges financial support from UTC’s Research Directorate through the Early-Career Researchers Program for the project \emph{``ANACONDAFUN – ANAlyses CONtinues et Discrétisations Avancées de FlUides Non-newtoniens''}.

\appendix
\section{Variable-exponent Lebesgue and Sobolev background}\label{app:ve-background}
{We collect here the variable-exponent background used in the error analysis. Some definitions and basic conventions are also stated in the companion convergence study~\cite{harnist-romain-conv}; they are repeated to make the present argument self-contained.}

\subsection{Lebesgue and Sobolev spaces with variable exponents}

Let $U\subset\RR^d$ be a bounded measurable set. On volume subsets of $\Om$, exponents are understood as restrictions of exponents defined on $\Om$. On mesh faces, including boundary faces, any exponent used in surface integrals is the restriction of the unique finite-valued continuous extension to $\overline\Om$ of an exponent in $\Plogp(\Om)$; intrinsic surface measure is used there. We write $\Plog(U)$ for the class of measurable exponents $p:U\to[1,\infty)$ satisfying the local log-H\"older condition: there exists $p_\cst{log}\in(0,\infty)$ such that, for all distinct $\b x,\b y\in U$,
\bae
 \left|\frac1{p(\b x)}-\frac1{p(\b y)}\right|
 \le \frac{p_\cst{log}}{\log(e+1/|\b x-\b y|)}.\label{eq:Plog:1}
\eae
\br[Unbounded domains]
When the underlying set is unbounded, local log-H\"older continuity is ordinarily combined with an additional condition controlling the exponent at infinity, involving a limit $p_\infty$ and the factor $\log(e+|\b x|)$; see \cite[Definition~4.1.1]{diening2011lebesgue}. This additional condition is not needed here, since $\Omega$ and all its subsets used below are bounded. No openness of $U$ is required for the definition of the exponent class itself; Sobolev and trace results are applied only on domains with the requisite geometry.
\er
{For $p\in\Plog(U)$ and a measurable $f : U \to \RR$, define the $L^{p(\cdot)}(U)$ modular and norm}
\bae
\rho_{L^{p(\cdot)}(U)}(f) &\coloneq \int_U |f(\b x)|^{p(\b x)}\,\cst{d}\b x, \\
\|f\|_{L^{p(\cdot)}(U)} &\coloneq \inf\{\lambda > 0 : \rho_{L^{p(\cdot)}(U)}(f/\lambda) \le 1 \}.
\eae
\br[Luxemburg norm]
{For constant $p(\cdot)$, the familiar relation is} $\|f\|_{L^{p}(U)} =  \rho_{L^{p}(U)}(f)^{1/p}$. {When the exponent varies, the appropriate quantity is the Luxemburg norm introduced above.} {More generally, each modular $\rho_\star(\cdot)$ is associated with the Luxemburg norm} 
\bae\label{eq:luxembourg}
\|{\cdot}\|_\star \eq \inf\{\lambda > 0 : \rho_\star({\cdot}/\lambda) \le 1 \}.
\eae 
\er
For a multi-index $\alpha=(\alpha_1,\ldots,\alpha_d)\in\mathbb N^d$, write $|\alpha|=\sum_{j=1}^d\alpha_j$ and $\partial_\alpha\eq\partial_{x_1}^{\alpha_1}\cdots\partial_{x_d}^{\alpha_d}$.
{For every $m \in \NN$ and measurable $f : U \to \RR$, the full and homogeneous $W^{m,p(\cdot)}(U)$ modulars are defined by}
\bae
\rho_{W^{m,p(\cdot)}(U)}(f) &\coloneq \sum_{|\alpha| \le m} \rho_{L^{p(\cdot)}(U)}(\partial_\alpha f) \quad \text{ and } \quad \rho_{\dot W^{m,p(\cdot)}(U)}(f) &\coloneq \sum_{|\alpha| = m} \rho_{L^{p(\cdot)}(U)}(\partial_\alpha f).
\eae
The Lebesgue-space definitions only require $U$ to be measurable; for the Sobolev-space definitions, $U$ is assumed to be open and weak derivatives are taken in the distributional sense.
{The space $L^{p(\cdot)}(U)$ comprises the measurable functions} $f : U \to \RR$ such that $\|f\|_{L^{p(\cdot)}(U)} < \infty$; {the Sobolev space} $W^{m,p(\cdot)}(U)$ {is defined analogously}.
{We shall also use the corresponding definitions for vector-valued functions} $\b f : U \to \RR^d$, applied componentwise; the resulting spaces are written in bold, e.g., $\tL^{p(\cdot)}(U)$.
{The following standard facts can be found, for example, in} \cite{diening2011lebesgue}.

{For scalar arguments, Young's generalized inequality states that, for all} $x,y \in [0,\infty)$, all $s\in[1,\infty)$ and $p,q\in(1,\infty)$ with $\frac{1}{s}=\frac{1}{p}+\frac{1}{q}$, and all $\varepsilon \in (0,\infty)$,
\bae\label{eq:gener:young:scal}
(xy)^s \le \varepsilon x^p + \varepsilon^{-q/p}y^q. 
\eae

\bpr
{Given $p(\cdot),q(\cdot),s(\cdot) \in \Plog(U)$ satisfying} $\frac{1}{s(\cdot)} = \frac{1}{p(\cdot)}+\frac{1}{q(\cdot)}$ in $U$, {and for any} $f \in L^{p(\cdot)}(U)$, $g \in L^{q(\cdot)}(U)$, {the following H\"older and Young estimates hold:}
\bse
\ba
\|fg\|_{L^{s(\cdot)}(U)} &\le 2\|f\|_{L^{p(\cdot)}(U)}\|g\|_{L^{q(\cdot)}(U)},\label{eq:holder}\\
\rho_{L^{s(\cdot)}(U)}(fg) &\le \rho_{L^{p(\cdot)}(U)}(f)+\rho_{L^{q(\cdot)}(U)}(g).\label{eq:young}
\ea
\ese
{A parameter-dependent form of Young's inequality is, for all} $\varepsilon \in (0,\infty)$,
\bae\label{eq:young:gener}
\rho_{L^{s(\cdot)}(U)}(fg) \le \varepsilon\rho_{L^{p(\cdot)}(U)}(f)+\rho_{L^{q(\cdot)}(U)}(\varepsilon^{-1/p(\cdot)}g).
\eae 
\epr

\subsection{Bounded variable-exponent class}

{For $p(\cdot) \in \Plog(U)$, we use the constants}
\bae
\sob_\petm = \essinf_U p \in [1,\infty) \quad \text{ and } \quad \sob_\petp = \esssup_U p \in [1,\infty].
\eae
The notation $\Plogm(U)$ (resp. $\Plogp(U)$) designates the class of $p(\cdot) \in \Plog(U)$ for which $\sob_\petm > 1$ (resp. $\sob_\petp < \infty$), and we define $\Plogpm(U) \eq \Plogm(U) \cap \Plogp(U)$.
For $p\in\Plogpm(U)$, the conjugate exponent and the auxiliary exponents are defined by
\bae\label{eq:exp:conj.lower.uper}
 p' \eq \frac{p}{p-1},\qquad \sobmax\coloneq\max(p,2),\qquad \sobmin\coloneq\min(p,2).
\eae
In this case $p',\sobmin,\sobmax\in\Plogpm(U)$.
We denote the conjugates of the endpoint values by $p'_\petm\coloneq(p_\petm)'$ and $p'_\petp\coloneq(p_\petp)'$ (and similarly, locally, by $p'_{T,-}\coloneq(p_{T,-})'$ and $p'_{T,+}\coloneq(p_{T,+})'$). Since taking the conjugate reverses the order, the essential extrema of $p'(\cdot)$ are instead $(p')_\petm=p'_\petp$ and $(p')_\petp=p'_\petm$, whenever $p_\petm>1$.
The following basic modular--norm relations hold for bounded exponents.
\bpr
For all $p(\cdot) \in \Plogp(U)$ and $f\in L^{p(\cdot)}(U)$, {where the first identity below is understood only for $f\ne0$}, we have 
\bse
\ba
&\hspace{-15mm}\rho_{L^{p(\cdot)}(U)}\left(\frac{f}{\|f\|_{L^{p(\cdot)}(U)}}\right) = 1,\label{eq:mod_over_norm_1}\\
\|f\|_{L^{p(\cdot)}(U)} \le 1 \quad &\Longleftrightarrow \quad \rho_{L^{p(\cdot)}(U)}(f) \le 1.\label{eq:equiv_norm_mod_1}
\ea
\ese
\epr
{The following statement is established in} \cite{diening2011lebesgue}.
\bpr 
{For any $p(\cdot) \in \Plogp(U)$ and $f \in L^{p(\cdot)}(U)$, the bounds read} 
\bse\label{eq:modular-norm-equiv}
\ba
\min(\|f\|_{L^{p(\cdot)}(U)}^{\sob_\petm},\|f\|_{L^{p(\cdot)}(U)}^{\sob_\petp}) &\le \rho_{L^{p(\cdot)}(U)}(f) \le \max(\|f\|_{L^{p(\cdot)}(U)}^{\sob_\petm},\|f\|_{L^{p(\cdot)}(U)}^{\sob_\petp}),\\
\min(\rho_{L^{p(\cdot)}(U)}(f)^{1/\sob_\petm},\rho_{L^{p(\cdot)}(U)}(f)^{1/\sob_\petp}) &\le \|f\|_{L^{p(\cdot)}(U)} \le \max(\rho_{L^{p(\cdot)}(U)}(f)^{1/\sob_\petm},\rho_{L^{p(\cdot)}(U)}(f)^{1/\sob_\petp}).
\ea
\ese
\epr
\subsection{{ Variable-exponent polynomial tools}}\label{sec:ve-tools}
{The polynomial estimates for variable exponents required in the discrete analysis are collected here.}
When a polynomial embedding is applied on a face, it is understood on the affine hull of that face, equipped with its intrinsic $(d-1)$-dimensional measure; the corresponding estimate follows from the same fixed-exponent polynomial argument in dimension $d-1$, under the induced uniform face-regularity assumptions. 

\subsubsection{{ Mean exponent equivalence and discrete Lebesgue embedding}}

The companion paper~\cite{harnist-romain-conv} extends to variable exponents the following polynomial Lebesgue embedding, cf. \cite[Lemma 5.1]{Di-Pietro.Droniou:17} and \cite[Remark 1.51]{Di-Pietro.Droniou:20}. Let $U\subset\RR^d$ be a measurable set with inradius $r_U \in (0,\infty)$ and diameter $h_U \in (0,\infty)$, such that
\bae\label{eq:mean.exp.regularity}
\frac{r_U}{h_U} \ge \kappa > 0.
\eae
Then, for every $p,q \in [1,\infty]$, all $k \in \NN$, and all $v \in \PP^k(U)$,
\bae\label{eq:lebesgue.emb:fixed}
\|v\|_{L^{p}(U)} \eqsim h_U^{\frac{d}{p}-\frac{d}{q}}\|v\|_{L^{q}(U)}
\eae
where the hidden constant, denoted $C_\cst{le}$, depends only on $k$, $d$, $p$, $q$, and $\kappa$.

\br[Dependency of the exponents]\label{rm:dependency_pq}
For finite exponents in fixed bounded intervals, the constant in~\eqref{eq:lebesgue.emb:fixed} may be chosen uniformly in these exponents. For suitably regular unions of simplices, the same constant can even be chosen independently of $p$ and $q$; see \cite[Remark~1.51]{Di-Pietro.Droniou:20}.
\er

\br[Alternative geometric assumptions]
\label{rm:assumption:lebesgue}
The results below extend to any geometric setting admitting a uniform fixed-exponent polynomial embedding of the form~\eqref{eq:lebesgue.emb:fixed}. In particular, the regular simplicial submeshes of Section~\ref{sec:discrete.setting} provide such an embedding for mesh cells, with constants governed by the polynomial-embedding constant; see \cite[Lemma~1.25]{Di-Pietro.Droniou:20}.
\er

The main tool is the following theorem.

\bt[Mean-exponent equivalence]\label{thm:mean_exp_equiv}
{Consider a bounded measurable $U\subset\Om\subset\RR^d$ satisfying} \eqref{eq:mean.exp.regularity} (or any assumption implying \eqref{eq:lebesgue.emb:fixed}, cf. Remark \ref{rm:assumption:lebesgue}).
Let $p \in \Plogp(U)$ and denote 
\bae\label{eq:mean.exp.pU}
p_U \eq \frac{1}{|U|}\int_U p(\b x)\,\cst{d}\b x.
\eae
Then, for all $\b x,\b y \in U$,
\bae\label{eq:mean.exp.diam.equiv}
h_U^{\frac{1}{p(\b x)}-\frac{1}{p_U}} \eqsim h_U^{\frac{1}{p(\b x)}-\frac{1}{p(\b y)}} \eqsim h_U^{p(\b x)-p_U} \eqsim h_U^{p(\b x)-p(\b y)} \eqsim 1. 
\eae
{In addition, for any $k \in \NN$ and $v \in \PP^k(U)$,}
\bae\label{eq:mean.exp.norm}
\|v\|_{L^{p(\cdot)}(U)} \eqsim \|v\|_{L^{p_U}(U)},
\eae
where all the hidden constants depend only on $\Om$, $k$, $d$, $p_\cst{log}$, $p_\petp$, and $C_\cst{le}$.
\et 

\br[{Dependence on $\Om$}]
{The dependence on $\Om$ is intrinsic because} the log-H\"older continuity controls the oscillation of the exponent but does not provide a uniform control of arbitrarily large diameters.
\er 

\noindent\emph{Proof.} See \cite{harnist-romain-conv}.

\bc[Discrete Lebesgue embedding]\label{cor:discrete-lebesgue-emb}
{Consider a bounded measurable $U\subset\Om\subset\RR^d$ satisfying} \eqref{eq:mean.exp.regularity} (or any assumption implying \eqref{eq:lebesgue.emb:fixed}, cf. Remark \ref{rm:assumption:lebesgue}).
Let $p(\cdot),q(\cdot) \in \Plogpm(U)$ and denote 
\bae
p_U \eq \frac{1}{|U|}\int_U p(\b x)\,\cst{d}\b x \quad \text{and} \quad q_U \eq \frac{1}{|U|}\int_U q(\b x)\,\cst{d}\b x.
\eae 
Then, for all $k \in \NN$, and all $v \in \PP^k(U)$,
\bae\label{eq:disc:leb:emb}
\|v\|_{L^{p(\cdot)}(U)} \eqsim h_U^{\frac{d}{p_U}-\frac{d}{q_U}}\|v\|_{L^{q(\cdot)}(U)},
\eae
where the hidden constant depends only on $\Om$, $k$, $d$, $p_\cst{log}$, $p_\petp$, $q_\cst{log}$, $q_\petp$, and $C_\cst{le}$.
\ec 

\noindent\emph{Proof.} See \cite{harnist-romain-conv}.

\section{Verification of the anisotropic Carreau--Yasuda example}\label{sec:properties.stress}

The following criterion is used to verify Assumption~\ref{ass:stress} for the anisotropic Carreau--Yasuda flux of Example~\ref{ex:Carreau--Yasuda;mat}.

\bl[Elementary two-point power inequality]\label{lem:two-point-power}
Let $p\in(1,\infty)$ and $A,B\ge0$ with $A+B>0$. Setting
$\zeta_p\coloneq 2^{(2-p)/p}$, one has
\bae\label{eq:two-point-power}
 \frac{A^{p-1}+B^{p-1}}{A+B}
 \begin{cases}
 \le \zeta_p(A^p+B^p)^{\frac{p-2}{p}},&1<p<2,\\[1mm]
 =1,&p=2,\\[1mm]
 \ge \zeta_p(A^p+B^p)^{\frac{p-2}{p}},&p>2.
 \end{cases}
\eae
\el
\bp
After dividing by $(A+B)^{p-2}$ and setting $s=A/(A+B)$, the assertion reduces to the fact that
\bae
 \frac{s^{p-1}+(1-s)^{p-1}}
 {(s^p+(1-s)^p)^{(p-2)/p}}
\eae
is symmetric about $s=1/2$ and has its maximum there when $1<p<2$ and its minimum there when $p>2$. Its value at $s=1/2$ is $2^{(2-p)/p}$.
\ep

\bt[{Jacobian bounds imply $(p,\delta)$-structure}]\label{thm:1d.power-framed}
Let $\stress : \RR^d \to \RR^d$ be a continuous function and assume that there exist $p \in (1,\infty)$, $\delta \in [0,\infty)$ and $\lambda_\cst{c},\lambda_\cst{m} \in (0,\infty)$ such that $\stress$ is differentiable on $\RR^d$ if $\delta>0$, and on $\RR^d\setminus\{\b 0\}$ if $\delta=0$, and
  \bse\label{eq:1d.power-framed:stress}
    \ba
  \lambda_{\min}(\GRAD\stress(\b\tau)) &\ge \lambda_\cst{m}(\delta^p+|\b\tau|^p)^\frac{p-2}{p},\\
      \lambda_{\max}(\GRAD\stress(\b\tau)) &\le \lambda_\cst{c}(\delta^p+|\b\tau|^p)^\frac{p-2}{p}.
    \ea
  \ese
where, for every $A \in \RR^{d\times d}$, we denote $\lambda_{\min}(A)$ as the smallest eigenvalue of $(A+A^\text{T})/2$ and $\lambda_{\max}(A)$ as the largest singular value of $A$ (equivalently, the square root of the largest eigenvalue of $A^\text{T}A$).
The Jacobian bounds are imposed wherever the derivative exists; in the case $\delta=0$ they are imposed for $\b\tau\neq\b 0$.
  Then, $\stress$ verifies the H\"older continuity and strong monotonicity
 \bse\label{eq:thm:stress}
 \ba
       |\stress(\b\tau)-\stress(\b\eta)| &\le \gamma_\cst{c}\lambda_\cst{c} {\left(\delta^p+|\b\tau|^p+|\b\eta|^p\right)^\frac{p-2}{p}}| \b\tau-\b\eta|,
    \label{eq:thm:stress:cont}
\\
      \left(\stress(\b\tau)-\stress(\b\eta)\right)\cdot\left(\b\tau-\b\eta\right) &\ge \gamma_\cst{m}\lambda_\cst{m}{\left(\delta^p+|\b\tau|^p+|\b\eta|^p\right)^\frac{p-2}{p}}|\b\tau-\b\eta|^{2},\label{eq:thm:stress:mono}
\ea
\ese
with
\bae
\gamma_\cst{c} \eq \frac{2^{(2-\sobmin)/p}}{\sobmin-1}
\quad \text{and} \quad
\gamma_\cst{m} \eq \frac{2^{(2-\sobmax)/p}}{\sobmax-1}.
\eae
\et

\bp
  (i) \emph{Strong monotonicity.}
First, for all $t \in [0,1]$, we set $\b\theta_t \eq \b\eta+t(\b\tau-\b\eta)$ and $\b\phi : [0,1] \to \RR^d$ by $\b\phi(t) = \stress(\b\theta_t)$. We obtain
  \bae
  \stress(\b\tau)-\stress(\b\eta) &= \b\phi(1)-\b\phi(0) = \int_0^1 \b\phi'(t)\,\cst{d}t = \int_0^1 \GRAD\stress(\b\theta_t)(\b\tau-\b\eta)\,\cst{d}t.
  \eae
  When $\delta=0$ and $1<p<2$, if the segment crosses the origin, the same identity follows by splitting the integral on the two sides of the crossing point and passing to the limit; the singularity $|t-t_0|^{p-2}$ is integrable since $p>1$.
  Therefore, if $p \le 2$, by convexity of $|\cdot|^p$, we have
  \bae
  \! (\stress(\b\tau)-\stress(\b\eta))\cdot (\b\tau-\b\eta) &\ge \left(\int_0^1 \lambda_{\min}(\GRAD\stress(\b\theta_t))\,\cst{d}t\right)|\b\tau-\b\eta|^2 \\
   &\ge \lambda_\cst{m}\left(\int_0^1 (\delta^p+|\b\theta_t|^p)^\frac{p-2}{p}\,\cst{d}t\right)|\b\tau-\b\eta|^2 \ge \lambda_\cst{m}(\delta^p+|\b\tau|^p+|\b\eta|^p)^\frac{p-2}{p}|\b\tau-\b\eta|^2.
  \eae

Assume now $p>2$. Set $A\coloneq(\delta^p+|\b\eta|^p)^{1/p}$ and $B\coloneq(\delta^p+|\b\tau|^p)^{1/p}$. On $\RR\times\RR^d$, endowed with the norm
$\|(s,\b z)\|_p=(|s|^p+|\b z|^p)^{1/p}$, the reverse triangle inequality yields
\bae
 (\delta^p+|\b\theta_t|^p)^{1/p}
 =\|(1-t)(\delta,\b\eta)+t(\delta,\b\tau)\|_p
 \ge |(1-t)A-tB|.
\eae
Hence, using Lemma~\ref{lem:two-point-power},
\bae
 \int_0^1(\delta^p+|\b\theta_t|^p)^{\frac{p-2}{p}}\,\cst dt
 &\ge \int_0^1 |(1-t)A-tB|^{p-2}\,\cst dt\\
 &=\frac{A^{p-1}+B^{p-1}}{(p-1)(A+B)}\ge \frac{2^{(2-p)/p}}{p-1}(A^p+B^p)^{\frac{p-2}{p}}\ge \frac{2^{(2-p)/p}}{p-1}
 (\delta^p+|\b\tau|^p+|\b\eta|^p)^{\frac{p-2}{p}}.
\eae
Inserting this estimate in the first line above proves \eqref{eq:thm:stress:mono} for $p>2$.

  \medskip
  (ii) \emph{H\"older continuity.}
If $p \ge 2$, by convexity of $|\cdot|^p$, we have
  \bae
   |\stress(\b\tau)-\stress(\b\eta)| &\le \left(\int_0^1 \lambda_{\max}(\GRAD\stress(\b\theta_t))\,\cst{d}t\right)|\b\tau-\b\eta| \\
   &\le \lambda_\cst{c}\left(\int_0^1 (\delta^p+|\b\theta_t|^p)^\frac{p-2}{p}\,\cst{d}t\right)|\b\tau-\b\eta| \le \lambda_\cst{c}(\delta^p+|\b\tau|^p+|\b\eta|^p)^\frac{p-2}{p}|\b\tau-\b\eta|.
  \eae
Assume now $1<p<2$ and keep $A$ and $B$ as above. Since $(p-2)/p<0$, the reverse triangle inequality gives
\bae
 \int_0^1(\delta^p+|\b\theta_t|^p)^{\frac{p-2}{p}}\,\cst dt
 &\le \int_0^1 |(1-t)A-tB|^{p-2}\,\cst dt=\frac{A^{p-1}+B^{p-1}}{(p-1)(A+B)}\\
 &\le \frac{2^{(2-p)/p}}{p-1}(A^p+B^p)^{\frac{p-2}{p}}\le \frac{2^{(2-p)/p}}{p-1}
 (\delta^p+|\b\tau|^p+|\b\eta|^p)^{\frac{p-2}{p}}.
\eae
Together with the Jacobian upper bound this proves \eqref{eq:thm:stress:cont}.
\ep

{
\bc[Non-symmetric anisotropic Carreau--Yasuda]\label{cor:Carreau--Yasuda:anisotropic}
Under the hypotheses of Example~\ref{ex:Carreau--Yasuda;mat}, the anisotropic flux \eqref{eq:Carreau--Yasuda:mat} satisfies Assumption~\ref{ass:stress}.
\ec

\bp
Fix $\b x$ and suppress the $\b x$-dependence. Set $\b y=\b A\b\tau$, $\b C=\mu\b I_d+\b K$, and $\phi(\b y)=(\delta^\alpha+|\b y|^\alpha)^{\frac{p-2}{\alpha}}$.
Since $\b B=\b A^\text{T}\b C\b A$, we have
\bae
\stress(\b\tau)=\b A^\text{T}\big(\phi(\b y)\b C\b y\big),
\qquad
D_{\b\tau}\stress(\b\tau)
=\b A^\text{T}D_{\b y}(\phi\b C\b y)\b A.
\eae
For $\b C=\mu\b I_d$, the symmetric matrix $D_{\b y}(\phi\mu\b y)$ has eigenvalues between
\bae
\mu\min(1,p-1)\phi
\qquad\text{and}\qquad
\mu\max(1,p-1)\phi.
\eae
The contribution of $\b K$ satisfies
\bae
\big\|D_{\b y}(\phi\b K\b y)\big\|_{\rm op}
\le \big(1+|p-2|\big)\|\b K\|_{\rm op}\phi.
\eae
Hence, if
\bae\label{eq:anisotropic:epsilon}
\varepsilon_0\Big(1+\max\{|p_\petm-2|,|p_\petp-2|\}\Big)
<\frac12\min\{1,p_\petm-1\},
\eae
then
\bae
\lambda_{\min}\!\left(\operatorname{sym}D_{\b y}(\phi\b C\b y)\right)
\ge \frac12\mu\min(1,p-1)\phi,
\eae
while $\|D_{\b y}(\phi\b C\b y)\|_{\rm op}\lesssim\mu\max(1,p-1)\phi$. The bounds on $\b A$ transfer these estimates to $D_{\b\tau}\stress$, and the equivalence
\bae
(\delta^\alpha+|\b A\b\tau|^\alpha)^{\frac{p-2}{\alpha}}
\eqsim
(\delta^p+|\b\tau|^p)^{\frac{p-2}{p}}
\eae
is uniform for $p\in[p_\petm,p_\petp]$, $\alpha\in[\alpha_\petm,\alpha_\petp]$, and the prescribed bounds on $\b A$. Theorem~\ref{thm:1d.power-framed} therefore yields Assumption~\ref{ass:stress}.
\ep
}

\begingroup\emergencystretch=2em\sloppy

\bibliographystyle{hho_doi_or_url}
\bibliography{hhollve}
\endgroup
\end{document}